\documentclass[12pt,final,reqno]{amsart}
\usepackage[square,sort,comma,numbers]{natbib}
\usepackage{amsmath}
\usepackage{amsthm}
\usepackage{amsfonts}  
\usepackage{amssymb}   
\usepackage{latexsym}
\usepackage{amscd}

\usepackage{hyperref} 

\DeclareMathAlphabet\mathoo{U}{eur}{b}{n}

 \DeclareMathOperator{\Ker}{Ker}
 \DeclareMathOperator{\Image}{Im}
 \DeclareMathOperator{\Real}{Re}
 \DeclareMathOperator{\Dom}{Dom}
 \DeclareMathOperator{\sgn}{sgn}
\newcommand\separ{\nobreak\hskip.1111em\mathpunct{}\nonscript\mkern-\thinmuskip{:}\hskip.3333emplus.0555em\relax}

\theoremstyle{plain}
\newtheorem{theorem}{Theorem}
\newtheorem{proposition}[theorem]{Proposition}
\newtheorem{corollary}[theorem]{Corollary}
\newtheorem{lemma}[theorem]{Lemma}

\theoremstyle{remark}
\newtheorem{remark}{Remark}
\newtheorem{example}{Example}

\renewcommand\:{\nobreak\hskip.1111em\mathpunct{}\nonscript\mkern-%
 \thinmuskip{:}\hskip.3333emplus.0555em\relax}

\begin{document}

\title[Analytic functional calculus for two operators]{Analytic functional calculus for two operators}

\author{V.G. Kurbatov}
 \address{Department of Mathematical Physics,
Voronezh State University\\
1, Universitetskaya Square, Voronezh 394036, Russia}
 \email{kv51@inbox.ru}

\author{I.V. Kurbatova}
\address{Department of Mathematics, Air Force Academy of the Ministry of Defense of the Russian Federation\\
54 ``A'', Starykh Bol'shevikov Str., Voronezh 394064, Russia}
 \curraddr{} \email{la\_soleil@bk.ru}

\author{M.N. Oreshina}
 \address{Department of Applied Mathematics,
Lipetsk State Technical University\\
30, Mos\-kov\-s\-kaya Str., Lipetsk 398600, Russia}
 \curraddr{} \email{M\_Oreshina@mail.ru}



\markright{Functional calculus for two operators}

\begin{abstract}
Properties of the mappings
\begin{align*}
C&\mapsto\frac1{(2\pi i)^2}\int_{\Gamma_1}\int_{\Gamma_2}f(\lambda,\mu)\,R_{1,\,\lambda}\,C\,
R_{2,\,\mu}\,d\mu\,d\lambda,\\
C&\mapsto\frac1{2\pi i}\int_{\Gamma}g(\lambda)R_{1,\,\lambda}\,C\,
R_{2,\,\lambda}\,d\lambda
\end{align*}
are discussed; here $R_{1,\,(\cdot)}$ and $R_{2,\,(\cdot)}$ are pseudo-resolvents, i.~e., resolvents of bounded, unbounded, or multivalued linear operators, and $f$ and $g$ are analytic functions. Several applications are considered: a representation of the impulse response of a second order linear differential equation with operator coefficients, a representation of the solution of the Sylvester equation, and an exploration of properties of the differential of the ordinary functional calculus.
\end{abstract}
\maketitle

\section{Introduction}
Let $A$ and $B$ be matrices of the sizes $n\times n$ and $m\times m$ respectively. The result $p(A,B)$ of the substitution of these matrices into a polynomial $p(\lambda,\mu)=\sum_{i,j=0}^Nc_{ij}\lambda^i\mu^j$ of \emph{two} variables is usually understood to be the mapping $C\mapsto\sum_{i,j=0}^Nc_{ij}A^iCB^j$ acting on $n\times m$-matrices, or to be the block matrix $\{a_{ij}B\}$. This article is devoted to extensions and applications of this construction.

First, the matrices $A$ and $B$ can be replaced by bounded linear operators acting in (infinite-dimensional) Banach spaces $X$ and $Y$ respectively. In this case, the number of interpretations of the object $p(A,B)$ increases. The most natural abstract interpretation is considering $p(A,B)$ as an operator acting in the completion of the algebraic tensor product~\cite{Defant-Floret,Schatten,Helemskii:eng} $X\otimes Y$ with respect to some cross-norm. This interpretation embraces many spaces of functions of two variables. For example~\cite{Defant-Floret,Schatten,Helemskii:eng},
$L_1[a,b]\overline{\otimes}_\pi L_1[c,d]$ is isometrically isomorphic to $L_1[a,b]\times[c,d]$, and $C[a,b]\overline{\otimes}_\varepsilon C[c,d]$ is isometrically isomorphic to $C[a,b]\times[c,d]$. But unfortunately, $L_\infty[a,b]\overline{\otimes}_\varepsilon L_\infty[c,d]$ is isomorphic only to a subspace of the space $L_\infty[a,b]\times[c,d]$. Another example, which can not be treated directly in terms of tensor products, is the interpretation of $p(A,B)$ as the transformation $C\mapsto\sum_{i,j=0}^Nc_{ij}A^iCB^j$ of the operators $C:\,Y\to X$. From the point of view of applications, the last example seems to be the most important. Therefore, in all cases, we call the operator that corresponds to $p(A,B)$ a \emph{transformator}; this term is conventional~\cite{Gohberg-Krein:eng} for the mappings of the type $C\mapsto\sum_{i,j=0}^Nc_{ij}A^iCB^j$ acting on operators $C$. In order to embrace the last example and some others, our treatment is based on the notion of an extended tensor product (Section~\ref{s:extended tensor product}) proposed in~\cite{Kurbatov-Kurbatova-IM-2015:eng}.

Second, one may replace the polynomial $p$ by an analytic function $f$. In this case, it is convenient to define $f(A,B)$ by means of a contour integral. For our main interpretation of $f(A,B)$ as a transformator acting on operators $C:\,Y\to X$, the relevant formula looks as follows:
\begin{equation}\label{e:*}
f(A,B)C=\frac1{(2\pi
i)^2}\int_{\Gamma_1}\int_{\Gamma_2}f(\lambda,\mu)(\lambda\mathbf1-A)^{-1}C
(\mu\mathbf1-B)^{-1}\,d\mu\,d\lambda.\tag{*}
\end{equation}
We call a correspondence of the type $f\mapsto f(A,B)$ that maps functions $f$ to transformators (operators) $f(A,B)$ a \emph{functional calculus}.

From the algebraic point of view, the functional calculus $f\mapsto f(A,B)$ possesses properties of the tensor product $\varphi_1\otimes\varphi_2$ of two ordinary functional calculi
\begin{equation*}
\varphi_1(f)=\frac1{2\pi i}\int_{\Gamma_1}f(\lambda)(\lambda
\mathbf1-A)^{-1}\,d\lambda,\qquad\varphi_2(f)=\frac1{2\pi i}\int_{\Gamma_2}f(\mu)(\mu
\mathbf1-B)^{-1}\,d\mu.
\end{equation*}
To emphasise this fact, we use the notation $\varphi_1\boxtimes\varphi_2$ and write $\bigl[(\varphi_1\boxtimes\varphi_2)f\bigr]C$ instead of $f(A,B)C$ when the transformator $f(A,B)$ acts in an extended tensor product. An important and nontrivial property of the transformation $\varphi_1\boxtimes\varphi_2$ is the spectral mapping theorem (Theorem~\ref{t:spectral map:R_1,R_2->C}).

Third, the basic properties of the functional calculus $\varphi_1\boxtimes\varphi_2$ are preserved when one replaces the resolvents  $R_\lambda=(\lambda\mathbf1-A)^{-1}$ of operators by \emph{pseudo-resolvents}~\cite{Hille-Phillips:eng}, i.~e. operator-valued functions $R_{(\cdot)}$ that satisfy the Hilbert identity
\begin{equation*}
R_\lambda-R_\mu=-(\lambda-\mu)R_\lambda R_\mu.
\end{equation*}
This generalization enables one to cover some additional examples. For example, a special case of a pseudo-resolvent is~\cite{Hille-Phillips:eng} the resolvent of an unbounded operator, and the most general example of a pseudo-resolvent is the resolvent of a linear relation or, in other terminology, a multivalued linear operator~\cite{Arens61,Baskakov2004:eng,Bichegkuev:eng,Cross,Favini-Yagi,Haase2006,Kurbatova-PMM09:eng}.

Moreover, in this article, we adhere to the point of view that a pseudo-resolvent is as a fundamental object as an operator (bounded, unbounded or multivalued) that generates it. The reason for this stems from the fact that when speaking of unbounded operators and linear relations we often actually work with their resolvents. For example, an unbounded operator is a generator of a strongly continuous or analytic semigroup if and only if~\cite{Hille-Phillips:eng,Yosida:eng,Dunford-Schwartz-I:eng} its resolvent satisfies a special estimate of the decay rate at infinity; in~\cite[Theorem 2.25]{Kato:eng} and~\cite[VIII.7]{Reed-Simon1:eng}, the natural convergence of unbounded operators is defined as the convergence of their resolvents in norm; and in~\cite{Reed-Simon:BAMS72,Reed-Simon:JFA73}, a function $f$ of unbounded operators $A$ and $B$ is defined as an (unbounded) operator $f(A,B)$ that possesses the following property: there exist sequences of bounded operators $A_n$ and $B_n$ such that the resolvents of $A_n$, $B_n$, and $f(A_n,B_n)$ converge in norm to the resolvents of $A$, $B$, and $f(A,B)$ respectively. Another argument (not used in this article) is that there is no analogue of unbounded and multivalued operators in Banach algebras, but, nevertheless, there are evident analogues of the resolvents of such operators.

This approach enables one to extend the notion of $f(A,B)$ to meromorphic functions $f$ (Theorem~\ref{t:mero f:2}): we define the result of the action of a meromorphic function on $A$ and $B$ to be a new pseudo-resolvent
and do not discuss which operator it is generated by.
Along the way, we answer (Corollary~\ref{c:mero f:2}) the question of the independence of the definition of $f(A,B)$ for unbounded $A$ and $B$ posed in~\cite{Reed-Simon:BAMS72,Reed-Simon:JFA73} (see the previous paragraph) from the choice of approximating sequences $A_n$ and $B_n$.

Many important applications are connected with the special cases of construction~\eqref{e:*} and their modifications. For example, it often occurs that the function $f$ depends on the difference or the sum of its arguments: the transformator $C\mapsto AC-CB$ generated by the function $f(\lambda,\mu)=\lambda-\mu$ is related to the Sylvester equation (Section~\ref{s:Sylvester}), and the transformator $C\mapsto e^{At}Ce^{Bt}$ generated by the function $f(\lambda,\mu)=e^{t(\lambda+\mu)}$ is connected with the stability theory of differential equations~\cite{Antoulas,Daletskii-Krein:eng}.

The version (Section~\ref{s:boxdot})
\begin{equation}\label{e:**}
f^{[1]}(A,B)C=\frac1{2\pi i}\int_{\Gamma}f(\lambda)(\lambda\mathbf1-A)^{-1}C
(\lambda\mathbf1-B)^{-1}\,d\lambda\tag{**}
\end{equation}
of functional calculus~\eqref{e:*} frequently occurs in applications; it involves functions $f$ of one variable. For example, expression~\eqref{e:**} with $f(\lambda)=e^{\lambda t}$ forms the principal part of the representation of a solution of the second order differential equation $\bigl(\frac d{dt}-A\bigr)C\bigl(\frac d{dt}-B\bigr)y=0$ (Section~\ref{s:pulse char}). Further (Section~\ref{s:df}), the differential of the ordinary functional calculus $A\mapsto f(A)$ at a point $A$ can be represented in the form (Theorem~\ref{t:differential})
\begin{equation*}
df(\cdot,A)=f^{[1]}(A,A).
\end{equation*}
It turns out that (Theorems~\ref{t:Theta}) $f^{[1]}(A,B)C$ coincides with~\eqref{e:*} provided $f^{[1]}$ is understood to be the divided difference
\begin{equation*}
f^{[1]}(\lambda,\mu)= \begin{cases}
\frac{f(\lambda)-f(\mu)}{\lambda-\mu}, & \text{if $\lambda\neq\mu$},\\
f'(\lambda), & \text{if $\lambda=\mu$},\\
0,& \text{if $\lambda=\infty$ or $\mu=\infty$}.
 \end{cases}
\end{equation*}

When choosing the level of generality of our exposition, we proceed from the following principles. First, in order
that a specialist in the classical operator theory can use the article we are trying to minimize explicit mention of Banach algebras and linear relations (multivalued operators) at least in main statements.
At the same time, when using the operator language, we aspire the maximal generality and, in particular, where possible, consider the case of an arbitrary pseudo-resolvent (and thereby, implicitly, the cases of unbounded operators and linear relations). Second, we are trying not to fall outside the framework of the theory of analytic functions of an operator and thus, for example, do not discuss issues related to the generators of semigroups. Third, as far as possible, we avoid the implicit use of operator pencils $\lambda\mapsto\lambda F-G$ with $F\neq\mathbf1$ because this approach leads to very cumbersome formulae. Finally, we restrict ourselves to the consideration of functions of two variables, suggesting that the generalization to three and more variables will not cause significant difficulties.

The literature on the subject under discussion is extremely extensive. Therefore, the bibliography can not be made comprehensive; the presented references reflect authors' tastes and interests. Many additional references can be found in the cited articles and books.

Sections~\ref{s:Banach algebras}--\ref{s:extended tensor product} outlines preliminary information. Here we recall and fix notation and the main facts in a convenient form. In Section~\ref{s:Banach algebras}, the terminology connected with Banach algebras and their properties is recalled. In Section~\ref{s:mathoo O}, the basic properties of algebras of analytic functions of one and two variables are described. In Section~\ref{s:pseudo-resolvents}, we discuss the notion of pseudo-resolvent and recall the construction of the functional calculus of analytic functions of one variable (Theorems~\ref{t:functional_calculas:infty in rho} and~\ref{t:functional_calculas:infty in sigma}) including the spectral mapping theorem (Theorem~\ref{t:spectral mappping}). In Section~\ref{s:extended tensor product}, the definition of the extended tensor product is given, the main examples are described, and the construction of the functional calculus of operator-valued analytic functions of one variable (Theorems~\ref{t:functional_calculas:infty in rho:R:otimes} and~\ref{t:functional_calculas:infty in sigma:R:otimes}) is recalled as well as the relevant spectral mapping theorem (Theorem~\ref{t:spectral map:Res}).

In Section~\ref{s:boxtimes}, we present the construction of the functional calculus~\eqref{e:*} of functions of two variables (Theorems~\ref{t:functional_calculas:infty in rho:R_1,R_2->C}, \ref{t:functional_calculas:infty in rho+sigma:R_1,R_2->C}, and~\ref{t:functional_calculas:infty in sigma:R_1,R_2->C}) and prove the corresponding spectral mapping theorem (Theorem~\ref{t:spectral map:R_1,R_2->C}). In Section~\ref{s:mero}, we extend these results to meromorphic functions. A well-known example of a meromorphic function of an operator is a polynomial of an unbounded operator (a polynomial has a pole at infinity, and the point at infinity belongs to the extended spectrum of an unbounded operator). This example shows that the result of applying a meromorphic function cannot be a bounded transformator; as a convenient tool for its description we use not the resulting object itself, but its resolvent, and we interpret the extended singular set of this resolvent as its spectrum (Theorem~\ref{t:mero f:2}).

In Section~\ref{s:boxdot}, we discuss modified variant~\eqref{e:**} of the functional calculus of functions of two variables. The connection between functional calculi~\eqref{e:*} and~\eqref{e:**} is established as well as some properties of functional calculus~\eqref{e:**}. The subsequent sections are devoted to applications. In Section~\ref{s:pulse char}, the pencil $\lambda\mapsto\lambda^2E+\lambda F+H$ of the second order is considered; it is induced by the equation $E\ddot y(t)+F\dot y(t)+Hy(t)=0$; we assume that the pencil admits a factorization, i.~e. the representation in the form of a product of two linear pencils. In such a case, the solution of the differential equation is expressed by a transformation of the kind~\eqref{e:**} (Theorem~\ref{t:pulse char:2}). In Section~\ref{s:Sylvester}, we discuss the properties of the transformator $Q:\,C\mapsto Z$ generated by the Sylvester equation $AZ-ZB=C$ (Theorem~\ref{t:Sylvester}). Finally, in Section~\ref{s:df}, it is shown that the differential of the ordinary functional calculus $A\mapsto f(A)$ is also a kind of transformator~\eqref{e:**} (Theorem~\ref{t:differential}).

\section{Banach algebras}\label{s:Banach algebras}
In this Section, we clarify the terminology connected with Banach algebras~\cite{Hille-Phillips:eng,Rudin:eng,Bourbaki_Theories_Spectrales:eng} and recall some of their properties.

In this article, all linear spaces and algebras are assumed to be complex.

Let $X$ and $Y$ be Banach spaces. We denote by $\mathoo B(X,Y)$ the set of all bounded linear operators
 $A\separ X\to Y$. When $X=Y$, we use the shorthand symbol $\mathoo B(X)$. The symbol $\mathbf1=\mathbf1_X$
stands for the identity operator. We adhere to the following notations: $X^*$ denotes the conjugate space of $X$; $\langle x,x^*\rangle$~denotes the value of the functional $x^*\in
X^*$ on $x\in X$, and $\langle x^{**},x^*\rangle$~is the value of
$x^{**}\in X^{**}$ on $x^*\in X^*$; $A^*$ denotes the conjugate operator of
$A\in\mathoo B(X,Y)$. The \emph{preconjugate} of an operator $A\in\mathoo B(Y^*,X^*)$
is an operator $A^0\in\mathoo B(X,Y)$ such that $(A^0)^*=A$.

The \emph{unit}~\cite{Hille-Phillips:eng,Rudin:eng,Bourbaki_Theories_Spectrales:eng} of an algebra $\mathoo B$ is an element $\mathbf1\in\mathoo B$ such that $\mathbf{1}A=A\mathbf1$ for all $A\in\mathoo B$. If an algebra $\mathoo B$ has a unit, it is called \emph{an algebra with a unit} or \emph{unital}.

A subset $\mathoo R$ of an algebra $\mathoo B$ is called a \emph{subalgebra} if $\mathoo R$ is stable under the algebraic operations (addition, scalar multiplication, and multiplication), i.~e. $A+B,\lambda A,AB\in\mathoo R$ for all $A,B\in\mathoo R$ and $\lambda\in\mathbb C$.
If the unit $\mathbf1$ of an algebra $\mathoo B$ belongs to its subalgebra $\mathoo R$, then $\mathoo R$ is called a \emph{subalgebra with a unit} or a \emph{unital subalgebra}.

Let $\mathoo B$ be a non-unital algebra. The set $\widetilde{\mathoo B}=\mathbb
C\oplus\mathoo B$ with the componentwise linear operations and the multiplication
$(\alpha,A)(\beta,B)=(\alpha\beta,\alpha B+\beta A+AB)$ is obviously an algebra with the unit
$(1,\mathbf0)$, where $\mathbf0$ is the zero of the algebra $\mathoo B$. The algebra $\widetilde{\mathoo B}$ is called the algebra derived from $\mathoo B$ by adjoining a unit or an algebra with an \emph{adjoint unit}.
The symbol $\mathbf1$ stands for the element $(1,\mathbf0)$, and the symbol $\alpha\mathbf1+A$ denotes the element $(\alpha,A)$. If $\mathoo B$ is a normed algebra, we set $\Vert(\alpha,A)\Vert=|\alpha|+\Vert A\Vert$.
Clearly, this formula defines a norm on $\widetilde{\mathoo B}$. It is also clear that
$\widetilde{\mathoo B}$ is complete provided that $\mathoo B$ is complete. If $\mathoo B$
is unital, then $\widetilde{\mathoo B}$ means the algebra $\mathoo B$ itself.

\begin{theorem}\label{t:Neumann} Let $\mathoo B$ be a unital Banach algebra and
$A,B\in\mathoo B$. If $A$ is invertible and
\begin{equation*}
\Vert B\Vert\cdot\Vert A^{-1}\Vert<1,
\end{equation*}
then the element $A-B$ is also invertible and
\begin{equation*}
(A-B)^{-1}=A^{-1}+A^{-1}BA^{-1}+A^{-1}BA^{-1}BA^{-1}
+A^{-1}BA^{-1}BA^{-1}BA^{-1}+\ldots\,.
\end{equation*}
In this case
\begin{align}
\Vert(A-B)^{-1}\Vert&\le\frac{\Vert A^{-1}\Vert}{1-\Vert B\Vert\cdot\Vert A^{-1}\Vert},\notag\\
\Vert(A-B)^{-1}-A^{-1}\Vert&\le\frac{\Vert B\Vert\cdot\Vert A^{-1}\Vert^{2}}{1-\Vert
B\Vert\cdot\Vert A^{-1}\Vert},\label{e:(A-B)-1-A-1}\\
\Vert(A-B)^{-1}-A^{-1}-A^{-1}BA^{-1}\Vert&\le\frac{\Vert B\Vert^2\cdot\Vert A^{-1}\Vert^{3}}{1-\Vert
B\Vert\cdot\Vert A^{-1}\Vert}.\notag
\end{align}
\end{theorem}
\begin{proof} We consider the series
\begin{equation*}
A^{-1}+A^{-1}BA^{-1}+A^{-1}BA^{-1}BA^{-1} +A^{-1}BA^{-1}BA^{-1}BA^{-1}+\ldots\,.
\end{equation*}
We represent this series in the form
$A^{-1}\bigl(\mathbf1+BA^{-1}+BA^{-1}BA^{-1} +BA^{-1}BA^{-1}BA^{-1}+\ldots\bigr).$
Since $\Vert BA^{-1}\Vert\le\Vert B\Vert\cdot\Vert A^{-1}\Vert<1$, the series converges absolutely. We denote its sum by~$C$. It is straightforward to verify that $C$ coincides with the inverse of $A-B$.

Estimates~\eqref{e:(A-B)-1-A-1} follow from the geometric series formula $\sum_{k=0}^\infty q^k=\frac1{1-q}$, $|q|<1$. For example, let us prove the second estimate:
\begin{equation*}
\Vert(A-B)^{-1}-A^{-1}\Vert\le\sum_{k=1}^\infty\Vert B\Vert^k\cdot\Vert
A^{-1}\Vert^{k+1}=\frac{\Vert B\Vert\cdot\Vert A^{-1}\Vert^{2}}{1-\Vert B\Vert\cdot\Vert
A^{-1}\Vert}.\qed
\end{equation*}
\renewcommand\qed{}
\end{proof}

Let $\mathoo B$ be a (nonzero) unital algebra and $A\in\mathoo B$. The set of all $\lambda\in\mathbb C$
such that the element $\lambda\mathbf1-A$ is not invertible is called the {\it spectrum}
of the element $A$ (in the algebra $\mathoo B$) and is denoted by the symbol $\sigma(A)$. The complement
$\rho(A)=\mathbb C\setminus\sigma(A)$ is called the \emph{resolvent set\/} of $A$.
The function
$
R_\lambda=(\lambda\mathbf1-A)^{-1}
$
is called the \emph{resolvent\/} of the element~$A$.

 \begin{proposition}[{\rm\cite[Theorem 4.8.1]{Hille-Phillips:eng}}]\label{p:resolvent proprties}
The resolvent $R_{(\cdot)}$ of any element $A\in\mathoo B$ satisfies the {\it Hilbert identity}
 \begin{equation}\label{e:Hilbert identity:alg}
R_\lambda-R_\mu=-(\lambda-\mu)R_\lambda R_\mu,\qquad\lambda,\mu\in\rho(A).
 \end{equation}
 \end{proposition}

\begin{proposition}[{\rm\cite[ch.1, \S~2, п.~5]{Bourbaki_Theories_Spectrales:eng},~\cite[Theorem 10.13]{Rudin:eng}}\bf]\label{p:spectra properties}
The spectrum of any element $A$ of a nonzero unital Banach algebra $\mathoo B$ is a compact and nonempty subset of the complex plane $\mathbb C$.
\end{proposition}

Let $\mathoo A$ and $\mathoo B$ be algebras. A mapping $\varphi\separ\mathoo
A\to\mathoo B$ is called~\cite{Bourbaki_Theories_Spectrales:eng} a {\it morphism of algebras} if
\begin{align*}
\varphi(A+B)&=\varphi(A)+\varphi(B),\\
\varphi(\alpha A)&=\alpha\varphi(A),\\
\varphi(AB)&=\varphi(A)\varphi(B).
\end{align*}
If, in addition, $\mathoo A$ and $\mathoo B$ are unital and
\begin{equation*}
\varphi(\mathbf 1_{\mathoo A})=\mathbf 1_{\mathoo B},
\end{equation*}
$\varphi$ is called a {\it morphism of unital algebras}.
If $\mathoo A$ and $\mathoo B$ are Banach algebras~\cite{Hille-Phillips:eng,Rudin:eng,Bourbaki_Theories_Spectrales:eng} and, in addition, the morphism $\varphi$ is continuous, then $\varphi$ is called a {\it morphism of Banach algebras}.

A unital subalgebra $\mathoo R$ of a unital algebra $\mathoo B$ is called~\cite[ch.~1,~\S~1.4]{Bourbaki_Theories_Spectrales:eng} {\it full} if every $B\in\mathoo R$ which is invertible in
$\mathoo B$ is also invertible in $\mathoo R$. Since the inverse is unique, the last definition is equivalent to the following: if for $B\in\mathoo R$ there exists
$B^{-1}\in\mathoo B$ such that $BB^{-1}=B^{-1}B=\mathbf 1$, then $B^{-1}\in\mathoo R$.

 \begin{example}\label{ex:preconjugate}
Let $X$ be a Banach space. The set $\mathoo B_0(X^*)$ of all operators that have a preconjugate is a full subalgebra of the algebra $\mathoo B(X^*)$.
 \end{example}

 \begin{proposition}[{\rm\cite[ch. 1, \S~2.5]{Bourbaki_Theories_Spectrales:eng}}]\label{p:closure of full}
The closure of a full subalgebra of a Banach algebra is also a full subalgebra.
The closure of the least full subalgebra of a Banach algebra that contains a set $M$ is the least full closed subalgebra that contains $M$.
 \end{proposition}

An algebra $\mathoo B$ is called \emph{commutative} if $AB=BA$ for all $A,B\in\mathoo B$.

A \emph{character} of a unital commutative algebra $\mathoo B$~\cite[ch.~1, \S~1.5]{Bourbaki_Theories_Spectrales:eng} is a morphism $\chi\:\mathoo B\to\mathbb C$ of unital algebras.
A \emph{character} of a commutative non-unital algebra $\mathoo B$~\cite[ch.~1, \S~1.5]{Bourbaki_Theories_Spectrales:eng} is a morphism of (non-unital) algebras $\chi\:\mathoo B\to\mathbb C$. If an algebra
$\mathoo B$ is non-unital, we denote by $\chi_0$ the character
$\chi_0\:\mathoo B\to\mathbb C$ that is equal to zero on all elements of
$\mathoo B$. We call $\chi_0$ the \emph{zero character}.
We stress that the character $\chi_0$ exists only if the algebra $\mathoo B$
is non-unital.

 \begin{proposition}\label{p:char and adj unit}
All characters of a non-unital algebra $\mathoo B$ are extendable uniquely to characters of the algebra
$\widetilde{\mathoo B}$ derived from $\mathoo B$ by adjoining a unit. The extension is defined by the formula
$\chi(\alpha\mathbf1+A)=\alpha+\chi(A)$. Conversely, the restriction of any character of the algebra
$\widetilde{\mathoo B}$ to $\mathoo B$ is a character of the algebra $\mathoo B$. In particular, the zero character
$\chi_0$ is the restriction of the character $\alpha\mathbf1+A\mapsto\alpha${\rm;} we will denote it by the same symbol $\chi_0$.
\end{proposition}
 \begin{proof}
The proof is obvious.
 \end{proof}

We denote by $X(\mathoo B)$ the set of all \emph{nonzero} characters of a commutative algebra $\mathoo B$ (unital or non-unital), and we denote by $\widetilde{X}(\mathoo B)$ the set of all characters of a commutative algebra $\mathoo
B$ (including the zero character $\chi_0$ if the algebra is non-unital). If an algebra $\mathoo B$ is unital, then $\widetilde{X}(\mathoo B)$ obviously coincides with $X(\mathoo B)$. The set $X(\mathoo
B)$ is called~\cite{Bourbaki_Theories_Spectrales:eng} the \emph{character space} of the algebra $\mathoo B$.

 \begin{theorem}[{\rm\cite[ch.~1, \S~3.3, Proposition 3]{Bourbaki_Theories_Spectrales:eng}}]\label{t:Gelfand}
Let $\mathoo B$ be a unital commutative Banach algebra. Then for all $A\in\mathoo B$,
\begin{equation*}
\sigma(A)=\{\,\chi(A)\:\chi\in \widetilde{X}(\mathoo B)\,\}.
\end{equation*}
 \end{theorem}

 \begin{corollary}[{\rm\cite[ch.~1, \S~3, Theorem 1]{Bourbaki_Theories_Spectrales:eng}}]\label{c:contin of characters}
Every character of a commutative Banach algebra is continuous{\rm;} namely, its norm is less than or equal to unity.
 \end{corollary}

 \begin{corollary}\label{c:continuity of spr}
In a unital commutative Banach algebra $\mathoo B$, the spectrum continuously depends on an element{\rm;} more precisely, if $A,B\in\mathoo B$ and $\Vert A-B\Vert<\varepsilon$, then $\sigma(B)$ is contained in the $\varepsilon$-neighbourhood of $\sigma(A)$.
 \end{corollary}

\section{Algebras of analytic functions}\label{s:mathoo O}
This Section is a preparation for a discussion of analytic functional calculi. Here we collect some preliminaries on algebras of analytic functions defined on subsets of $\overline{\mathbb C}$ and~$\overline{\mathbb C}^2$.

We denote by $\overline{\mathbb C}$ the one-point compactification $\mathbb C\cup\{\infty\}$ of the complex plane
$\mathbb C$, and we denote by $\overline{\mathbb C}^2$ the Cartesian product $\overline{\mathbb C}\times\overline{\mathbb C}$.

 \begin{proposition}\label{p:U times V}
Let $\sigma_1,\sigma_2\subseteq\overline{\mathbb C}$ be closed sets, and let an open set
$W\subseteq\overline{\mathbb C}^2$ contains $\sigma_1\times\sigma_2$.
Then there exist open sets $U,V\subseteq\overline{\mathbb C}$ such that
$\sigma_1\times\sigma_2\subseteq U\times V\subseteq W$.
 \end{proposition}
 \begin{proof}
For an arbitrary $\lambda\in\sigma_1$, we consider the set $\{\lambda\}\times\sigma_2$. We consider a finite cover of $\{\lambda\}\times\sigma_2$ by the sets of the form
$U_i\times V_i$, where $U_i,V_i\subseteq\overline{\mathbb C}$ are open and $U_i\times
V_i\subseteq W$. We put $\widetilde{U}=\cap_iU_i$ and $\widetilde{V}=\cup_iV_i$. It is clear that the set
$\widetilde{U}\times\widetilde{V}\subseteq W$ also covers the set
$\{\lambda\}\times\sigma_2$, namely,
$\{\lambda\}\subseteq\widetilde{U}$, $\sigma_2\subseteq\widetilde{V}$.

Further, we cover every subset of the form $\{\lambda\}\times\sigma_2$ of the set
$\sigma_1\times\sigma_2$ by a set of the form
$\widetilde{U}\times\widetilde{V}\subseteq W$. We choose a finite subcover $\{\,\widetilde{U}_k\times\widetilde{V}_k\,\}$ and put
$U=\cup_k\widetilde{U}_k$, $V=\cap_k\widetilde{V}_k$. Obviously,
$\sigma_1\times\sigma_2\subseteq U\times V\subseteq W$.
 \end{proof}

 \begin{proposition}\label{p:U times V 2}
Let $U_1,U_2\subseteq\overline{\mathbb C}$ be open sets. Then for every
compact set $N\subset U_1\times U_2$ there exist compact sets $N_1\subseteq U_1$ and $N_2\subseteq U_2$ such that $N\subseteq N_1\times N_2$.
 \end{proposition}
 \begin{proof}
It is sufficient to take for $N_1$ the image of the set $N$ under the projection $(\lambda,\mu)\mapsto\lambda$ onto the first coordinate, and for $N_2$ the image of the set $N$ under the projection $(\lambda,\mu)\mapsto\mu$ onto the second coordinate.
 \end{proof}

Let $K$ be a closed subset of $\overline{\mathbb C}^2$ or $\overline{\mathbb C}$ and $\mathoo B$ be a unital Banach algebra. We denote by $\mathoo O(K,\mathoo B)$ the set of all analytic\footnote{A function $f:\,U\subset\overline{\mathbb C}\to\mathoo B$ is called \emph{analytic at infinity} if $f$ can be expanded in a power series $f(\lambda)=\sum_{k=0}^\infty\frac{c_k}{\lambda^k}$ in a neighbourhood of infinity.
A function $f:\,U\subset\overline{\mathbb C}^2\to\mathbb C$ is called \emph{analytic at} $(\infty,\infty)$ if $f$ can be expanded in a power series $f(\lambda,\mu)=\sum_{k,m=0}^\infty\frac{c_{km}}{\lambda^k\mu^m}$ in a neighbourhood of $(\infty,\infty)$.
A function $f:\,U\subset\overline{\mathbb C}^2\to\mathbb C$ is called \emph{analytic at} $(\lambda_*,\infty)$, $\lambda_*\in\mathbb C$ if $f$ can be expanded in a power series $f(\lambda,\mu)=\sum_{k,m=0}^\infty\frac{c_{km}(\lambda-\lambda_*)^k}{\mu^m}$ in a neighbourhood of $(\lambda_*,\infty)$.}~\cite{Shabat2:eng,Greene-Krantz} functions $f\separ U\to\mathoo B$, where $U$ is an open neighbourhood of the set $K$ (it is implied that the neighbourhood $U$ may depend on $f$). Two functions $f_1\separ U_1\to\mathoo B$ and $f_2\separ U_2\to\mathoo B$ are called {\it equivalent} if there exists an open neighbourhood $U\subseteq U_1\cap U_2$ of the set $K$ such that $f_1$ and $f_2$ coincide on $U$, i.~e.~$f_1(\lambda)=f_2(\lambda)$ for all
$\lambda\in U$. It can be easily shown that this is really an equivalence relation.
Thus, strictly speaking, elements of $\mathoo O(K,\mathoo B)$ are classes of equivalent functions. The notation $\mathoo O(K,\mathbb C)$ is abbreviated to $\mathoo O(K)$.

\begin{proposition}\label{e:O is an algebra:Z2}
The set $\mathoo O(K,\mathoo B)$ is an algebra with respect to pointwise operations with the unit $u(\lambda)=1$, $\lambda\in U\supset K$.
\end{proposition}

\begin{proposition}\label{p:spec in O:Z2}
{\rm(a)} For $f\in\mathoo O(K,\mathoo B)$, the following conditions are equivalent: the function $f$ is invertible in the algebra $\mathoo O(K,\mathoo B)${\rm;} the element $f(\lambda)\in\mathoo B$ is invertible at all points $\lambda\in K${\rm;} the element $f(\lambda)$ is invertible at all points $\lambda\in U$, where $U\supset K$ is some open set. {\rm(b)} The spectrum of a function $f\in\mathoo O(K,\mathoo B)$ in the algebra $\mathoo O(K,\mathoo B)$ is given by the formula
\begin{equation*}
\bigcup_{\lambda\in K}\sigma\bigl(f(\lambda)\bigr).
\end{equation*}
\end{proposition}

We recall the definition of the natural topology on the algebra $\mathoo O(K,\mathoo B)$~\cite[ch.~1, \S~4.1]{Bourbaki_Theories_Spectrales:eng}.

For each open set $U\supset K$, we denote by $\mathoo O(U,\mathoo B)$ the linear space of all
analytic functions $f\separ U\to\mathoo B$. We endow $\mathoo O(U,\mathoo B)$ with the \emph{topology of compact convergence}~\cite[ch.~X, \S~3.6]{Bourbaki_Topologie_Generale-10:eng}, \cite[ch. III,
\S~3]{Schaefer:eng}. A fundamental system of neighbourhoods of zero in this topology is formed by the sets
$T(N,\delta)=\{\,f\separ\Vert f(N)\Vert<\delta\,\}$, where $N\subset U$ is compact and
$\delta>0$; clearly, when $N$ enlarges, the neighbourhood $T(N,\delta)$ shrinks; therefore, it is enough to consider
only those sets $N$, the interior of which contains $K$. There are evident canonical
mappings $g_U:\,\mathoo O(U,\mathoo B)\to\mathoo O(K,\mathoo B)$.
The mappings $g_U$ are not always injective. Nevertheless, by misuse of language, we will regard
$\mathoo O(U,\mathoo B)$ as subspaces of the space $\mathoo O(K,\mathoo B)$.

We endow $\mathoo O(K,\mathoo B)$ with the inductive topology~\cite[ch. II, \S~6]{Schaefer:eng} induced by the mappings $g_U$ (one may restrict himself to a decreasing sequence of open sets $U\supset K$).
A fundamental system of neighbourhoods of zero in $\mathoo O(K,\mathoo B)$ consists of all balanced, absorbent, and convex sets $W\subseteq\mathoo O(K,\mathoo B)$ such that the inverse image $g_U^{-1}(W)$ is a neighbourhood of zero in $\mathoo O(U,\mathoo B)$. Thus, for all $U\supset N\supset K$ such that the interior of the compact set
$N$ contains $K$, the inverse image  $g_U^{-1}(W)$ must contain the set
$T(N,\delta)=\{\,f\separ \Vert f(N)\Vert<\delta\,\}\subset\mathoo O(U,\mathoo B)$.

We recall~\cite[ch.~2, Theorem~6.1]{Schaefer:eng} that a linear mapping $\varphi\separ
\mathoo O(K,\mathoo B)\to\mathbb E$, where $\mathbb E$ is a Banach space, is continuous if and only if all the compositions
$\varphi\circ g_U\separ \mathoo O(U,\mathoo B)\to\mathbb E$,
where $U\supset K$, are continuous, i.~e. for any neighbourhood $W\subseteq\mathbb E$ of zero there
exist a compact set $N\subset U$ and a number
$\delta>0$ such that the interior of $N$ contains $K$ and $\varphi\circ g_U\bigl(T(N,\delta)\bigr)\subseteq W$. We note that since
$\mathbb E$ is a Banach space, it is sufficient to restrict ourselves to the consideration of $\varepsilon$-neighbourhoods of zero for $W$. Below, by misuse of language, we denote $\varphi\circ g_U$ by the abbreviated symbol $\varphi$.

 \begin{proposition}[{\rm\cite[Proposition 1.3]{Ichinose75}}]\label{p:dense}
Let $U_1\subseteq\overline{\mathbb C}$ and $U_2\subseteq\overline{\mathbb C}$ be open sets. Then the natural image of the algebraic tensor product $\mathoo O(U_1)\otimes\mathoo O(U_2)$ is everywhere dense in $\mathoo O(U_1\times U_2)$.
 \end{proposition}

\section{Pseudo-resolvents and functional calculus}\label{s:pseudo-resolvents}
We call a mapping that converts functions to operators (or transformators) a \emph{functional calculus}.
Of course, the most interesting are functional calculi that possess special properties (for example, morphisms of algebras).

A \emph{pseudo-resolvent} is a function that takes values in a Banach algebra and satisfies the Hilbert identity~\eqref{e:Hilbert identity:alg}, like a resolvent. Every pseudo-resolvent generates a functional calculus (Theorems~\ref{t:functional_calculas:infty in rho} and~\ref{t:functional_calculas:infty in
sigma}) which is a morphism of algebras and possesses the property of preserving the spectrum (Theorem~\ref{t:spectral mappping}).

Let $\mathoo B$ be a Banach algebra and $U\subseteq\mathbb C$ be a subset.
A function (family) $\lambda\mapsto
R_\lambda$ defined on $U$ and taking values in $\mathoo B$ is called~\cite[ch.~5, \S~2]{Hille-Phillips:eng} a \emph{pseudo-resolvent} if it satisfies the \emph{Hilbert identity}
 \begin{equation}\label{e:Hilbert identity:pseudo}
R_\lambda-R_\mu=-(\lambda-\mu)R_\lambda R_\mu,\qquad\lambda,\mu\in U.
 \end{equation}
A pseudo-resolvent is called~\cite[p.~103]{Baskakov2004:eng} \emph{maximal}\label{text:def of max pseudo-resolvent} if it cannot be extended to a larger set with the preservation of the Hilbert identity~\eqref{e:Hilbert identity:pseudo}. Below
(Theorem~\ref{t:pseudoresovent}) we will see that every pseudo-resolvent can be extended to a unique maximal one. The domain $\rho(R_{(\cdot)})$ of the maximal extension of a pseudo-resolvent $R_{(\cdot)}$ is called a \emph{regular
set} of the original pseudo-resolvent. The complement
$\sigma(R_{(\cdot)})$ of the regular set $\rho(R_{(\cdot)})$ is called~\cite{Arendt},
\cite[p.~103]{Baskakov2004:eng} the \emph{singular set}.

\begin{example}\label{ex:pseudo-resolvent}
The examples of pseudo-resolvents are: (a) the resolvent of an element of a unital Banach algebra (Proposition~\ref{p:res of A is max res}); (b) a constant function $\lambda\mapsto N$, where $N\in\mathoo B$ is an arbitrary element whose square equals zero (Proposition~\ref{c:pseudo in an infinite point}); (c) in particular, the identically zero function; (d) the resolvent of a closed linear operator~\cite[Theorem 5.8.1]{Hille-Phillips:eng}; (e) the resolvent of a linear relation~\cite{Arens61,Baskakov2004:eng,Bichegkuev:eng,Cross,Favini-Yagi,Haase2006,Kurbatova-PMM09:eng}; this example is the most general, because every pseudo-resolvent is a resolvent of some linear relation~\cite[Theorem 5.2.4]{Baskakov2004:eng}, \cite[Proposition A.2.4]{Haase2006}; (f) direct sums of pseudo-resolvents from the previous examples.
 \end{example}

A simple example of a maximal pseudo-resolvent is given in the following proposition.

 \begin{proposition}[\mbox{\rm\cite[Proposition 17]{Kurbatov-Kurbatova-IM-2015:eng}}]\label{p:res of A is max res}
The resolvent of an arbitrary element $A\in\mathoo B$ is a maximal pseudo-resolvent, i.~e. it cannot be extended to a set larger than $\rho(A)$ with the preservation of the Hilbert identity~\eqref{e:Hilbert
identity:pseudo}.
 \end{proposition}

We note that identity~\eqref{e:Hilbert identity:pseudo} can be equivalently written in the form
(here $\mathbf1$ is an adjoint unit if the original algebra has no unit)
 \begin{equation}\label{e:Hilbert identity':alg}
R_\lambda\bigl(\mathbf1+(\lambda-\mu)R_\mu\bigr)=R_\mu.
 \end{equation}

Below in this Section, we will adjoin a unit to $\mathoo B$ when it has no unit.

 \begin{proposition}[\mbox{\rm\cite[Corollary 1 of Theorem 5.8.4]{Hille-Phillips:eng}}]\label{p:mathbf1+(lambda-mu)}
Let $R_\lambda,R_\mu\in\mathoo B$ be two commuting elements that satisfy identity~\eqref{e:Hilbert identity:pseudo}. Then the element $\mathbf1+(\lambda-\mu)R_\mu\in\widetilde{\mathoo B}$ is necessarily invertible.
 \end{proposition}

 \begin{theorem}[\mbox{\rm\cite[Theorem 5.8.6]{Hille-Phillips:eng}}]\label{t:pseudoresovent}
Every pseudo-resolvent whose domain contains at least one point $\mu\in\mathbb C$ can be extended to a maximal pseudo-resolvent{\rm;} this extension is unique. The domain of the maximal extension is the set of all $\lambda\in\mathbb C$ such that the element
$\mathbf1+(\lambda-\mu)R_\mu$ is invertible in $\widetilde{\mathoo B}$. This extension can be defined by the formula
 \begin{equation}\label{e: restore res}
R_\lambda=R_\mu\bigl(\mathbf1+(\lambda-\mu)R_\mu\bigr)^{-1}=\bigl(\mathbf1+(\lambda-\mu)R_\mu\bigr)^{-1}R_\mu.
 \end{equation}
 \end{theorem}

We will denote the original pseudo-resolvent and its continuation to a maximal pseudo-resolvent by the same symbol
$R_{(\cdot)}$. Moreover, we will generally assume that all pseudo-resolvents under consideration are already extended to maximal pseudo-resolvents.

\begin{corollary}[{\rm\cite[Theorem 5.8.2]{Hille-Phillips:eng}, \cite[ch.~6, \S~1]{Cross}}]\label{c:holom Rlambda} The domain of a maximal pseudo-resolvent is an open set and the maximal pseudo-resolvent is an analytic function {\rm(}with values in $\mathoo B${\rm)}.\footnote{We accept that the domain of an analytic function may be disconnected.} More precisely, in a neighbourhood of any point $\mu\in\rho(R_{(\cdot)})$, the maximal pseudo-resolvent admits the power series expansion
\begin{equation*}
R_\lambda=\sum_{n=0}^\infty(\mu-\lambda)^nR_\mu^{n+1}.
\end{equation*}
 \end{corollary}

 \begin{proposition}[{\rm\cite[Theorem 5.8.3]{Hille-Phillips:eng}}]\label{p:pseudo is res:alg}
A pseudo-resolvent $R_{(\cdot)}$ in a Banach algebra $\mathoo B$ is a resolvent of some
element $A\in\mathoo B$ if and only if $\mathoo B$ is unital and the element
$R_\mu$ is invertible for at least one {\rm(}and, consequently, for all{\rm)}
$\mu\in\rho(R_{(\cdot)})$.
A pseudo-resolvent $R_{(\cdot)}$ in $\mathoo B(X)$, where $X$ is a Banach space, is a resolvent of some
unbounded operator $A:\,D(A)\subset X\to X$ if and only if the operator $R_\mu:\,X\to\Image R_\mu$ is invertible for at least one {\rm(}and, consequently, for all{\rm)} $\mu\in\rho(R_{(\cdot)})$.
In this case $A=\lambda\mathbf1-(R_{\lambda})^{-1}$ for all $\lambda\in\rho(R_{(\cdot)})$.
 \end{proposition}

In a similar way, a linear relation can also be recovered from the value of its resolvent at one point. Thus, the resolvent contains all information about a linear relation or an operator that generates it. On the other hand, the conditions on unbounded operators and linear relations are often imposed in terms of their resolvents (the nonemptiness of the resolvent set, the estimate of decay rate at infinity etc.). Besides, functions of linear relations and unbounded operators are often expressed directly via their resolvents. For this reason, the resolvent can be considered as a more fundamental object than an operator or relation that generates it. This is the viewpoint we adhere to in this article.

\medskip
We fix a pseudo-resolvent $R_{(\cdot)}$. We denote by $\mathoo B_R$ the
smallest closed subalgebra of the algebra $\mathoo B$ that contains all elements $R_\lambda$,
$\lambda\in\rho(R_{(\cdot)})$, of the extension of the pseudo-resolvent $R_{(\cdot)}$ to a maximal pseudo-resolvent.

 \begin{proposition}[{\rm\cite[Proposition 21]{Kurbatov-Kurbatova-IM-2015:eng}}]\label{p:B_R is commutative}
The algebra $\mathoo B_R$ coincides with the closure of the linear span of the family of all elements
$R_\lambda$, $\lambda\in\rho(R_{(\cdot)})$, and is commutative.
 \end{proposition}

If the algebra $\mathoo B_R$ does not contain the unit of the algebra $\mathoo B$ (this is certainly the case
if $\mathoo B$ is not unital), then we will, as usual, denote by $\widetilde{\mathoo B}_R$
the algebra $\mathoo B_R$ with an adjoint unit from $\mathoo B$\footnote{Adjoining $\mathbf1\in\mathoo B$ to $\mathoo B_R$ we obtain a closed subalgebra because the sum of a closed subspace and a one-dimensional subspace is a closed subspace.} (or the algebra $\widetilde{\mathoo B}$ with an adjoint unit). If $\mathoo B_R$ contains the unit of the algebra $\mathoo B$, the symbol
$\widetilde{\mathoo B}_R$ is understood to be $\mathoo B_R$.

 \begin{proposition}[{\rm\cite[Theorem 22]{Kurbatov-Kurbatova-IM-2015:eng}}]\label{p:B_R=B(R)}
The subalgebra $\widetilde{\mathoo B}_R$ is commutative and full.
 \end{proposition}

 \begin{proposition}\label{p:chi_infty}
If a character $\chi$ of the algebra $\mathoo B_R$ equals zero at least at one element $R_\mu$, $\mu\in\rho(R_{(\cdot)})$, then it is identically equal to zero on $\mathoo B_R$, i.~e. coincides with $\chi_0$.
 \end{proposition}
 \begin{proof}
The proof follows from formula~\eqref{e: restore res} and the description of $\mathoo B_R$
{\rm(}Proposition~\ref{p:B_R is commutative}{\rm)} as the closure of the linear span of the family $R_\lambda$, $\lambda\in\rho(R_{(\cdot)})$.
 \end{proof}

We say that a sequence of maximal pseudo-resolvents $R_{n,\,(\cdot)}$
\emph{converges}\label{page:conv of Res} to a maximal pseudo-resolvent $R_{(\cdot)}$ if there exists a point
$\mu\in\mathbb C$ such that all the pseudo-resolvents $R_{n,\,(\cdot)}$ are defined at $\mu$ (for $n$ sufficiently large) and the sequence $R_{n,\,\mu}$ converges to $R_{\mu}$ in norm, cf.~\cite[Theorem 2.25]{Kato:eng}. The following lemma shows that this definition does not depend on the choice of the point $\mu\in\mathbb C$.

 \begin{lemma}\label{l:limit of pseudo-resolvents}
Let a sequence $R_{n,\,(\cdot)}$ of maximal pseudo-resolvents converge to a maximal pseudo-resolvent $R_{(\cdot)}$ at a point $\mu\in\mathbb C$ {\rm(}it is assumed that the pseudo-resolvents $R_{n,\,\mu}$ are defined at the point $\mu$ for all $n$ large enough{\rm)}. Then for any point
$\lambda\in\rho(R_{(\cdot)})$, the sequence $R_{n,\,\lambda}$ is defined
for all $n$ large enough and converges to $R_{\lambda}$ in norm.

Moreover, given a compact set $\Gamma\subset\rho\bigl(R_{(\cdot)}\bigr)$, the elements $R_{n,\,\lambda}$ are defined
for $n$ large enough at all $\lambda\in\Gamma$ and converges to $R_{\lambda}$ uniformly with respect to $\lambda\in\Gamma$.
 \end{lemma}

 \begin{proof}
Let $R_{n,\,\mu}$ converge to $R_{\mu}$. By Theorem~\ref{t:pseudoresovent}, the element
$\mathbf1+(\lambda-\mu)R_\mu$ is invertible for all $\lambda\in\Gamma$. Since the function $\lambda\mapsto\mathbf1+(\lambda-\mu)R_\mu$ is continuous, from Theorem~\ref{t:Neumann} it follows that
\begin{equation*}
\min_{\lambda\in\Gamma}\bigl\Vert\bigl(\mathbf1+(\lambda-\mu)R_\mu\bigr)^{-1}\bigr\Vert>0.
\end{equation*}
Since $R_{n,\,\mu}$ converges to $R_{\mu}$, the sequence $\lambda\mapsto\mathbf1+(\lambda-\mu)R_{n,\,\mu}$ converges to $\lambda\mapsto\mathbf1+(\lambda-\mu)R_\mu$ uniformly with respect to $\lambda\in\Gamma$. Therefore, again by Theorem~\ref{t:Neumann}, the elements $\mathbf1+(\lambda-\mu)R_{n,\,\mu}$ are invertible for all $\lambda\in\Gamma$ provided $n$ is large enough; in this case, by estimate~\eqref{e:(A-B)-1-A-1}, the inverses also converge uniformly. It remains to apply formula~\eqref{e: restore res}.
 \end{proof}

We note that the limit of a sequence of resolvents of bounded operators can be the resolvent of a non-bounded operator, see~\cite[Lemma 7]{Reed-Simon:JFA73}.

 \begin{proposition}[{\rm\cite[Theorem 5.9.2]{Hille-Phillips:eng}}]\label{c:pseudo in an infinite point}
Let a pseudo-resolvent $R_{(\cdot)}$ admit an analytic continuation in a neighbourhood of the point $\infty$.\footnote{We recall~\cite[p.~107]{Greene-Krantz} that the possibility of an analytic continuation of $f$ in a neighbourhood of the point $\infty$ is equivalent to the existence of a bounded analytic continuation of $f$ in a deleted neighbourhood of~$\infty$.}
Then there exist elements $P,A,N\in\mathoo B_R$ such that
\begin{equation*}
N^{2}=\mathbf0,\qquad P^{2}=P,\qquad AP=PA=A,\qquad NP=PN=\mathbf0
\end{equation*}
and the expansion of the pseudo-resolvent into the Laurent series with centre $\infty$ has the form
\begin{equation}\label{e:sigma_R:infty}
R_\lambda=-N+\frac P\lambda+\frac
A{\lambda^2}+\frac{A^2}{\lambda^3}+\frac{A^3}{\lambda^4}+\ldots\,.
\end{equation}
 \end{proposition}

We call the \emph{extended regular set} $\bar\rho(R_{(\cdot)})\subseteq\overline{\mathbb{C}}$ of a pseudo-resolvent $R_{(\cdot)}$ (in the algebra $\mathoo B$) either the regular set $\rho(R_{(\cdot)})$ or the union $\rho(R_{(\cdot)})\cup\{\infty\}$; more precisely, we add the point $\infty$ to $\bar\rho(R_{(\cdot)})$ if the algebra $\mathoo B$ is unital, the regular set $\rho(R_{(\cdot)})$ contains a (deleted) neighbourhood of $\infty$, and $\lim_{\lambda\to\infty}\lambda
R_\lambda=\mathbf1$. We call the \emph{extended singular set} of the pseudo-resolvent the complement $\bar\sigma(R_{(\cdot)})=\overline{\mathbb{C}}\setminus\bar\rho(R_{(\cdot)})$ of the extended regular set.

 \begin{proposition}\label{p:1 in B_R}
The following properties of a maximal pseudo-resolvent are equivalent:
 \begin{itemize}
 \item [\rm(a)] $\infty\in\bar\rho(R_{(\cdot)})${\rm;}
 \item [\rm(b)] the maximal pseudo-resolvent is the resolvent of some element $A\in\mathoo B_R$ {\rm(}see Proposition~\ref{p:pseudo is res:alg}{\rm)}{\rm;}
 \item [\rm(c)] the algebra $\mathoo B$ is unital and the subalgebra $\mathoo B_R$ contains the unit of the algebra $\mathoo B$.
 \end{itemize}
 \end{proposition}
 \begin{proof} The equivalence of (a) and (b) is proved in~\cite[Proposition 23]{Kurbatov-Kurbatova-IM-2015:eng}.

Let assumption (b) be fulfilled, i.~e. $R_\lambda=(\lambda\mathbf1-A)^{-1}$, where $A\in\mathoo
B_R$. Then, by virtue of Theorem~\ref{t:Neumann}, the power series expansion
\begin{equation*}
(\lambda\mathbf1-A)^{-1}=\frac{\mathbf1}\lambda+\frac{A}{\lambda^2}+\frac{A^2}{\lambda^3}+\ldots\,
\end{equation*}
holds in a neighbourhood of infinity, which shows that $\mathbf1\in\mathoo B_R$, i.~e., assumption (c) holds.

Let assumption (c) be fulfilled, i.~e. the subalgebra $\mathoo B_R$ contain the unit of the algebra $\mathoo B$. There are no identically zero characters on a unital commutative algebra, because $\chi(\mathbf1)=1$.
Therefore, by Proposition~\ref{p:chi_infty}, $\chi(R_\mu)\neq0$ for all $\chi\in\widetilde{X}(\mathoo B_R)$ and $\mu\in\rho(R_{(\cdot)})$.
Hence, by Theorem~\ref{t:Gelfand}, all values
$R_\mu$ of the pseudo-resolvent are invertible. Then, by virtue of Proposition~\ref{p:pseudo is res:alg},
the pseudo-resolvent is the resolvent of some element $A\in\mathoo B_R$, i.~e. assumption (b) holds.
 \end{proof}

Below in this Section, we assume that $X$ is a Banach space and we are given a maximal pseudo-resolvent $R_{(\cdot)}$ in $\mathoo B(X)$.

Let $\sigma$ and $\Sigma$ be two disjoint closed subsets of $\overline{\mathbb C}$. A contour $\Gamma$ is called~\cite[ch. V, \S~5.2]{Hille-Phillips:eng} an \emph{oriented envelope} of the set $\sigma$ \emph{with respect} to the set $\Sigma$ if $\Gamma$ is an oriented boundary of an open set $U$ that contains $\sigma$ and is disjoint from $\Sigma$.
Thus, $\Gamma$ surrounds the set $\sigma$ in the counterclockwise direction and surrounds the set $\Sigma$ in the clockwise direction.

\begin{figure}[htb]
\unitlength=1.mm
\begin{center}
\begin{picture}(50,30)
\put(22.5,17){\circle*{10}} \put(21,14){\circle*{10}} \put(18,17){\circle*{10}}
\put(20,15){\oval(34,22)} \put(37,15){\vector(0,1){3}} \put(32,15){$\Gamma$}
\put(20,8){$\sigma$} \put(40,12){\boxed{$\phantom{aa}f\phantom{aa}$}}
\end{picture}\hfil
\begin{picture}(65,30)
\put(27,19){\circle*{30}} \put(26,16){\circle*{30}} \put(25,14){\circle*{30}}
\put(23,17){\circle*{30}} \put(50,15){\circle{20}} \put(43,14){\vector(0,1){2}} \put(39,
12){$\Gamma$} \put(25,7){$\bar\sigma$} \put(50,12){\boxed{$f$}}
\end{picture}
\caption{The contour $\Gamma$ is an oriented envelope of the set $\sigma$
with respect to $\infty$ and the set \boxed{$f$} (left); the contour
$\Gamma$ is an oriented envelope of the set $\bar\sigma$ and the point $\infty$
with respect to the set \boxed{$f$} (right). }\label{f:Gamma0}
\end{center}
 \end{figure}
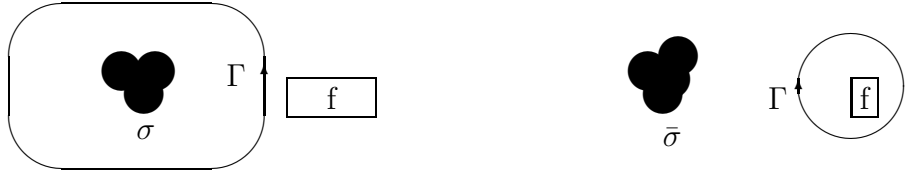

\begin{theorem}\label{t:functional_calculas:infty in rho}
Assume that $\infty\notin\bar\sigma(R_{(\cdot)})$. We define the mapping\/ $\varphi\separ\mathoo
O\bigl(\sigma(R_{(\cdot)})\bigr)\to\mathoo B_R$ by the formula
\begin{equation}\label{e:functional calculus_1}
\varphi(f)=\frac1{2\pi i}\int_\Gamma f(\lambda)R_\lambda\,d\lambda,
\end{equation}
where $\Gamma$ {\rm(}see left fig.~{\rm\ref{f:Gamma0})} is an oriented envelope of the singular set $\sigma(R_{(\cdot)})$ with respect to the point $\infty$ and the complement of the domain of the function $f$. We assert that $\varphi$ is a continuous morphism of unital algebras.

The morphism $\varphi$ maps the function $u(\lambda)=1$ to the identity operator $\mathbf1_X$. The function $v_1(\lambda)=\lambda$ is mapped by $\varphi$ to the operator $A\in\mathoo B(X)$ that generates the maximal pseudo-resolvent $R_{(\cdot)}$ in accordance with Proposition~{\rm\ref{p:1 in B_R}}{\rm;} and the function $r_{\lambda_0}(\lambda)=\frac1{\lambda_0-\lambda}$, where $\lambda_0\in\rho(R_{(\cdot)})$, is mapped by $\varphi$ to $R_{\lambda_0}$.
\end{theorem}

 \begin{proof} The proof is analogous to that of the theorem on analytic functional calculus for bounded operators~\cite[ch.~1, \S~4, Theorem 3]{Bourbaki_Theories_Spectrales:eng}, \cite[Theorem 5.2.5]{Hille-Phillips:eng}, \cite[Theorem~10.27]{Rudin:eng}.
\end{proof}

When it is desirable to stress that the functional  calculus $\varphi$ considered in Theorem~\ref{t:functional_calculas:infty in rho} is generated by the resolvent of an operator $A\in\mathoo B(X)$, we will use the notation $R_{A,\,\lambda}$ instead of $R_\lambda$, the notation $\varphi_A$ instead of $\varphi$, and the notation $f(A)$ instead of $\varphi(f)$.

\begin{theorem}\label{t:functional_calculas:infty in sigma}
Assume that $\infty\in\bar\sigma(R_{(\cdot)})$. We define the mapping\/ $\varphi\separ \mathoo
O\bigl(\bar\sigma(R_{(\cdot)})\bigr)\to\widetilde{\mathoo B}_R$ by the formula
\begin{equation}\label{e:functional calculus:infty in sigma}
\varphi(f)=\frac1{2\pi i}\int_\Gamma f(\lambda)R_\lambda\,d\lambda+f(\infty)\mathbf1,
\end{equation}
where $\Gamma$ is an oriented envelope of the extended singular set $\bar\sigma(R_{(\cdot)})$
with respect to the complement of the domain of $f$ {\rm(}see the right fig.~{\rm\ref{f:Gamma0})}. We assert that $\varphi$ is a continuous morphism of unital algebras.

The morphism $\varphi$ maps the function $u(\lambda)=1$ to the identity operator $\mathbf1$. The function
$r_{\lambda_0}(\lambda)=\frac1{\lambda_0-\lambda}$, where $\lambda_0\in\rho(R_{(\cdot)})$,
is mapped by $\varphi$ to $R_{\lambda_0}$.
\end{theorem}

 \begin{proof}
The proof is analogous to that of the theorem on analytic functional calculus for unbounded operators~\cite[Theorem~5.11.2]{Hille-Phillips:eng}.
\end{proof}

When it is desirable to stress that the functional  calculus $\varphi$ considered in Theorem~\ref{t:functional_calculas:infty in sigma} is generated by a pseudo-resolvent $R_{(\cdot)}$, we will use the notation $\varphi_{R_{(\cdot)}}$ instead of $\varphi$ and the notation $f(R_{(\cdot)})$ instead of $\varphi(f)$.

A unified notation for formulae~\eqref{e:functional calculus_1} and~\eqref{e:functional calculus:infty in sigma} is suggested in~\cite{Hille-Phillips:eng}:
\begin{equation*}
\varphi(f)=\frac1{2\pi i}\int_\Gamma f(\lambda)R_\lambda\,d\lambda+\delta f(\infty)\mathbf1,
\end{equation*}
where $\delta=0$ if $\Gamma$ does not enclose $\infty$, and $\delta=1$ if $\Gamma$ encloses $\infty$.

The following theorem is a version of the spectral mapping theorem for the case of pseudo-resolvents.

\begin{theorem}\label{t:spectral mappping}
For any function $f\in\mathoo O\bigl(\bar\sigma(R_{(\cdot)})\bigr)$ we have the equality
$$
\sigma_{\widetilde{\mathoo B}}\bigl(\varphi(f)\bigr)=\sigma_{\widetilde{\mathoo
B}_R}\bigl(\varphi(f)\bigr)=\{\,f(\lambda)\:\lambda\in\bar\sigma(R_{(\cdot)})\,\}.
$$
\end{theorem}
\begin{proof}
The proof is analogous to that of the spectral mapping theorem for unbounded operators~\cite[теорема 5.12.1]{Hille-Phillips:eng} and for linear relations~\cite[теорема 5.2.17]{Baskakov2004:eng}.
\end{proof}

An analytic functional calculus for bounded operators was created in~\cite{Dunford43a,Dunford43b,Taylor43}.
It was carried over to unbounded operators in~\cite{Dunford39,Hille-Phillips:eng} and to linear relations in~\cite{Arens61,Baskakov2004:eng,Cross,Favini-Yagi,Haase2006}.

\section{Extended tensor products}\label{s:extended tensor product}
The notion of an extended tensor product is a generalization of the notion of a completion of an algebraic tensor product with respect to a uniform cross-norm~\cite{Grothendieck,Schatten,Defant-Floret,Diestel-Fourie-Swart,Helemskii:eng}.
It enables one to extend some constructions which are natural for the usual tensor products to supplementary applications. We recall an example which is the most important for our applications. It is known (see, e.~g.,~\cite{Graham}) that in the case of finite-dimensional Banach spaces $X$ and $Y$, the space $\mathoo{B}(Y,X)$ can be identified with the tensor product $X\otimes Y^*$. If $X$ and $Y$ are infinite-dimensional, then $X\otimes Y^*$ corresponds only to the subspace of $\mathoo{B}(Y,X)$ consisting of operators that have a finite-dimensional image. Therefore the completion of $X\otimes Y^*$ with respect to any reasonable norm cannot coincide with the whole $\mathoo{B}(Y,X)$. Nevertheless, $\mathoo{B}(Y,X)$ can be represented (see example~\ref{ex:TP of cat}(e) below) as an extended tensor product $X\boxtimes Y^*$ which enables one to treat it almost as a usual tensor product.
The exposition in this Section is based on~\cite{Kurbatov-Kurbatova-IM-2015:eng}.

We denote by $X\otimes Y$ the usual tensor product of linear spaces
$X$ and $Y$. In order to distinguish $X\otimes Y$ from its extensions, we call $X\otimes Y$ an
\emph{algebraic tensor product}.

Let $X$ and $Y$ be Banach spaces. A norm $\alpha(\cdot)=\Vert\cdot\Vert_\alpha$
on $X\otimes Y$ is called a \emph{cross-norm} if
\begin{equation*}
\Vert x\otimes y\Vert_\alpha=\Vert x\Vert\cdot\Vert y\Vert
\end{equation*}
for all $x\in X$ and $y\in Y$. We denote by $X\overline{\otimes}_\alpha\,Y$ the completion
of the tensor product $X\otimes Y$ by the cross-norm $\alpha$.

Every element $v^*=\sum_{l=1}^m x_l^*\otimes y_l^*\in X^*\otimes Y^*$ induces the linear functional
\begin{equation}\label{e:nat imbedding j}
v^*\separ\sum_{k=1}^n x_k\otimes y_k\mapsto\sum_{l=1}^m\sum_{k=1}^n\langle
x_k,x_l^*\rangle\cdot\langle y_k,y_l^*\rangle
\end{equation}
on the space $X\otimes Y$. We define the norm $\alpha^*$ on $X^*\otimes Y^*$, \emph{conjugate to the cross-norm} $\alpha$, by the formula
\begin{equation*}
\Vert v^*\Vert_{\alpha^*}=\sup\{\,|\langle v,v^*\rangle|\separ v\in X\otimes Y,\,\Vert
v\Vert_\alpha\le1\,\}.
\end{equation*}
A cross-norm $\alpha$ is called $*$-\emph{uniform} if $\alpha^*$ is finite and is a cross-norm.

The space $\mathoo B(X)\otimes\mathoo B(Y)$ has the natural structure of an algebra. Every element $T=\sum_{l=1}^m A_l\otimes
B_l\in\mathoo B(X)\otimes\mathoo B(Y)$ induces the linear operator
\begin{equation*}
T\separ\sum_{k=1}^n x_k\otimes y_k\mapsto\sum_{l=1}^m\sum_{k=1}^n(A_lx_k)\otimes (B_ly_k)
\end{equation*}
in $X\otimes Y$. A cross-norm $\alpha$ on the space $X\otimes Y$ induces the norm $\tilde\alpha$ of the operator $T\in\mathoo B(X)\otimes\mathoo B(Y)$ by the formula
\begin{equation*}
\Vert T\Vert_{\tilde\alpha}=\sup\{\,\Vert Tv\Vert\separ v\in X\otimes Y,\,\Vert
v\Vert_\alpha\le1\,\}.
\end{equation*}
A cross-norm $\alpha$ is called~\cite{Schatten} \emph{uniform} if $\tilde\alpha$
is finite and is a cross-norm. Every uniform cross-norm is $*$-uniform, see~\cite{Simon-1973}.

\medskip
Let $X$ and $Y$ be Banach spaces. We call an \emph{extended tensor product}~\cite{Kurbatov-Kurbatova-IM-2015:eng} of $X$
and $Y$ a collection consisting of three objects: a Banach space $X\boxtimes Y$ (which we briefly refer to as
the extended tensor product) and two (not necessarily closed) full unital subalgebras $\mathoo B_0(X)$ and $\mathoo B_0(Y)$ of the algebras $\mathoo B(X)$ and $\mathoo B(Y)$ respectively that satisfy assumptions (A), (B), and (C) listed below.
 \begin{itemize}
 \item [\rm(A)] We are given a linear mapping
$j$ from the algebraic tensor product $X\otimes Y$ to $X\boxtimes Y$. In the sequel, we denote
$j(x\otimes y)$ by the symbol $x\boxtimes y$. It is assumed that
\begin{equation}\label{e:C2}
\Vert x\boxtimes y\Vert_{X\boxtimes Y}=\Vert x\Vert_{X}\cdot\Vert y\Vert_{Y}
\end{equation}
for all $x\in X$ and $y\in Y$.
 \item [\rm(B)] We are given a linear mapping $J$
from the algebraic tensor product $X^*\otimes Y^*$ to $(X\boxtimes Y)^*$. In the sequel, we denote
$J(x^*\otimes\,y^*)$ by the symbol $x^*\boxtimes\,y^*$. It is assumed that
\begin{equation}\label{e:D1}
\langle x\boxtimes y,x^*\boxtimes\,y^*\rangle=\langle x,x^*\rangle\langle y,y^*\rangle
\end{equation}
for all $x^*\in X^*$, $y^*\in Y^*$, $x\in X$, and $y\in Y$, and
\begin{equation}\label{e:D2}
\Vert x^*\boxtimes y^*\Vert_{(X\boxtimes Y)^*}=\Vert x^*\Vert_{X^*}\cdot\Vert
y^*\Vert_{Y^*}
\end{equation}
for all $x^*\in X^*$ and $y^*\in Y^*$.
 \item [\rm(C)] We are given a morphism $\mathfrak J$ of unital algebras from the algebraic tensor product $\mathoo B_0(X)\otimes\mathoo B_0(Y)$ to $\mathoo B(X\boxtimes Y)$. In the sequel, we denote $\mathfrak J(A\otimes B)$ by the symbol
$A\boxtimes B$. It is assumed that
\begin{equation}\label{e:C1}
\bigl(A\boxtimes B\bigr)(x\boxtimes y)=(Ax)\boxtimes(By)
\end{equation}
for all $A\in\mathoo B_0(X)$, $B\in\mathoo B_0(Y)$, $x\in X$, and $y\in Y$, and
\begin{equation}\label{e:D3}
\bigl(A\boxtimes B\bigr)^*(x^*\boxtimes y^*)=(A^*x^*)\boxtimes(B^*y^*)
\end{equation}
for all $A\in\mathoo B_0(X)$, $B\in\mathoo B_0(Y)$, $x^*\in X^*$, and $y^*\in Y^*$, and
\begin{equation}\label{e:A otimes B le}
\Vert A\boxtimes B\Vert_{\mathoo B(X\boxtimes Y)}=\Vert A\Vert_{\mathoo B(X)}\cdot\Vert B\Vert_{\mathoo B(Y)}
\end{equation}
for all $A\in\mathoo B_0(X)$ and $B\in\mathoo B_0(Y)$.
 \end{itemize}

 \begin{example}\label{ex:TP of cat}
We recall~\cite{Kurbatov-Kurbatova-IM-2015:eng} some examples of extended tensor products.

(a)
Let $\alpha$ be a cross-norm on an algebraic tensor product $X\otimes Y$. We take for $X\boxtimes Y$ the completion $X\overline{\otimes}_\alpha\,Y$ of the space $X\otimes Y$ with respect to the cross-norm $\alpha$, and we take for $\mathoo B_0(X)$ and $\mathoo B_0(Y)$ the algebras $\mathoo B(X)$ and $\mathoo B(Y)$ respectively. In such a case, assumption~\eqref{e:D2} means that the cross-norm $\alpha$ is $*$-\emph{uniform}, and assumption~\eqref{e:A otimes B
le} means that the cross-norm $\alpha$ is \emph{uniform}.

(b) Let $X$ and $Y$ be Banach spaces. We denote by $\mathoo
K(X,Y)$ the Banach space of all bilinear forms $K\separ X\times Y\to\mathbb C$
that are bounded with respect to the norm $\Vert K\Vert=\sup\{\,|K(x,y)|\separ\Vert
x\Vert\le1,\,\Vert y\Vert\le1\,\}$. In order to represent $\mathoo K(X,Y)$ as an extended
tensor product $X^*\boxtimes Y^*$, we take for $\mathoo B_0(X^*)$ and $\mathoo
B_0(Y^*)$ the subalgebras of algebras $\mathoo B(X^*)$ and $\mathoo B(Y^*)$ consisting of all operators that have a preconjugate.
We define the mappings $j$, $J$, and $\mathfrak J$ by the rules (extended by linearity)
\begin{align*}
[x^*\boxtimes y^*](x,y)&=\langle x,x^*\rangle\langle y,y^*\rangle,\\
\langle x^{**}\boxtimes y^{**},K\rangle
&=\overline{K}(x^{**},y^{**}),\\
\bigl[(A\boxtimes B)K\bigr](x,y)&=K(A^0x,B^0y),
\end{align*}
where $\overline{K}$ is the canonical extension~\cite{Arens51} of $K$ to $X^{**}\times Y^{**}$.

(c) Let $X$ and $Y$ be Banach spaces, and $X\overline{\otimes}_\alpha\,Y$ be a completion of the space
$X\otimes Y$ with respect to a uniform cross-norm $\alpha$.
The conjugate space $(X\overline{\otimes}_\alpha\,Y)^*$ can be regarded as an
extended tensor product $X^*\boxtimes Y^*$ if one takes for $\mathoo B_0(X^*)$
and $\mathoo B_0(Y^*)$ the subalgebras of the algebras $\mathoo B(X^*)$ and $\mathoo B(Y^*)$ consisting of all operators that have a preconjugate. We notice that this example is a generalization of the previous one, since $\mathoo K(X,Y)\cong(X\overline{\otimes}_\pi\,Y)^*$,
where $\pi$~ is the largest cross-norm~\cite{Defant-Floret,Schatten,Helemskii:eng}.

We define $j\separ X^*\otimes Y^*\to X^*\boxtimes
Y^*=(X\overline{\otimes}_\alpha\,Y)^*$ as the canonical embedding~\eqref{e:nat imbedding j}.

Next, we define $J\separ X^{**}\otimes Y^{**}\to (X^*\boxtimes
Y^*)^{*}=(X\overline{\otimes}_\alpha\,Y)^{**}$. To this end, we assign to each functional
$w^*\in X^*\boxtimes Y^*=(X\overline{\otimes}_\alpha\,Y)^{*}$ the bilinear form
$K_{w^*}(x,y)=\langle x\otimes\,y,w^*\rangle$ on $X\times Y$.
For $\sum_{k=1}^n x_k^{**}\otimes y_k^{**}\in X^{**}\otimes Y^{**}$, we set
\begin{equation*}
\Bigl\langle J\Bigl(\sum_{k=1}^n x_k^{**}\otimes
y_k^{**}\Bigr),w^{*}\Bigr\rangle=\sum_{k=1}^n \overline{K_{w^*}}(x_k^{**},y_k^{**}),
\end{equation*}
where $\overline{K_{w^*}}$ is the canonical extension of the bilinear form $K_{w^*}$ to
$X^{**}\times Y^{**}$.

We define the operator $\mathfrak J\bigl(\sum_{k=1}^n A_k\otimes B_k\bigr)\in\mathoo
B\bigl((X\overline{\otimes}_\alpha\,Y)^*\bigr)$ as the conjugate of the operator
$\sum_{k=1}^n A_k^0\otimes B_k^0\separ X\overline{\otimes}_\alpha\,Y\to
X\overline{\otimes}_\alpha\,Y$.

(d) Let $X=L_\infty[a,b]$ and $Y=L_\infty[c,d]$. By example (c), the space
$L_\infty[a,b]\times[c,d]$ can be regarded as the extended tensor product
$L_\infty[a,b]\boxtimes L_\infty[c,d]$ (we recall that the space
$L_\infty[a,b]$ is conjugate of the space $L_1[a,b]$). We notice that one should take for $\mathoo B_0(X)$ and
$\mathoo B_0(Y)$ the subalgebras of the algebras $\mathoo B(X)$ and $\mathoo B(Y)$
consisting of all operators that have a preconjugate.

(e) Let $X$ and $Y$ be Banach spaces. We represent the space
$\mathoo B(Y,X)$ as an extended tensor product $X\boxtimes Y^*$. To this end, we take for
$\mathoo B_0(X)$ the whole algebra $\mathoo B(X)$ and we take for $\mathoo
B_0(Y^*)$ the subalgebra of the algebra $\mathoo B(Y^*)$ consisting of all operators that have a
preconjugate. We define the mappings $j$, $J$, and $\mathfrak J$ by the rules (extended by linearity)
\begin{align*}
(x\boxtimes y^*)y&=x\langle y,y^*\rangle,\\
\langle U,x^*\boxtimes\,y^{**}\rangle&=\langle y^{**},U^*x^*\rangle,\\
(A\boxtimes B)U&=AUB^0.
\end{align*}
Note that in this example the subalgebra $\mathoo B_0(Y^*)$ can be thought of as $\mathoo B(Y)$, but the action
of $\mathoo B(Y)$ on $U\in\mathoo B(Y,X)$ should be understood as contravariant, i.~e.,
$$(A_1\boxtimes B_1)\bigl((A_2\boxtimes B_2)U\bigr)=A_1A_2UB_2B_1.$$
 \end{example}

Below in this Section, we assume that we are given an extended tensor product $X\boxtimes Y$ of Banach spaces $X$
and $Y$, and a pseudo-resolvent $R_{(\cdot)}$ in the algebra $\mathoo B_0(Y)$.

\begin{theorem}[{\rm\cite[Theorem 26]{Kurbatov-Kurbatova-IM-2015:eng}}]\label{t:functional_calculas:infty in rho:R:otimes}
Assume that $\infty\notin\bar\sigma(R_{(\cdot)})$. We define the mapping\/ $\Phi\separ\mathoo O\bigl(\sigma(R_{(\cdot)}),\mathoo
B_0(X)\bigr)\to\mathoo B(X\boxtimes Y)$ by the formula
\begin{equation}\label{e:Phi(f):infty in bar rho}
\Phi(F)=\frac1{2\pi i}\int_\Gamma F(\lambda)\boxtimes R_\lambda\,d\lambda,
\end{equation}
where $\Gamma$ is an oriented envelope of the singular set
$\sigma(R_{(\cdot)})$ of the pseudo-resolvent
with respect to the point $\infty$ and the complement of the domain of $F$. We assert that $\Phi$ is a continuous morphism of unital algebras.

For all $A\in\mathoo B_0(X)$ and $h\in\mathoo O\bigl(\sigma(R_{(\cdot)})\bigr)$the morphism $\Phi$ maps the function $F(\lambda)=Ah(\lambda)$ to the operator
$A\boxtimes\varphi(h)$, where $\varphi$ is defined as in Theorem~{\rm\ref{t:functional_calculas:infty in rho}}.
\end{theorem}

We stress that the function $F$ in~\eqref{e:Phi(f):infty in bar rho} takes its values in $\mathoo
B_0(X)$, but not in $\mathbb C$.

 \begin{theorem}[{\rm\cite[Theorem 27]{Kurbatov-Kurbatova-IM-2015:eng}}]\label{t:functional_calculas:infty in sigma:R:otimes}
Assume that $\infty\in\bar\sigma(R_{(\cdot)})$. We define the mapping\/ $\Phi\separ\mathoo O\bigl(\bar\sigma(R_{(\cdot)}),\mathoo
B_0(X)\bigr)\to\mathoo B(X\boxtimes Y)$ by the formula
\begin{equation}\label{e:Phi(f):infty in bar sigma}
\Phi(F)=\frac1{2\pi i}\int_\Gamma F(\lambda)\boxtimes
R_\lambda\,d\lambda+F(\infty)\boxtimes\mathbf1,
\end{equation}
where $\Gamma$ is an oriented envelope of the extended singular set $\bar\sigma(R_{(\cdot)})$ of the
pseudo-resolvent with respect to the complement of the domain of $F$. We assert that $\Phi$ is a continuous morphism of unital algebras.

For all $A\in\mathoo B_0(X)$ and $h\in\mathoo O\bigl(\sigma(R_{(\cdot)})\bigr)$the morphism $\Phi$ maps the function $F(\lambda)=Ah(\lambda)$ to the operator
$A\boxtimes\varphi(h)$, where $\varphi$ is defined as in Theorem~{\rm\ref{t:functional_calculas:infty in sigma}}.
 \end{theorem}

 \begin{theorem}[{\rm\cite[Theorem 41]{Kurbatov-Kurbatova-IM-2015:eng}}]\label{t:1:Res}
Let $F\in\mathoo O\bigl(\bar\sigma(R_{(\cdot)}),\mathoo B_0(X)\bigr)$.
We define the operator $\Phi(F)$ by formula~\eqref{e:Phi(f):infty in bar rho} if $\infty\notin\bar\sigma(R_{(\cdot)})${\rm;} and
we define the operator $\Phi(F)$ by formula~\eqref{e:Phi(f):infty in bar sigma} if $\infty\in\bar\sigma(R_{(\cdot)})$.
We assert that the operator $\Phi(F)\separ X\boxtimes Y\to X\boxtimes Y$ is not invertible if and only if for some $\lambda\in\bar\sigma(R_{(\cdot)})$ the operator $F(\lambda)\in\mathoo B_0(X)$ is not invertible.
 \end{theorem}

 \begin{theorem}[{\rm\cite[Theorem 42]{Kurbatov-Kurbatova-IM-2015:eng}}]\label{t:spectral map:Res}
Let $F\in\mathoo O\bigl(\bar\sigma(R_{(\cdot)}),\mathoo B_0(X)\bigr)$.
We define the operator $\Phi(F)$ by formula~\eqref{e:Phi(f):infty in bar rho} if $\infty\notin\bar\sigma(R_{(\cdot)})${\rm;} and
we define the operator $\Phi(F)$ by formula~\eqref{e:Phi(f):infty in bar sigma} if $\infty\in\bar\sigma(R_{(\cdot)})$.
We assert that the spectrum of the operator $\Phi(F)\separ X\boxtimes\,Y\to X\boxtimes\,Y$ is given by the formula
\begin{equation*}
\sigma\bigl[\Phi(F)\bigr]=\bigcup_{\lambda\in\bar\sigma(R_{(\cdot)})}\sigma\bigl(F(\lambda)\bigr).
\end{equation*}
 \end{theorem}

\section{Functional calculus $\varphi_1\boxtimes\varphi_2$}\label{s:boxtimes}
In this Section, we discuss the product $\varphi_1\boxtimes\varphi_2$ of functional calculi $\varphi_1$ and $\varphi_2$ that were defined in Section~\ref{s:pseudo-resolvents}; it acts in the extended tensor product $X\boxtimes Y$. Keeping in mind the space $\mathoo B(Y,X)$ (see Example 3(e)) as the main example of an extended tensor product, we call \emph{transformators} operators acting in $X\boxtimes Y$.

In this Section, we assume that we are given an extended tensor product $X\boxtimes Y$ of Banach spaces $X$ and $Y$, and we are given pseudo-resolvents $R_{1,\,(\cdot)}$ and $R_{2,\,(\cdot)}$ in the algebras $\mathoo B_0(X)$ and $\mathoo B_0(Y)$ respectively.

 \begin{theorem}\label{t:functional_calculas:infty in rho:R_1,R_2->C}
Assume that $\infty\notin\bar\sigma(R_{1,\,(\cdot)})$
and $\infty\notin\bar\sigma(R_{2,\,(\cdot)})$. We define the mapping
$\varphi_1\boxtimes\varphi_2\:\mathoo
O\bigl(\sigma(R_{1,\,(\cdot)})\times\sigma(R_{2,\,(\cdot)})\bigr)\to\mathoo
B(X\boxtimes Y)$ by the formula
\begin{equation*}
\bigl(\varphi_1\boxtimes\varphi_2\bigr)f=\frac1{(2\pi
i)^2}\int_{\Gamma_1}\int_{\Gamma_2}f(\lambda,\mu)R_{1,\,\lambda}\boxtimes
R_{2,\,\mu}\,d\mu\,d\lambda,
\end{equation*}
where $\Gamma_1$ and $\Gamma_2$ are oriented envelopes of the singular sets
$\sigma(R_{1,\,(\cdot)})$ and $\sigma(R_{2,\,(\cdot)})$ with respect to the point $\infty$ and the complements
$\overline{\mathbb C}\setminus U_1$ and $\overline{\mathbb C}\setminus U_2${\rm;} here $U_1\times U_2$ is an open neighbourhood of the set $\sigma(R_{1,\,(\cdot)})\times\sigma(R_{2,\,(\cdot)})$ that lies in the domain of the function $f$ {\rm(}see Proposition~\ref{p:U times V}{\rm)}. We assert that $\varphi_1\boxtimes\varphi_2$ is a continuous morphism of unital algebras.

For all $g\in\mathoo O\bigl(\sigma(R_{1,\,(\cdot)})\bigr)$ and $h\in\mathoo
O\bigl(\sigma(R_{2,\,(\cdot)})\bigr)$ the morphism
$\varphi_1\boxtimes\varphi_2$ maps the function $f(\lambda,\mu)=g(\lambda)h(\mu)$ to the transformator $\varphi_1(g)\boxtimes\varphi_2(h)$, where
$\varphi_1$ and $\varphi_2$ are scalar functional calculi
{\rm(}Theorem~\ref{t:functional_calculas:infty in rho}{\rm)} generated by the pseudo-resolvents $R_{1,\,(\cdot)}$ and $R_{2,\,(\cdot)}$.
 \end{theorem}

 \begin{proof}
The proof is analogous to that of Theorem~\ref{t:functional_calculas:infty in sigma:R_1,R_2->C}, see below.
 \end{proof}

\begin{theorem}\label{t:functional_calculas:infty in rho+sigma:R_1,R_2->C}
Assume that $\infty\notin\bar\sigma(R_{1,\,(\cdot)})$, but $\infty\in\bar\sigma(R_{2,\,(\cdot)})$.
We define the mapping $\varphi_1\boxtimes\varphi_2\:\mathoo
O\bigl(\sigma(R_{1,\,(\cdot)})\times\bar\sigma(R_{2,\,(\cdot)})\bigr)\to\mathoo
B(X\boxtimes Y)$ by the formula
\begin{align*}
\bigl(\varphi_1\boxtimes\varphi_2\bigr)f=\frac1{(2\pi
i)^2}\int_{\Gamma_1}\int_{\Gamma_2}f(\lambda,\mu)R_{1,\,\lambda}\boxtimes
R_{2,\,\mu}\,d\mu\,d\lambda+\frac1{2\pi
i}\int_{\Gamma_1}f(\lambda,\infty)R_{1,\,\lambda}\boxtimes\mathbf 1_Y\,d\lambda,
\end{align*}
where $\Gamma_1$ is an oriented envelope of the singular set $\sigma(R_{1,\,(\cdot)})$ with respect to the point $\infty$ and the complement $\overline{\mathbb C}\setminus U_1$, and $\Gamma_2$ is an oriented envelope of the extended singular set $\bar\sigma(R_{2,\,(\cdot)})$ with respect to the complement $\overline{\mathbb C}\setminus U_2${\rm;} here $U_1\times U_2$ is an open neighbourhood of the set $\sigma(R_{1,\,(\cdot)})\times\bar\sigma(R_{2,\,(\cdot)})$ that lies in the domain of the function $f$ {\rm(}see Proposition~\ref{p:U times V}{\rm)}. We assert that $\varphi_1\boxtimes\varphi_2$ is a continuous morphism of unital algebras.

For all $g\in\mathoo O\bigl(\sigma(R_{1,\,(\cdot)})\bigr)$ and $h\in\mathoo O\bigl(\bar\sigma(R_{2,\,(\cdot)})\bigr)$ the morphism
$\varphi_1\boxtimes\varphi_2$ maps the function $f(\lambda,\mu)=g(\lambda)h(\mu)$ to the transformator $\varphi_1(g)\boxtimes\varphi_2(h)$, where $\varphi_1$ and $\varphi_2$ are scalar functional calculi
{\rm(}Theorems~\ref{t:functional_calculas:infty in rho} and~\ref{t:functional_calculas:infty in sigma}{\rm)} generated by the pseudo-resolvents $R_{1,\,(\cdot)}$ and $R_{2,\,(\cdot)}$.
\end{theorem}
 \begin{proof}
The proof is analogous to that of Theorem~\ref{t:functional_calculas:infty in sigma:R_1,R_2->C}, see below.
 \end{proof}

 \begin{theorem}\label{t:functional_calculas:infty in sigma:R_1,R_2->C}
Assume that $\infty\in\bar\sigma(R_{1,\,(\cdot)})$ and $\infty\in\bar\sigma(R_{2,\,(\cdot)})$. We define the mapping  $\varphi_1\boxtimes\varphi_2\:\mathoo
O\bigl(\bar\sigma(R_{1,\,(\cdot)})\times\bar\sigma(R_{2,\,(\cdot)})\bigr)\to\mathoo
B(X\boxtimes Y)$ by the formula
\begin{align*}
\bigl(\varphi_1\boxtimes\varphi_2\bigr)f&=\frac1{(2\pi
i)^2}\int_{\Gamma_1}\int_{\Gamma_2}f(\lambda,\mu)R_{1,\,\lambda}\boxtimes
R_{2,\,\mu}\,d\mu\,d\lambda+\frac1{2\pi i}\int_{\Gamma_1}f(\lambda,\infty)R_{1,\,\lambda}\boxtimes
\mathbf1\,d\lambda\\
&+\frac1{2\pi i}\int_{\Gamma_2}f(\infty,\mu)\mathbf1\boxtimes
R_{2,\,\mu}\,d\mu+f(\infty,\infty)\mathbf 1_{X\boxtimes Y},
\end{align*}
where $\Gamma_1$ and $\Gamma_2$ are oriented envelopes of the singular sets
$\bar\sigma(R_{1,\,(\cdot)})$ and $\bar\sigma(R_{2,\,(\cdot)})$ with respect to the complements
$\overline{\mathbb C}\setminus U_1$ and $\overline{\mathbb C}\setminus U_2${\rm;} here $U_1\times U_2$ is an open neighbourhood of the set $\bar\sigma(R_{1,\,(\cdot)})\times\bar\sigma(R_{2,\,(\cdot)})$ that lies in the domain of the function $f$ {\rm(}see Proposition~\ref{p:U times V}{\rm)}. We assert that $\varphi_1\boxtimes\varphi_2$ is a continuous morphism of unital algebras.

For all $g\in\mathoo O\bigl(\bar\sigma(R_{1,\,(\cdot)})\bigr)$ and $h\in\mathoo O\bigl(\bar\sigma(R_{2,\,(\cdot)})\bigr)$ the morphism $\varphi_1\boxtimes\varphi_2$ maps the function $f(\lambda,\mu)=g(\lambda)h(\mu)$ to the transformator $\varphi_1(g)\boxtimes\varphi_2(h)$, where $\varphi_1$ and $\varphi_2$ are scalar functional calculi
{\rm(}Theorem~\ref{t:functional_calculas:infty in sigma}{\rm)} generated by the
pseudo-resolvents $R_{1,\,(\cdot)}$ and $R_{2,\,(\cdot)}$.
 \end{theorem}

 \begin{proof}
For each $\mu\in U_2$ we consider the operator
\begin{equation}\label{e:g=half f:2}
G(\mu)=\varphi_2\bigl(f(\cdot,\mu)\bigr)=\frac1{2\pi i}\int_{\Gamma_1}f(\lambda,\mu)R_{1,\,\lambda}\,d\lambda+f(\infty,\mu)\mathbf1_X.
\end{equation}
By Theorem~\ref{t:functional_calculas:infty in sigma}, for any fixed
$\mu\in U_2$ the correspondence $f\mapsto G(\mu)$ preserves the three operations: addition, scalar multiplication, and multiplication. We change the interpretation: formula~\eqref{e:g=half f:2} defines a mapping $f\mapsto G$ from $\mathoo
O\bigl(\sigma(R_{1,\,(\cdot)})\times\sigma(R_{2,\,(\cdot)})\bigr)$ to $\mathoo
O\bigl(\sigma(R_{2,\,(\cdot)}),\mathoo B_0(X)\bigr)$. Since the three operations in $\mathoo
O\bigl(\sigma(R_{2,\,(\cdot)}),\mathoo B_0(X)\bigr)$ are understood in the pointwise sense, it follows that the correspondence $f\mapsto G$ is a morphism of algebras.

In accordance with Theorem~\ref{t:functional_calculas:infty in sigma:R:otimes} we put
\begin{equation}\label{e:f via g}
\begin{split}
\Phi_1(G)&=\frac1{2\pi
i}\int_{\Gamma_2}G(\mu)\boxtimes R_{2,\,\mu}\,d\mu+G(\infty)\boxtimes\mathbf1_Y\\
&=\frac1{2\pi
i}\int_{\Gamma_2}\Bigl(\frac1{2\pi i}\int_{\Gamma_1}f(\lambda,\mu)R_{1,\,\lambda}\,d\lambda+f(\infty,\mu)\mathbf1_X\Bigr)\boxtimes R_{2,\,\mu}\,d\mu\\
&+\Bigl(\frac1{2\pi i}\int_{\Gamma_1}f(\lambda,\infty)R_{1,\,\lambda}\,d\lambda+f(\infty,\infty)\mathbf1_X\Bigr)\boxtimes\mathbf1_Y\\
&=\frac1{2\pi
i}\int_{\Gamma_2}\Bigl(\frac1{2\pi i}\int_{\Gamma_1}f(\lambda,\mu)R_{1,\,\lambda}\,d\lambda\Bigr)\boxtimes R_{2,\,\mu}\,d\mu\\
&+\frac1{2\pi
i}\int_{\Gamma_2}\Bigl(f(\infty,\mu)\mathbf1_X\Bigr)\boxtimes R_{2,\,\mu}\,d\mu\\
&+\Bigl(\frac1{2\pi i}\int_{\Gamma_1}f(\lambda,\infty)R_{1,\,\lambda}\,d\lambda+f(\infty,\infty)\mathbf1_X\Bigr)\boxtimes\mathbf1_Y\\
&=\frac1{(2\pi i)^2}\int_{\Gamma_1}\int_{\Gamma_2}f(\lambda,\mu)R_{1,\,\lambda}\boxtimes R_{2,\,\mu}\,d\mu\,d\lambda+
\frac1{2\pi i}\int_{\Gamma_2}f(\infty,\mu)\mathbf1_X\boxtimes R_{2,\,\mu}\,d\mu\\
&+\frac1{2\pi i}\int_{\Gamma_1}f(\lambda,\infty)R_{1,\,\lambda}\boxtimes\mathbf1_Y\,d\lambda+f(\infty,\infty)\mathbf1_X\boxtimes\mathbf1_Y.
\end{split}
\end{equation}
By Theorem~\ref{t:functional_calculas:infty in sigma:R:otimes},
the correspondence $G\mapsto\Phi_1(G)$ also preserves the three operations. Clearly, the mapping $\varphi_1\boxtimes\varphi_2$ from the formulation of the theorem is the composition of the correspondences $f\mapsto G$ and $G\mapsto\Phi_1(G)$, and, by what has been proved, is a morphism of algebras.

The continuity is evident.

The second statement is verified by direct calculations.
 \end{proof}

When it is desirable to stress that in Theorems~\ref{t:functional_calculas:infty in rho:R_1,R_2->C}, \ref{t:functional_calculas:infty in rho+sigma:R_1,R_2->C}, and~\ref{t:functional_calculas:infty in sigma:R_1,R_2->C}, the functional calculus $\varphi_1\boxtimes\varphi_2$ is generated by pseudo-resolvents $R_{1,\,(\cdot)}$ and $R_{2,\,(\cdot)}$, we will use the notation $\varphi_{R_{1,\,(\cdot)}}\boxtimes\varphi_{R_{2,\,(\cdot)}}$ instead of $\varphi_1\boxtimes\varphi_2$, and we will use the notation $f(R_{1,\,(\cdot)},\,R_{2,\,(\cdot)})$ instead of $(\varphi_1\boxtimes\varphi_2)(f)$.
If the pseudo-resolvents $R_{1,\,(\cdot)}$ and $R_{2,\,(\cdot)}$ are generated by the operators $A$ and $B$ (see Proposition~{\rm\ref{p:1 in B_R}}), we will use the notations $\varphi_{A}\boxtimes\varphi_{B}$ and $f(A,B)$.

In order to present the definitions of $\varphi_1\boxtimes\varphi_2$ from Theorems~\ref{t:functional_calculas:infty in rho:R_1,R_2->C}, \ref{t:functional_calculas:infty in rho+sigma:R_1,R_2->C}, and~\ref{t:functional_calculas:infty in sigma:R_1,R_2->C} in a unified form, it is convenient to use the notation
\begin{equation}\label{e:delta1,delta2}
\begin{split}
 (\varphi_1\boxtimes\varphi_2)(f)&=\frac1{(2\pi
i)^2}\int_{\Gamma_1}\int_{\Gamma_2}f(\lambda,\mu)R_{1,\,\lambda}\boxtimes
R_{2,\,\mu}\,d\mu\,d\lambda\\
&+\delta_2\frac1{2\pi i}\int_{\Gamma_1}f(\lambda,\infty)R_{1,\,\lambda}\boxtimes
\mathbf1\,d\lambda\\
&+\delta_1\frac1{2\pi i}\int_{\Gamma_2}f(\infty,\mu)\mathbf1\boxtimes
R_{2,\,\mu}\,d\mu+\delta_1\delta_2f(\infty,\infty)\mathbf1\boxtimes\mathbf1,
 \end{split}
\end{equation}
where $\delta_i=1$ if $\Gamma_i$ encloses $\infty$, and $\delta_i=0$ in the opposite case, $i=1,2$.

We enumerate the results of the action of $\varphi_1\boxtimes\varphi_2$ on some frequently encountered functions.

 \begin{corollary}\label{c:functional_calculas:infty in rho:R_1,R_2->C}
Under the assumptions of Theorems~\ref{t:functional_calculas:infty in rho:R_1,R_2->C}, \ref{t:functional_calculas:infty in rho+sigma:R_1,R_2->C}, and \ref{t:functional_calculas:infty in sigma:R_1,R_2->C}, the morphism $\varphi_1\boxtimes\varphi_2$ maps the function $u(\lambda,\mu)=1$ to the unit $\mathbf1\boxtimes\mathbf1$ of the algebra
$\mathoo B(X\boxtimes Y)${\rm;} the function
$r_{1,\lambda_0}(\lambda,\mu)=\frac1{\lambda_0-\lambda}$, where
$\lambda_0\in\rho(R_{1,\,(\cdot)})$, is mapped by the morphism $\varphi_1\boxtimes\varphi_2$ to the transformator $R_{1,\,\lambda_0}\boxtimes\mathbf 1_Y${\rm;} the function
$r_{2,\mu_0}(\lambda,\mu)=\frac1{\mu_0-\mu}$, where $\mu_0\in\rho(R_{2,\,(\cdot)})$,
is mapped by the morphism $\varphi_1\boxtimes\varphi_2$ to the transformator $\mathbf 1_X\boxtimes
R_{2,\,\mu_0}${\rm;} the function
$r_{\lambda_0,\mu_0}(\lambda,\mu)=\frac1{(\lambda_0-\lambda)(\mu_0-\mu)}$, where
$\lambda_0\in\rho(R_{1,\,(\cdot)})$ and $\mu_0\in\rho(R_{2,\,(\cdot)})$, is mapped by the morphism $\varphi_1\boxtimes\varphi_2$ to the transformator $R_{1,\,\lambda_0}\boxtimes
R_{2,\,\mu_0}$.

Under the assumptions of Theorems~\ref{t:functional_calculas:infty in rho:R_1,R_2->C} and \ref{t:functional_calculas:infty in rho+sigma:R_1,R_2->C}, the morphism $\varphi_1\boxtimes\varphi_2$ maps the function
$c_1(\lambda,\mu)=\lambda$ to the transformator $A\boxtimes\mathbf 1_Y$, where $A$ is the operator that generates
the maximal pseudo-resolvent $R_{1,\,(\cdot)}$ in accordance with Proposition~{\rm\ref{p:1 in B_R}}.

Under the assumptions of Theorem~\ref{t:functional_calculas:infty in rho:R_1,R_2->C}, the morphism $\varphi_1\boxtimes\varphi_2$ maps the function $c_2(\lambda,\mu)=\mu$ to the transformator $\mathbf 1_X\boxtimes B$, where $B$ is the operator that generates
the maximal pseudo-resolvent $R_{1,\,(\cdot)}$ in accordance with Proposition~{\rm\ref{p:1 in B_R}}{\rm;} the function $r_{\nu_0}(\lambda,\mu)=\frac1{\nu_0-\lambda\mp\mu}$ is mapped by the morphism
$\varphi_1\boxtimes\varphi_2$ to the transformator
$(\nu_0\mathbf1\boxtimes\mathbf1-A\boxtimes\mathbf1\mp\mathbf1\boxtimes B)^{-1}$ provided $\nu_0\notin\sigma(A)\pm\sigma(B)$.
 \end{corollary}
 \begin{proof}
We restrict ourselves to proving the last statement. Clearly, the function
$(\lambda,\mu)\mapsto\nu_0-\lambda\mp\mu$ is mapped by the morphism
$\varphi_1\boxtimes\varphi_2$ to the transformator
$\nu_0\mathbf1\boxtimes\mathbf1-A\boxtimes\mathbf1\mp\mathbf1\boxtimes B$. Since
$\varphi_1\boxtimes\varphi_2$ is a morphism of algebras, the reciprocal function
is mapped to the inverse transformator.
 \end{proof}

\begin{theorem}\label{t:sp map of comp}
Let $g\in\mathoo
O\bigl(\bar\sigma(R_{1,\,(\cdot)})\times\bar\sigma(R_{2,\,(\cdot)})\bigr)$ and $f\in\mathoo  O\bigl(g\bigl(\bar\sigma(R_{1,\,(\cdot)})\times\bar\sigma(R_{2,\,(\cdot)})\bigr)\bigr)$. Then the transformator
$(\varphi_1\boxtimes\varphi_2)(f\circ g)$
is the function $f$ of the transformator $(\varphi_1\boxtimes\varphi_2)(g)${\rm:}
\begin{equation*}
(\varphi_1\boxtimes\varphi_2)(f\circ g)=\frac1{2\pi i}\int_{\Gamma_3}f(\nu)\bigl(\nu\mathbf1\boxtimes\mathbf1-(\varphi_1\boxtimes\varphi_2)(g)\bigr)^{-1}\,d\nu,
\end{equation*}
where $\Gamma_3$ is an oriented envelope of the spectrum $\sigma\bigl((\varphi_1\boxtimes\varphi_2)(g)\bigr)$.
\end{theorem}

 \begin{proof} We have ($\delta_1,\delta_2=0,1$)
\begin{align*}
(\varphi_1\boxtimes\varphi_2)&(f\circ g)=\frac1{(2\pi
i)^2}\int_{\Gamma_1}\int_{\Gamma_2}f\bigl(g(\lambda,\mu)\bigr)R_{1,\,\lambda}\boxtimes
R_{2,\,\mu}\,d\mu\,d\lambda\\
&+\delta_2\frac1{2\pi i}\int_{\Gamma_1}f\bigl(g(\lambda,\infty)\bigr)R_{1,\,\lambda}\boxtimes
\mathbf1\,d\lambda\\
&+\delta_1\frac1{2\pi i}\int_{\Gamma_2}f\bigl(g(\infty,\mu)\bigr)\mathbf1\boxtimes
R_{2,\,\mu}\,d\mu+\delta_1\delta_2f\bigl(g(\infty,\infty)\bigr)\mathbf 1_{X\boxtimes Y}\\
&=\frac1{(2\pi
i)^2}\int_{\Gamma_1}\int_{\Gamma_2}\biggl[\frac1{2\pi i}\int_{\Gamma_3}\frac{f(\nu)}{\nu-g(\lambda,\mu)}\,d\nu\biggr]R_{1,\,\lambda}\boxtimes
R_{2,\,\mu}\,d\mu\,d\lambda\\
&+\delta_2\frac1{2\pi i}\int_{\Gamma_1}\biggl[\frac1{2\pi i}\int_{\Gamma_3}\frac{f(\nu)}{\nu-g(\lambda,\infty)}\,d\nu\biggr]R_{1,\,\lambda}\boxtimes
\mathbf1\,d\lambda\\
&+\delta_1\frac1{2\pi i}\int_{\Gamma_2}\biggl[\frac1{2\pi i}\int_{\Gamma_3}\frac{f(\nu)}{\nu-g(\infty,\mu)}\,d\nu\biggr]\mathbf1\boxtimes
R_{2,\,\mu}\,d\mu+\delta_1\delta_2f\bigl(g(\infty,\infty)\bigr)\mathbf 1_{X\boxtimes Y}\\
\intertext{(here we interchange the order of integration)}
&=\frac1{2\pi i}\int_{\Gamma_3}f(\nu)\biggl[\frac1{(2\pi
i)^2}\int_{\Gamma_1}\int_{\Gamma_2}\frac{R_{1,\,\lambda}\boxtimes
R_{2,\,\mu}}{\nu-g(\lambda,\mu)}\,d\mu\,d\lambda\biggr]\,d\nu\\
&+\delta_2\frac1{2\pi i}\int_{\Gamma_3}f(\nu)\biggl[\frac1{2\pi i}\int_{\Gamma_1}\frac{R_{1,\,\lambda}\boxtimes
\mathbf1}{\nu-g(\lambda,\infty)}\,d\lambda\biggr]\,d\nu\\
&+\delta_1\frac1{2\pi i}\int_{\Gamma_3}f(\nu)\biggl[\frac1{2\pi i}\int_{\Gamma_2}\frac{\mathbf1\boxtimes
R_{2,\,\mu}}{\nu-g(\infty,\mu)}\,d\mu\biggr]\,d\nu+\delta_1\delta_2f\bigl(g(\infty,\infty)\bigr)\mathbf 1_{X\boxtimes Y}\\
\intertext{(further, by Theorems~\ref{t:functional_calculas:infty in rho:R_1,R_2->C}, \ref{t:functional_calculas:infty in rho+sigma:R_1,R_2->C}, and \ref{t:functional_calculas:infty in sigma:R_1,R_2->C}, it follows that)}
&=\frac1{2\pi i}\int_{\Gamma_3}f(\nu)\biggl[\bigl(\nu\mathbf1\boxtimes\mathbf1-(\varphi_1\boxtimes\varphi_2)(g)\bigr)^{-1}-\delta_2\frac1{2\pi i}\int_{\Gamma_1}\frac{R_{1,\,\lambda}\boxtimes\mathbf1}{\nu-g(\lambda,\infty)}\,d\lambda\\
&-\delta_1\frac1{2\pi i}\int_{\Gamma_2}\frac{\mathbf1\boxtimes R_{2,\,\mu}}{\nu-g(\infty,\mu)}\,d\mu-\delta_1\delta_2g(\infty,\infty)\mathbf 1_{X\boxtimes Y}\biggr]\,d\nu\\
&+\delta_2\frac1{2\pi i}\int_{\Gamma_3}f(\nu)\biggl[\frac1{2\pi i}\int_{\Gamma_1}\frac{R_{1,\,\lambda}\boxtimes
\mathbf1}{\nu-g(\lambda,\infty)}\,d\lambda\biggr]\,d\nu\\
&+\delta_1\frac1{2\pi i}\int_{\Gamma_3}f(\nu)\biggl[\frac1{2\pi i}\int_{\Gamma_2}\frac{\mathbf1\boxtimes
R_{2,\,\mu}}{\nu-g(\infty,\mu)}\,d\mu\biggr]\,d\nu+\delta_1\delta_2f\bigl(g(\infty,\infty)\bigr)\mathbf 1_{X\boxtimes Y}\\
&=\frac1{2\pi i}\int_{\Gamma_3}f(\nu)\bigl(\nu\mathbf1\boxtimes\mathbf1-(\varphi_1\boxtimes\varphi_2)(g)\bigr)^{-1}\,d\nu.\qed
\end{align*}
\renewcommand\qed{}
 \end{proof}

 \begin{corollary}\label{c:functional_calculas:infty in rho:R_1,R_2->C:lambda-mu}
Let $A\in\mathoo B(X)$ and $B\in\mathoo B(Y)$. Let $f\in\mathoo{O}\bigl(\sigma(A)\pm\sigma(B)\bigr)$. Then
\begin{equation*}
\frac1{(2\pi i)^2}\int_{\Gamma_1}\int_{\Gamma_2}f(\lambda\pm\mu)R_{A,\,\lambda}\boxtimes
R_{B,\,\mu}\,d\mu\,d\lambda
=\frac1{2\pi i}\int_{\Gamma_3}f(\nu)(\nu\mathbf1\boxtimes\mathbf1-A\boxtimes\mathbf1\mp\mathbf1\boxtimes B)^{-1}\,d\nu,
\end{equation*}
where $\Gamma_3$ is an oriented envelope of $\sigma(A)\pm\sigma(B)$.
 \end{corollary}

 \begin{proof}
This is a special case of Theorem~\ref{t:sp map of comp} for $g(\lambda,\mu)=\lambda\pm\mu$.
 \end{proof}

 \begin{example}\label{ex:tensor sum}
Let $A\in\mathoo{B}(X)$ and $B\in\mathoo{B}(Y)$.
By Corollary~\ref{c:functional_calculas:infty in rho:R_1,R_2->C:lambda-mu} and the formula $e^{\lambda t}e^{\mu t}=e^{(\lambda+\mu)t}$, one has (cf.~\cite{Benzi-Simoncini:1501.07376,Graham}, \cite[Theorem 10.9]{Higham08})
\begin{equation*}
 e^{At}\boxtimes e^{Bt}=e^{(A\boxtimes\mathbf1+\mathbf1\boxtimes B)t}.
\end{equation*}
 \end{example}

We proceed to the discussion of spectral mapping theorems.

 \begin{theorem}\label{t:1:Res:R_1,R_2->C}
Let $f\in\mathoo O\bigl(\bar\sigma(R_{1,\,(\cdot)})\times\bar\sigma(R_{2,\,(\cdot)})\bigr)$.
Then the transformator $\bigl(\varphi_1\boxtimes\varphi_2\bigr)(f)\separ X\boxtimes Y\to
X\boxtimes Y$ is not invertible if and only if $f(\lambda,\mu)=0$ for at least one couple of points $\lambda\in\bar\sigma(R_{1,\,(\cdot)})$ and $\mu\in\bar\sigma(R_{2,\,(\cdot)})$.
 \end{theorem}
 \begin{proof}
For each $\mu\in U_2$, we consider operator~\eqref{e:g=half f:2}. By Theorem~\ref{t:spectral mappping}, the following statement holds: the operator $G(\mu)\separ X\to X$  is not invertible if and only if $f(\lambda,\mu)=0$ for at least one
$\lambda\in\bar\sigma(R_{1,\,(\cdot)})$. Further, by Theorem~\ref{t:1:Res}, operator~\eqref{e:f via g} is not invertible if and only if $G(\mu)$ is not invertible for at least one $\mu\in\bar\sigma(R_{2,\,(\cdot)})$.
Combining (in the opposite order) all these results, we arrive at the desired statement.
 \end{proof}

 \begin{theorem}\label{t:spectral map:R_1,R_2->C}
Let $f\in\mathoo O\bigl(\bar\sigma(R_{1,\,(\cdot)})\times\bar\sigma(R_{2,\,(\cdot)})\bigr)$.
Then the spectrum of the transformator $(\varphi_1\boxtimes\varphi_2)f\separ X\boxtimes Y\to X\boxtimes Y$
is given by the formula
\begin{equation*}
\sigma\bigl((\varphi_1\boxtimes\varphi_2)f\bigr)=
\{\,f(\lambda,\mu):\,\lambda\in\bar\sigma(R_{1,\,(\cdot)}),\,\mu\in\bar\sigma(R_{2,\,(\cdot)})\,\}.
\end{equation*}
 \end{theorem}

 \begin{proof}
We take an arbitrary $\nu\in\mathbb C$. By the definition of the spectrum, the number $\nu$ belongs to the set
$\sigma\bigl((\varphi_1\boxtimes\varphi_2)f\bigr)$ if and only if the transformator
$\nu\mathbf1_{X\boxtimes Y}-\bigl(\varphi_1\boxtimes\varphi_2\bigr)(f)$
is not invertible.

We denote by $u$ the function from $\mathoo
O\bigl(\bar\sigma(R_{1,\,(\cdot)})\times\bar\sigma(R_{2,\,(\cdot)})\bigr)$
that identically equals $1$. By Theorems~\ref{t:functional_calculas:infty in rho:R_1,R_2->C},
\ref{t:functional_calculas:infty in rho+sigma:R_1,R_2->C}, and \ref{t:functional_calculas:infty in sigma:R_1,R_2->C}, we have
\begin{equation*}
\bigl(\varphi_1\boxtimes\varphi_2\bigr)(u)=\mathbf 1_{X\boxtimes Y},
\end{equation*}
whence
\begin{equation*}
\nu\mathbf1_{X\boxtimes Y}-\bigl(\varphi_1\boxtimes\varphi_2\bigr)(f)=\bigl(\varphi_1\boxtimes\varphi_2\bigr)(\nu
u-f).
\end{equation*}
We apply Theorem~\ref{t:1:Res:R_1,R_2->C}: the transformator
$\bigl(\varphi_1\boxtimes\varphi_2\bigr)(\nu u-f)$ is not invertible if and only if $\nu u(\lambda,\mu)-f(\lambda,\mu)=0$ for some
$\lambda\in\bar\sigma(R_{1,\,(\cdot)})$ and $\mu\in\bar\sigma(R_{2,\,(\cdot)})$ or, in other words,
$\nu\in\{\,f(\lambda,\mu)\separ\lambda\in\bar\sigma(R_{1,\,(\cdot)}),\,\mu\in\bar\sigma(R_{2,\,(\cdot)})\,\}$.
 \end{proof}

We denote by $\mathoo B_{R_{1},\,R_{2}}$ the closure in $\mathoo B\bigl(\mathoo B(X,Y)\bigr)$ of the set of all transformators $(\varphi_1\boxtimes\varphi_2)f$, where $f\in\mathoo O\bigl(\bar\sigma(R_{1,\,(\cdot)})\times\bar\sigma(R_{2,\,(\cdot)})\bigr)$.

\begin{corollary}\label{c:f(A,B) is full}
The set $\mathoo B_{R_{1},\,R_2}$ is a full commutative subalgebra of the algebra $\mathoo B\bigl(\mathoo B(X,Y)\bigr)$ of all transformators acting in $\mathoo B(X,Y)$.
\end{corollary}
\begin{proof}
Clearly, the image under $\varphi_1\boxtimes\varphi_2$ of the unital commutative algebra $\mathoo O\bigl(\bar\sigma(R_{1,\,(\cdot)})\times\bar\sigma(R_{2,\,(\cdot)})\bigr)$ is a unital commutative subalgebra.

Let the transformator $(\varphi_1\boxtimes\varphi_2)f$ be invertible. By Theorem~\ref{t:spectral map:R_1,R_2->C}, this means that $f(\lambda,\mu)\neq0$ for some $\lambda\in\bar\sigma(R_{1,\,(\cdot)})$ and $\mu\in\bar\sigma(R_{2,\,(\cdot)})$. Clearly, the inverse of $(\varphi_1\boxtimes\varphi_2)f$ is the transformator $(\varphi_1\boxtimes\varphi_2)\frac1f$.

It remains to apply Proposition~\ref{p:closure of full}.
\end{proof}

In~\cite{Schechter69}, an analogue of Theorem~\ref{t:functional_calculas:infty in rho:R_1,R_2->C} was proved in the tensor product of Banach spaces for bounded operators and a polynomial function $f$. In~\cite[Theorem 2.4]{Ichinose75}, an analogue of Theorem~\ref{t:functional_calculas:infty in rho:R_1,R_2->C} was proved in the tensor product of Banach spaces for bounded operators and an arbitrary analytic function $f${\rm;} in~\cite[Theorem 2.4]{Ichinose75}, an analogue of Theorem~\ref{t:functional_calculas:infty in sigma:R_1,R_2->C} was also proved
for unbounded operators and analytic functions. See also the initial version~\cite{Ichinose70} of article~\cite{Ichinose75}.

An analogue of Theorem~\ref{t:sp map of comp} for matrices was proved in~\cite[Theorem 4.4]{Kressner:10.7153/oam-08-23}.

There are several versions of Theorem~\ref{t:spectral map:R_1,R_2->C} in tensor products of Banach spaces.
It was shown in~\cite{Brown-Pearcy} that the spectrum of the tensor product $A\otimes B$ of two bounded operators acting in a Hilbert space is the set $\sigma(A)\times\sigma(B)$. For functions $f$ of the form $f(\lambda,\mu)=g(\lambda)h(\mu)$, Theorem~\ref{t:spectral map:R_1,R_2->C} was proved in~\cite{Lumer-Rosenblum}{\rm;} for polynomial functions $f$ of two variables, a version of Theorem~\ref{t:spectral map:R_1,R_2->C} was proved in~\cite[Theorem 3.3]{Harte73}, see also~\cite{Harte14}{\rm;} another equivalent version was proved in~\cite[Theorem 3.4]{Embry-Rosenblum}. In~\cite[Theorem 3.2]{Ichinose75}, it was proved an analogue of Theorem~\ref{t:spectral map:R_1,R_2->C} for unbounded operators and analytic functions, see also~\cite{Ichinose70}.
A modern version of Theorem~\ref{t:spectral map:R_1,R_2->C} for matrices can be found in~\cite[Lemma 4.1]{Kressner:10.7153/oam-08-23}.

A functional calculus for the transformator $A\otimes\mathbf1-\mathbf1\otimes B$ was first described in~\cite{Rosenblum56:Duke}.
Functions of the transformator $A\otimes\mathbf1\pm\mathbf1\otimes B$ are also investigated in~\cite{Benzi-Simoncini:1501.07376, Benzi-Simoncini:1503.02615,Gil'2006:IEOT}.

\section{Meromorphic functional calculus}\label{s:mero}
A meromorphic function of a bounded operator is an unbounded operator or a linear relation (provided a pole of the function is contained in the spectrum). According to our approach, we identify such an object with its resolvent.

Let $U$ be an open subset of $\overline{\mathbb C}^2$ and $f:\,U\to\overline{\mathbb C}$.
The function $f$ is called~\cite[ch.~IV, \S~15.43]{Shabat2:eng}
\emph{meromorphic} if: (i) $f$ is analytic on a set $U\setminus M$, where $M$ is a nowhere dense closed subset of $U$, (ii) $f$ cannot be analytically continued to any point of $M$, (iii) for any point $\zeta\in M$ there exist a connected neighborhood $V$ of $\zeta$ and an analytic function $q_\zeta:\,V\to\mathbb C$ such that the function $p_\zeta=f\cdot q_\zeta$ is analytic in $V\cap(U\setminus M)$ and can be extended analytically into $V$, and $q_\zeta$ equals zero only on $V\cap M$. Clearly, $q_\zeta(\zeta)=0$.
The set $M$ is called the \emph{polar set of the function~$f$}. It consists of points of two types: if $p_\zeta(\zeta)\neq0$ (and so $\lim_{z\to\zeta}f(z)=\infty$), then $\zeta$ is called a \emph{pole}; if $p_\zeta(\zeta)=0$, then $\zeta$ is called a \emph{point of indeterminacy}. In any neighbourhood of a point of indeterminacy, the function $f$ takes any value from $\mathbb C$~\cite{Shabat2:eng}. For example, for the function $f(\lambda,\mu)=\lambda\mu$, the points of indeterminacy are $(0,\infty)$ and $(\infty,0)$, for the function $f(\lambda,\mu)=\frac\lambda\mu$, the points of indeterminacy are $(0,0)$ and $(\infty,\infty)$, and for the function $f(\lambda,\mu)=\lambda-\mu$, the point of indeterminacy is $(\infty,\infty)$.

Assume that we are given an extended tensor product $X\boxtimes Y$ of Banach spaces $X$ and $Y$, and we are given pseudo-resolvents
$R_{1,\,(\cdot)}$ and $R_{2,\,(\cdot)}$ in the algebras $\mathoo B_0(X)$ and $\mathoo B_0(Y)$ respectively.

We consider a function $f$ that is meromorphic in a neighbourhood $U\subseteq\overline{\mathbb C}^2$ of the set $\bar\sigma(R_{1,\,(\cdot)})\times\bar\sigma(R_{2,\,(\cdot)})$ and has no points of indeterminacy in $U$. We consider the subset
\begin{equation*}
f\bigl(\bar\sigma(R_{1,\,(\cdot)}),\bar\sigma(R_{2,\,(\cdot)})\bigr)
=\{\,f(\lambda,\mu):\,(\lambda,\mu)\in\bar\sigma(R_{1,\,(\cdot)})\times\bar\sigma(R_{2,\,(\cdot)})\,\}
\end{equation*}
of the set $\overline{\mathbb C}$. The set $f\bigl(\bar\sigma(R_{1,\,(\cdot)}),\bar\sigma(R_{2,\,(\cdot)})\bigr)$ is compact, being the image under the continuous function $f$ of the compact set $\bar\sigma(R_{1,\,(\cdot)})\times\bar\sigma(R_{2,\,(\cdot)})$.

\begin{lemma}\label{l:cup:OFC:2}
For any $\nu\notin f\bigl(\bar\sigma(R_{1,\,(\cdot)}),\bar\sigma(R_{2,\,(\cdot)})\bigr)$ the set
\begin{equation}\label{e:union:OFC:2}
(\overline{\mathbb C}^2\setminus U)\cup\{\,(\lambda,\mu)\in U:\,f(\lambda,\mu)=\nu\,\}
\end{equation}
is closed in $\overline{\mathbb C}^2$ and does not intersect $\bar\sigma(R_{1,\,(\cdot)})\times\bar\sigma(R_{2,\,(\cdot)})$.
Moreover, for any closed set $W\subseteq\overline{\mathbb C}\setminus f\bigl(\bar\sigma(R_{1,\,(\cdot)}),\bar\sigma(R_{2,\,(\cdot)})\bigr)$, the set
\begin{equation}\label{e:union:OFM:2}
(\overline{\mathbb C}^2\setminus U)\cup\{\,(\lambda,\mu)\in U:\,f(\lambda,\mu)\in W\,\}
\end{equation}
is closed in $\overline{\mathbb C}^2$ and does not intersect $\bar\sigma(R_{1,\,(\cdot)})\times\bar\sigma(R_{2,\,(\cdot)})$.
\end{lemma}

\begin{proof}
The set $\overline{\mathbb C}^2\setminus U$ is closed, being a complement of an open one. The set $\{\,(\lambda,\mu)\in U:f(\lambda,\mu)\in W\,\}=f^{-1}(W)$ is closed in $U$, being the inverse image of the closed set $W$ under the continuous function $f$. This means that limit points of the set $f^{-1}(W)=\{\,(\lambda,\mu)\in U:\,f(\lambda,\mu)\in W\,\}$ either belongs to $f^{-1}(W)$ or to the complement of $\overline{\mathbb C}^2\setminus U$. Thus, set~\eqref{e:union:OFM:2} is closed.

We show that set~\eqref{e:union:OFM:2} is disjoint from $\bar\sigma(R_{1,\,(\cdot)})\times\bar\sigma(R_{2,\,(\cdot)})$. Actually, if $f(\lambda,\mu)=\nu\in W$ and $(\lambda,\mu)\in\bar\sigma(R_{1,\,(\cdot)})\times\bar\sigma(R_{2,\,(\cdot)})$, then $\nu\in f\bigl(\bar\sigma(R_{1,\,(\cdot)}),\bar\sigma(R_{2,\,(\cdot)})\bigr)$, which contradicts the assumption. If $(\lambda,\mu)\notin U$, then $(\lambda,\mu)\notin\bar\sigma(R_{1,\,(\cdot)})\times\bar\sigma(R_{2,\,(\cdot)})$ by the definition of $U$.
\end{proof}

For all $\nu\in\mathbb C\setminus f\bigl(\bar\sigma(R_{1,\,(\cdot)}),\bar\sigma(R_{2,\,(\cdot)})\bigr)$, we set
\begin{equation}\label{e:def of S_nu:2}
\begin{split}
S_\nu&=\frac1{(2\pi
i)^2}\int_{\Gamma_1}\int_{\Gamma_2}\frac{1}{\nu-f(\lambda,\mu)}R_{1,\,\lambda}\boxtimes
R_{2,\,\mu}\,d\mu\,d\lambda\\
&+\delta_2\frac1{2\pi i}\int_{\Gamma_1}\frac{1}{\nu-f(\lambda,\infty)}R_{1,\,\lambda}\boxtimes
\mathbf1\,d\lambda\\
&+\delta_1\frac1{2\pi i}\int_{\Gamma_2}\frac{1}{\nu-f(\infty,\mu)}\mathbf1\boxtimes
R_{2,\,\mu}\,d\mu+\frac{\delta_1\delta_2}{\nu-f(\infty,\infty)}\mathbf1\boxtimes\mathbf1,
\end{split}
\end{equation}
where $\Gamma_i$ is an oriented envelope of the spectrum $\sigma(R_{i,(\cdot)})$; $\delta_i=1$ if $\Gamma_i$ encloses $\infty$, and $\delta_i=0$ in the opposite case; $i=1,2$.
By Lemma~\ref{l:cup:OFC:2}, the function $h_\nu(\lambda,\mu)=\frac{1}{\nu-f(\lambda,\mu)}$ belongs to $\mathoo O\bigl(\bar\sigma(R_{1,\,(\cdot)})\times\bar\sigma(R_{2,\,(\cdot)})\bigr)$. Therefore $S_\nu$ can be regarded as the image under the morphism $\varphi_1\boxtimes\varphi_2$ of the function $h_\nu$:
\begin{equation*}
S_\nu=(\varphi_1\boxtimes\varphi_2)h_\nu.
\end{equation*}

We denote by $S_\nu\bigl(R_{1,\,(\cdot)},R_{2,\,(\cdot)}\bigr)$ transformator~\eqref{e:def of S_nu:2} generated by the pseudo-resolvents $R_{1,\,(\cdot)}$ and $R_{2,\,(\cdot)}$, and we call
$S_\nu$ the \emph{resolvent of the function $f$ of $R_{1,\,\lambda}$ and $R_{2,\,\lambda}$}.

 \begin{theorem}\label{t:mero f:2}
Let the function $f$ be meromorphic in an open neighbourhood $U\subseteq\overline{\mathbb C}^2$ of the set $\bar\sigma(R_{1,\,(\cdot)})\times\bar\sigma(R_{2,\,(\cdot)})$ and have no points of indeterminacy in $U$.
Then the family
\begin{equation*}
S_\nu,\qquad\nu\notin f\bigl(\bar\sigma(R_{1,\,(\cdot)}),\bar\sigma(R_{2,\,(\cdot)})\bigr),
\end{equation*}
defined by formula~\eqref{e:def of S_nu:2} is a maximal pseudo-resolvent. In particular,
\begin{equation}\label{s:spr of S nu:2}
\bar\sigma(S_{(\cdot)})=f\bigl(\bar\sigma(R_{1,\,(\cdot)}),\bar\sigma(R_{2,\,(\cdot)})\bigr).
\end{equation}
 \end{theorem}

Equality~\eqref{s:spr of S nu:2} can be considered as an analogue of the spectral mapping theorem.

\begin{proof}
We show that, on the set $\mathbb C\setminus f\bigl(\bar\sigma(R_{1,\,(\cdot)}),\bar\sigma(R_{2,\,(\cdot)})\bigr)$, the Hilbert identity holds:
\begin{equation}\label{e:Hilbert:S:2}
S_{\nu_1}-S_{\nu_2}=-(\nu_1-\nu_2) S_{\nu_1} S_{\nu_2},\qquad\nu_1,\nu_2\notin f\bigl(\bar\sigma(R_{1,\,(\cdot)}),\bar\sigma(R_{2,\,(\cdot)})\bigr).
\end{equation}
We note that
\begin{equation*}
\frac{1}{\nu_1-f(\lambda,\mu)}-\frac{1}{\nu_2-f(\lambda,\mu)}=-\frac{\nu_1-\nu_2}{\bigl(\nu_1-f(\lambda,\mu)\bigr){\bigl(\nu_2-f(\lambda,\mu)\bigr)}}.
\end{equation*}
Applying the morphism $\varphi_1\boxtimes\varphi_2$ to this identity we arrive at the Hilbert identity~\eqref{e:Hilbert:S:2}.

We verify that the pseudo-resolvent $S_{(\cdot)}$ is maximal.
The validity of the Hilbert identity implies that
\begin{equation*}
\sigma(S_{(\cdot)})\subseteq f\bigl(\bar\sigma(R_{1,\,(\cdot)}),\bar\sigma(R_{2,\,(\cdot)})\bigr).
\end{equation*}
To prove the reverse inclusion, we fix an auxiliary point $\nu\in\mathbb C\setminus f\bigl(\bar\sigma(R_{1,\,(\cdot)}),\bar\sigma(R_{2,\,(\cdot)})\bigr)$. By Theorem~\ref{t:pseudoresovent}, the pseudo-resolvent $S_{(\cdot)}$ can be extended to points $\eta\in\mathbb C$ in which the transformator $\mathbf1+(\eta-\nu)S_\nu$ is invertible. By Theorems~\ref{t:functional_calculas:infty in rho:R_1,R_2->C},
\ref{t:functional_calculas:infty in rho+sigma:R_1,R_2->C},~\ref{t:functional_calculas:infty in sigma:R_1,R_2->C}, and \ref{t:spectral map:R_1,R_2->C}, and formula~\eqref{e:def of S_nu:2} we have
\begin{equation*}
\sigma(S_\nu)
=\Bigl\{\frac1{\nu-f(\lambda,\mu)}:\,(\lambda,\mu)\in\bar\sigma(R_{1,\,(\cdot)})\times\bar\sigma(R_{2,\,(\cdot)})\Bigr\}.
\end{equation*}
It follows that
\begin{align*}
\sigma(\mathbf1+(\eta-\nu)S_\nu)&=\Bigl\{1+\frac{\eta-\nu}{\nu-f(\lambda,\mu)}:\,(\lambda,\mu)\in\bar\sigma(R_{1,\,(\cdot)})\times\bar\sigma(R_{2,\,(\cdot)})\Bigr\}\\
&=\Bigl\{\frac{\eta-f(\lambda,\mu)}{\nu-f(\lambda,\mu)}:\,(\lambda,\mu)\in\bar\sigma(R_{1,\,(\cdot)})\times\bar\sigma(R_{2,\,(\cdot)})\Bigr\}.
\end{align*}
From this formula it is seen that the transformator $\mathbf1+(\eta-\nu)S_\nu$ is invertible if and only if
\begin{equation*}
0\notin\Bigl\{\frac{\eta-f(\lambda,\mu)}{\nu-f(\lambda,\mu)}:\,(\lambda,\mu)\in\bar\sigma(R_{1,\,(\cdot)})\times\bar\sigma(R_{2,\,(\cdot)})\Bigr\},
\end{equation*}
which is equivalent to
\begin{equation*}
\eta\notin\bigl\{f(\lambda,\mu):\,(\lambda,\mu)\in\bar\sigma(R_{1,\,(\cdot)})\times\bar\sigma(R_{2,\,(\cdot)})\bigr\}.
\end{equation*}
Thus, from $\eta\notin f\bigl(\bar\sigma(R_{1,\,(\cdot)}),\bar\sigma(R_{2,\,(\cdot)})\bigr)$ it follows that $\eta\notin\sigma(S_{(\cdot)})$.

It remains to analyze the case $\nu=\infty$.

We assume that $\infty\notin f\bigl(\bar\sigma(R_{1,\,(\cdot)}),\bar\sigma(R_{2,\,(\cdot)})\bigr)$. Then, since the set $f\bigl(\bar\sigma(R_{1,\,(\cdot)}),\bar\sigma(R_{2,\,(\cdot)})\bigr)$ is closed, a neighbourhood $W\subseteq\mathbb C$ of infinity is also disjoint from $f\bigl(\bar\sigma(R_{1,\,(\cdot)}),\bar\sigma(R_{2,\,(\cdot)})\bigr)$. Without loss of generality, we may assume that the neighbourhood $W$ is closed. By definition, $S_\nu$ is defined for all $\nu\in W\setminus\{\infty\}$; besides, by Lemma~\ref{l:cup:OFC:2}, we may assume that the contours $\Gamma_1$ and $\Gamma_2$ in~\eqref{e:def of S_nu:2} do not depend on $\nu\in W\setminus\{\infty\}$. We calculate the limit (see Theorems~\ref{t:functional_calculas:infty in rho:R_1,R_2->C},
\ref{t:functional_calculas:infty in rho+sigma:R_1,R_2->C}, and \ref{t:functional_calculas:infty in sigma:R_1,R_2->C}):
\begin{align*}
\lim_{\nu\to\infty}\nu S_\nu&=\lim_{\nu\to\infty}\biggl(\frac1{(2\pi
i)^2}\int_{\Gamma_1}\int_{\Gamma_2}\frac{\nu}{\nu-f(\lambda,\mu)}R_{1,\,\lambda}\boxtimes
R_{2,\,\mu}\,d\mu\,d\lambda\\
&+\delta_2\frac1{2\pi i}\int_{\Gamma_1}\frac{\nu}{\nu-f(\lambda,\infty)}R_{1,\,\lambda}\boxtimes
\mathbf1\,d\lambda\\
&+\delta_1\frac1{2\pi i}\int_{\Gamma_2}\frac{\nu}{\nu-f(\infty,\mu)}\mathbf1\boxtimes
R_{2,\,\mu}\,d\mu+\frac{\delta_1\delta_2\nu}{\nu-f(\infty,\infty)}\mathbf1\boxtimes\mathbf1\biggr)\\
&=\lim_{\nu\to\infty}(\varphi_1\boxtimes\varphi_2)\Bigl(\frac{\nu}{\nu-f(\cdot,\cdot)}\Bigr)
=\lim_{\nu\to\infty}(\varphi_1\boxtimes\varphi_2)u=\mathbf1\boxtimes\mathbf1,
\end{align*}
because the functions $\frac{\nu}{\nu-f(\cdot,\cdot)}$ converge to $u$ as $\nu\to\infty$ uniformly on $\Gamma_1\times\Gamma_2$; here $u(\lambda,\mu)=1$.
Consequently, $\infty\notin\bar\sigma(S_{(\cdot)})$.

Conversely, let $\infty\notin\bar\sigma(S_{(\cdot)})$. This means that $S_{(\cdot)}$ \emph{is defined} in a deleted neighbourhood $W$ of infinity and
\begin{equation}\label{e:lim S_nu:2}
\begin{split}
\lim_{\nu\to\infty}\nu S_\nu&=\lim_{\nu\to\infty}\biggl(\frac1{(2\pi
i)^2}\int_{\Gamma_1}\int_{\Gamma_2}\frac{\nu}{\nu-f(\lambda,\mu)}R_{1,\,\lambda}\boxtimes
R_{2,\,\mu}\,d\mu\,d\lambda\\
&+\delta_2\frac1{2\pi i}\int_{\Gamma_1}\frac{\nu}{\nu-f(\lambda,\infty)}R_{1,\,\lambda}\boxtimes
\mathbf1\,d\lambda\\
&+\delta_1\frac1{2\pi i}\int_{\Gamma_2}\frac{\nu}{\nu-f(\infty,\mu)}\mathbf1\boxtimes
R_{2,\,\mu}\,d\mu+\frac{\delta_1\delta_2\nu}{\nu-f(\infty,\infty)}\mathbf1\boxtimes\mathbf1\biggr)=\mathbf1\boxtimes\mathbf1.
\end{split}
\end{equation}
By the definition of $S_{(\cdot)}$, we have that $W\cap f\bigl(\bar\sigma(R_{1,\,(\cdot)}),\bar\sigma(R_{2,\,(\cdot)})\bigr)=\varnothing$. We show that $\infty\notin f\bigl(\bar\sigma(R_{1,\,(\cdot)}),\bar\sigma(R_{2,\,(\cdot)})\bigr)$.

Assuming the contrary, let $f(\lambda_*,\mu_*)=\infty$ for a point $(\lambda_*,\mu_*)\in\bar\sigma(R_{1,\,(\cdot)})\times\bar\sigma(R_{2,\,(\cdot)})$. Then, by
Theorem~\ref{t:spectral map:R_1,R_2->C}, $0\in\bar\sigma(S_\nu)$ for any $\nu\in W$. We also have $0\in\bar\sigma(\nu S_\nu)$ for $\nu\in W$. By Corollaries~\ref{c:f(A,B) is full} and \ref{c:continuity of spr}, it follows that
\begin{equation*}
0\in\bar\sigma\bigl(\lim_{\nu\to\infty}\nu S_\nu\bigr),
\end{equation*}
which contradicts~\eqref{e:lim S_nu:2}.
\end{proof}

Corollary~\ref{c:mero f:2}, see below, answers in the affirmative to the question posed in~\cite{Reed-Simon:BAMS72,Reed-Simon:JFA73} about the independence of the definition of $f(A,B)$ for unbounded operators $A$ and $B$ from the choice of sequences of bounded operators $A_n$ and $B_n$, the resolvents of which converge to the resolvents of $A$ and $B$ respectively.

\begin{corollary}\label{c:mero f:2}
Let the sequences of pseudo-resolvents $R_{n,\,1,\,(\cdot)}$ and $R'_{n,\,1,\,(\cdot)}$ converge\footnote{See the definition on p.~\pageref{page:conv of Res}.} to the same pseudo-resolvent $R_{1,\,(\cdot)}$, and let the sequences of pseudo-resolvents $R_{n,\,2,\,(\cdot)}$ and $R'_{n,\,2,\,(\cdot)}$ converge to the same pseudo-resolvent $R_{2,\,(\cdot)}$. Then both the sequence $S_\nu\bigl(R_{n,\,1,\,(\cdot)},R_{n,\,2,\,(\cdot)}\bigr)$ and the sequence  $S_\nu\bigl(R'_{n,\,1,\,(\cdot)},R'_{n,\,2,\,(\cdot)}\bigr)$ converge to the pseudo-resolvent $S_\nu\bigl(R_{1,\,(\cdot)},R_{2,\,(\cdot)}\bigr)$.
\end{corollary}
\begin{proof}
We make use of definition~\eqref{e:def of S_nu:2}. By Lemma~\ref{l:limit of pseudo-resolvents}, $R_{n,\,1,\,(\cdot)}$ and $R'_{n,\,1,\,(\cdot)}$ converge to $R_{1,\,(\cdot)}$ uniformly on $\Gamma_1$, and $R_{n,\,2,\,(\cdot)}$ and $R'_{n,\,2,\,(\cdot)}$ converge to $R_{2,\,(\cdot)}$ uniformly on $\Gamma_2$. From formula~\eqref{e:def of S_nu:2} it is seen that $S_\nu\bigl(R_{n,\,1,\,(\cdot)},R_{n,\,2,\,(\cdot)}\bigr)$ and $S_\nu\bigl(R'_{n,\,1,\,(\cdot)},R'_{n,\,2,\,(\cdot)}\bigr)$  converge to $S_\nu\bigl(R_{1,\,(\cdot)},R_{2,\,(\cdot)}\bigr)$.
\end{proof}

The theory of meromorphic functions of one operator had its origin in the polynomial functional calculus for unbounded operators constructed in~\cite{Taylor51}, see an exposition in~\cite[ch.~VII, \S~9]{Dunford-Schwartz-I:eng}. Meromorphic functional calculus of one operator was constructed in~\cite{Haase2006}. A spectral mapping theorem for a polynomial of a linear relation was proved in~\cite[Theorem VI.5.4]{Cross}.

Polynomial functions of two unbounded operators were defined in \cite[Theorem 3.4]{Ichinose75}; in particular, a spectral mapping theorem was established, see~\cite[Theorem 3.13]{Ichinose75}. Other analogues of the spectral mapping theorem for analytic functions of unbounded operators (including polynomials) were obtained in~\cite[Theorem 1]{Reed-Simon:BAMS72} and \cite[Theorem 4]{Reed-Simon:JFA73}.

\section{Functional calculus $\varphi_1\boxdot\varphi_2$}\label{s:boxdot}
In this Section, we assume that we are given an extended tensor product $X\boxtimes Y$ of Banach spaces $X$ and $Y$, and we are given pseudo-resolvents $R_{1,\,(\cdot)}$ and $R_{2,\,(\cdot)}$ in the algebras $\mathoo B_0(X)$ and $\mathoo B_0(Y)$ respectively.

We define the mapping $\varphi_1\boxdot\varphi_2$ acting on functions $f\in\mathoo O\bigl(\bar\sigma(R_{1,\,(\cdot)})\cup\bar\sigma(R_{2,\,(\cdot)})\bigr)$ of one variable by the formula
\begin{equation}\label{e: def of boxdot}
(\varphi_1\boxdot\varphi_2)f=\frac1{2\pi i}\int_{\Gamma}f(\lambda)R_{1,\,\lambda}\boxtimes R_{2,\,\lambda}\,d\lambda,
\end{equation}
where $\Gamma$ is an oriented envelope of the union $\bar\sigma(R_{1,\,(\cdot)})\cup\bar\sigma(R_{2,\,(\cdot)})$
of the extended singular sets with respect to the complement of the domain of $f$.

Let $U\subseteq\overline{\mathbb{C}}$ be an open set and $f:\,U\to\mathbb C$ be an analytic function. We call the \emph{divided difference}~\cite{Gelfond:eng,Jordan} of the function $f$ the function $f^{[1]}:\,U\times
U\to\mathbb C$ defined by the formula
\begin{equation}\label{e:f[1]}
f^{[1]}(\lambda,\mu)= \begin{cases}
\frac{f(\lambda)-f(\mu)}{\lambda-\mu}, & \text{if $\lambda\neq\mu$},\\
f'(\lambda), & \text{if $\lambda=\mu$},\\
0,& \text{if $\lambda=\infty$ or $\mu=\infty$}.
 \end{cases}
\end{equation}

\begin{example}\label{ex:f[1]}
We give examples of divided differences of some functions:
\begin{align*}
v_1^{[1]}(\lambda,\mu)&=1,&&\text{where }v_1(\lambda)=\lambda,\\
v_2^{[1]}(\lambda,\mu)&=\lambda+\mu,&&\text{where }v_2(\lambda)=\lambda^2,\\
v_n^{[1]}(\lambda,\mu)&=\lambda^{n-1}+\lambda^{n-2}\mu+\dots+\mu^{n-1},&&\text{where }v_n(\lambda)=\lambda^n,\\
v_{1/2}^{[1]}(\lambda,\mu)&=\frac1{\sqrt\lambda+\sqrt\mu},&&\text{where }v_{1/2}(\lambda)=\sqrt\lambda,\\
r_1^{[1]}(\lambda,\mu)&=\dfrac1{(\lambda_0-\lambda)(\lambda_0-\mu)},&&\text{where }r_1(\lambda)=\frac1{\lambda_0-\lambda},\\
r_n^{[1]}(\lambda,\mu)&=-\dfrac{\frac1{(\lambda_0-\lambda)^n}-\frac1{(\lambda_0-\mu)^n}}{(\lambda_0-\lambda)-(\lambda_0-\mu)}=\\
&=\dfrac{v_n^{[1]}(\lambda_0-\lambda,\lambda_0-\mu)}{(\lambda_0-\lambda)^n(\lambda_0-\mu)^n},&&\text{where }r_n(\lambda)=\frac1{(\lambda_0-\lambda)^n}.
\end{align*}

The Taylor series for the divided difference of a function $f$ at a point $(\lambda_0,\lambda_0)$ has the form
\begin{equation*}
f^{[1]}(\lambda,\mu)=\sum_{n=0}^\infty\frac{f^{(n+1)}(\lambda_0)}{(n+1)!}v_{n+1}^{[1]}(\lambda-\lambda_0,\mu-\lambda_0)
=\sum_{n=0}^\infty\frac{f^{(n+1)}(\lambda_0)}{(n+1)!}\sum_{i=0}^n(\lambda-\lambda_0)^{n-i}(\mu-\lambda_0)^i,
\end{equation*}
where $v_n(\lambda)=\lambda^n$.
In particular, for $\exp_t(\lambda)=e^{\lambda t}$ and $\exp_t^{(1)}(\lambda)=\lambda e^{\lambda t}$ we have
\begin{align*}
\exp_t^{[1]}(\lambda,\mu)&=\sum_{n=0}^\infty\frac{t^n}{(n+1)!}v_{n+1}^{[1]}(\lambda,\mu)
=\sum_{n=0}^\infty\frac{t^n}{(n+1)!}\sum_{i=0}^n\lambda^{n-i}\mu^i,
\\
\exp_t^{(1)\,[1]}(\lambda,\mu)&=\sum_{n=0}^\infty\frac{t^n}{n!}v_{n+1}^{[1]}(\lambda,\mu)
=\sum_{n=0}^\infty\frac{t^n}{n!}\sum_{i=0}^n\lambda^{n-i}\mu^i.
\end{align*}
\end{example}

\begin{proposition}\label{e:f[1] is analytic}
Let $U\subseteq\overline{\mathbb{C}}$ be an open set and $f:\,U\to\mathbb C$ be an analytic function. Then the function $f^{[1]}$ is analytic in $U\times U$.
\end{proposition}
\begin{proof}
The analyticity at a finite point $(\lambda,\mu)$, $\lambda\neq\mu$, is evident. The analyticity at the points of the form $(\lambda,\infty)$ and $(\infty,\mu)$, where $\lambda,\mu\in\mathbb C$, is also evident.

We expand $f$ in the Taylor series about a finite point $\lambda_0\neq\infty$:
\begin{equation*}
f(\lambda)=\sum_{n=0}^\infty c_n(\lambda-\lambda_0)^n.
\end{equation*}
It follows that for $\lambda\neq\mu$ close to $\lambda_0$ one has
\begin{equation*}
f^{[1]}(\lambda,\mu)=\sum_{n=1}^\infty c_nv_{n}^{[1]}(\lambda-\lambda_0,\mu-\lambda_0),
\end{equation*}
where $v_n^{[1]}(\lambda,\mu)=\lambda^{n-1}+\lambda^{n-2}\mu+\dots+\mu^{n-1}$. This series determines an analytic function in a neighbourhood of the point $(\lambda_0,\lambda_0)$. Clearly, $f^{[1]}(\lambda_0,\lambda_0)=f'(\lambda_0)$.

We expand $f$ in the Laurent series with centre $\infty$:
\begin{equation*}
f(\lambda)=\sum_{n=0}^\infty\frac{c_n}{\lambda^n}.
\end{equation*}
This formula shows that for $\lambda\neq\mu$ close to $\infty$ one has
\begin{equation*}
f^{[1]}(\lambda,\mu)=-\sum_{n=1}^\infty c_n\frac{v_{n}^{[1]}(\lambda,\mu)}{\lambda^n\mu^n},
\end{equation*}
where $v_n^{[1]}(\lambda,\mu)=\lambda^{n-1}+\lambda^{n-2}\mu+\dots+\mu^{n-1}$. This series determines an analytic function in a neighbourhood of the point $(\infty,\infty)$. Clearly, $f^{[1]}(\infty,\infty)=0$.
\end{proof}

 \begin{theorem}\label{t:Theta}
Let $f\in\mathoo O\bigl(\bar\sigma(R_{1,\,(\cdot)})\cup\bar\sigma(R_{2,\,(\cdot)})\bigr)$. Then\footnote{Strictly speaking, in this formula, $f^{[1]}$ is understood to be the canonical projection of $f^{[1]}\in\mathoo O\bigl[\bigl(\bar\sigma(R_{1,\,(\cdot)})\cup\bar\sigma(R_{2,\,(\cdot)})\bigr)\times\bigl(\bar\sigma(R_{1,\,(\cdot)})\cup\bar\sigma(R_{2,\,(\cdot)})\bigr)\bigr]$  into $\mathoo O\bigl(\bar\sigma(R_{1,\,(\cdot)})\times\bar\sigma(R_{2,\,(\cdot)})\bigr)$.}
\begin{equation*}
(\varphi_1\boxdot\varphi_2)f=(\varphi_1\boxtimes\varphi_2)f^{[1]}.
\end{equation*}
The spectrum of the transformator
$(\varphi_1\boxdot\varphi_2)f\separ X\boxtimes Y\to X\boxtimes Y$
is given by the formula
\begin{equation*}
\sigma\bigl((\varphi_1\boxdot\varphi_2)f\bigr)
=\bigl\{\,f^{[1]}(\lambda,\mu):\,\lambda\in\bar\sigma(R_{1,\,(\cdot)}),\,\mu\in\bar\sigma(R_{2,\,(\cdot)})\,\bigr\}.
\end{equation*}
 \end{theorem}

 \begin{figure}[thb]
\unitlength=1.mm
\begin{center}
\begin{picture}(65,25)
\put(-30,20){\circle*{30}}\put(-20,20){\circle*{30}}\put(-10,20){\circle*{30}}\put(00,20){\circle*{30}}\put(10,20){\circle*{30}}\put(20,20){\circle*{30}}
\put(47,19){\circle*{30}}\put(46,16){\circle*{30}}\put(45,14){\circle*{30}}\put(43,17){\circle*{30}}
 \put(70,15){\circle{20}}\put(70,15){\oval(34,22)} \put(63,14){\vector(0,1){2}} \put(53,14){\vector(0,1){2}}
 \put(59,12){$\Gamma_1$} \put(88,12){$\Gamma_2$} \put(45,7){$\sigma_1$} \put(9,14){$\bar\sigma_2$} \put(68,12){\boxed{$f$}}
\end{picture}
\caption{The contours $\Gamma_1$ and $\Gamma_2$ from the proof of Theorem~\ref{t:Theta}. The localization of the complement of the domain of $f$ is marked by \boxed{$f$}}\label{f:Gamma2}
\end{center}
 \end{figure}
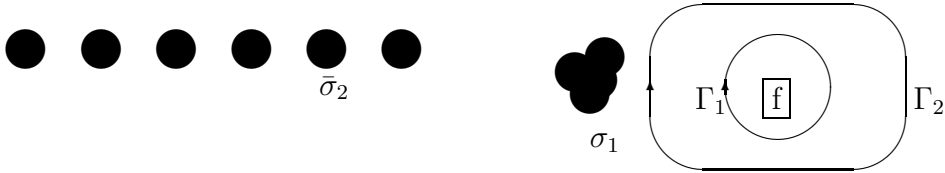

 \begin{proof}
We take contours $\Gamma_1$ and $\Gamma_2$ such that the both are oriented envelopes of $\bar\sigma(R_{1,\,(\cdot)})\cup\bar\sigma(R_{2,\,(\cdot)})$ with respect to the complement of the domain of the function $f$,
and $\Gamma_2$ lies outside of $\Gamma_1$ (so that
$\lambda-\mu$ does not vanish for $\lambda\in\Gamma_1$ and $\mu\in\Gamma_2$), see fig.~\ref{f:Gamma2}.
We make use of the definition:
\begin{align*}
(\varphi_1\boxtimes\varphi_2)f^{[1]}&=\frac1{(2\pi
i)^2}\int_{\Gamma_1}\int_{\Gamma_2}f^{[1]}(\lambda,\mu)R_{1,\,\lambda}\boxtimes
R_{2,\,\mu}\,d\mu\,d\lambda\\
&+\delta_2\frac1{2\pi i}\int_{\Gamma_1}f^{[1]}(\lambda,\infty)R_{1,\,\lambda}\boxtimes
\mathbf1\,d\lambda\\
&+\delta_1\frac1{2\pi i}\int_{\Gamma_2}f^{[1]}(\infty,\mu)\mathbf1\boxtimes
R_{2,\,\mu}\,d\mu+\delta_1\delta_2f^{[1]}(\infty,\infty)\mathbf1\boxtimes\mathbf1\\
&=\frac1{(2\pi
i)^2}\int_{\Gamma_1}\int_{\Gamma_2}\frac{f(\lambda)-f(\mu)}{\lambda-\mu}R_{1,\,\lambda}\boxtimes R_{2,\,\mu}\,d\mu\,d\lambda.
\end{align*}
We represent the last integral as the sum of two iterated integrals:
\begin{align}
&\phantom{+\;\,}\frac1{2\pi i}\int_{\Gamma_1}R_{1,\,\lambda}\boxtimes\Bigl(\frac{f(\lambda)}{2\pi
i}\int_{\Gamma_2}\frac{1}{\lambda-\mu}R_{2,\,\mu}\,d\mu\Bigr)\,d\lambda\label{e:integral 1}\\
&+\frac1{2\pi i}\int_{\Gamma_2}\Bigl(\frac{f(\mu)}{2\pi
i}\int_{\Gamma_1}\frac{1}{\mu-\lambda}R_{1,\,\lambda}\,d\lambda\Bigr)\boxtimes R_{2,\,\mu}\,d\mu.\label{e:integral
2}
\end{align}
By the Cauchy integral formula, for the internal integral in~\eqref{e:integral 1} we have
\begin{equation*}
\frac{1}{2\pi i}\int_{\Gamma_2}\frac{1}{\lambda-\mu}R_{2,\,\mu}\,d\mu=R_{2,\,\lambda},
\end{equation*}
and by the Cauchy integral theorem, for the internal integral in~\eqref{e:integral 2} we have
\begin{equation*}
\frac{1}{2\pi i}\int_{\Gamma_1}\frac{1}{\mu-\lambda}R_{1,\,\lambda}\,d\lambda=0,
\end{equation*}
since, by the assumption, the contour $\Gamma_1$ does not surround the singularities of the function
$\lambda\mapsto\frac{1}{\mu-\lambda}$, $\mu\in\Gamma_2$, and the pseudo-resolvent $\lambda\mapsto
R_{1,\,\lambda}$. Thus, the original integral takes the form
\begin{equation*}
\frac1{2\pi i}\int_{\Gamma_1}f(\lambda)R_{1,\,\lambda}\boxtimes R_{2,\,\lambda}\,d\lambda=\bigl[(\varphi_1\boxdot\varphi_2)f\bigr].
\end{equation*}

The second formula follows from~\ref{t:spectral map:R_1,R_2->C}.
 \end{proof}

\begin{theorem}\label{t:Theta via Q:0}
Let $A\in\mathoo{B}(X)$, $B\in\mathoo{B}(Y)$, and $f\in\mathoo O\bigl(\sigma(A)\cup\sigma(B)\bigr)$. Then
\begin{equation*}
\varphi_A(f)\boxtimes\mathbf1-\mathbf1\boxtimes\varphi_B(f)=\bigl[(\varphi_A\boxdot\varphi_B)f\bigr](A\boxtimes\mathbf1-\mathbf1\boxtimes B),
\end{equation*}
where the functional calculi $\varphi_A$ and $\varphi_B$ are constructed by $A$ and $B$ respectively.
\end{theorem}
 \begin{proof}
The proof follows from the identity
\begin{equation*}
f(\lambda)-f(\mu)=f^{[1]}(\lambda,\mu)(\lambda-\mu)
\end{equation*}
and Theorems~\ref{t:Theta} and \ref{t:functional_calculas:infty in rho:R_1,R_2->C}.
 \end{proof}

In the following corollary, we describe a representation for the increment of an analytic function.

\begin{corollary}\label{c:f(A)-f(B)}
Let $A,B\in\mathoo{B}(X)$ and $f\in\mathoo O\bigl(\sigma(A)\cup\sigma(B)\bigr)$. Then
\begin{equation*}
f(A)-f(B)
=\frac1{2\pi i}\int_\Gamma f(\lambda)(\lambda\mathbf1-A)^{-1}(A-B)(\lambda\mathbf1-B)^{-1}\,d\lambda,
\end{equation*}
where the functional calculi $\varphi_A$ and $\varphi_B$ are constructed by $A$ and $B$ respectively, and $\varphi_A\boxtimes\varphi_B$ acts in the extended tensor product $\mathoo B(X,X)$, see Example~\ref{ex:TP of cat}(e).
\end{corollary}
For the function $f=\exp_t$, this formula was found in~\cite[p. 978]{van_Loan77}.
\begin{proof}
We apply the formula from Theorem~\ref{t:Theta via Q:0} to the operator $C=\mathbf1$, assuming that $X=Y$. We have (taking into account that $C=\mathbf1$)
\begin{gather*}
(A\boxtimes\mathbf1-\mathbf1\boxtimes B)C=AC-CB=A-B,\\
\bigl[(\varphi_A\boxdot\varphi_B)f\bigr](A-B)
=\frac1{2\pi i}\int_\Gamma f(\lambda)(\lambda\mathbf1-A)^{-1}(A-B)(\lambda\mathbf1-B)^{-1}\,d\lambda,\\
\bigl(\varphi_A(f)\boxtimes\mathbf1-\mathbf1\boxtimes\varphi_B(f)\bigr)C=\varphi_A(f)C-C\varphi_B(f)=f(A)C-Cf(B)=f(A)-f(B).\qed
\end{gather*}
\renewcommand\qed{}
\end{proof}

One of the primary ideas~\cite{Higham08,Frommer-Simoncini,Kurbatov-Oreshina04,Kurbatov-Kurbatova-EMJ-2012,Mathias93,Dieci-Papini} of approximate calculation of an analytic function $f$ of an operator or a pseudo-resolvent consists in an approximation of $f$ by a polynomial or a rational function. In the case of $(\varphi_1\boxdot\varphi_2)f$, for applying this idea it is necessary to be able to calculate $\varphi_1\boxdot\varphi_2$ at least of monomials and elementary rational functions. Formulae of this kind are presented in Corollary~\ref{c:ratio appr} below.

 \begin{corollary}\label{c:ratio appr}
If the pseudo-resolvents $R_{1,\,(\cdot)}$ and $R_{2,\,(\cdot)}$ are generated by the operators $A$ and $B$ respectively, then
\begin{equation*}
\bigl(\varphi_A\boxdot\varphi_B\bigr)(v_n)=A^{n-1}\boxtimes\mathbf{1}+A^{n-2}\boxtimes B+\dots+\mathbf{1}\boxtimes B^{n-1},\qquad\text{where }v_n(\lambda)=\lambda^n.
\end{equation*}

If
$\lambda_0\in\rho(R_{1,\,(\cdot)})\cap\rho(R_{2,\,(\cdot)})$, then
\begin{equation*}
\bigl(\varphi_1\boxdot\varphi_2\bigr)(r_1)=R_{1,\,\lambda_0}\boxtimes R_{2,\,\lambda_0},\qquad\text{where }r_1(\lambda)=\frac1{\lambda_0-\lambda}.
\end{equation*}
If, in addition, the extended singular sets of the pseudo-resolvents $R_{1,\,(\cdot)}$ and $R_{2,\,(\cdot)}$ are disjoint, then
\begin{equation*}
\bigl(\varphi_1\boxdot\varphi_2\bigr)(r_n)=-\bigl(R_{1,\,\lambda_0}^n\boxtimes\mathbf1-\mathbf1\boxtimes R_{2,\,\lambda_0}^n\bigr)\bigl(R_{1,\,\lambda_0}\boxtimes\mathbf1-\mathbf1\boxtimes R_{2,\,\lambda_0}\bigr)^{-1},\;\;\text{where }r_n(\lambda)=\frac1{(\lambda_0-\lambda)^n}.
\end{equation*}
If, in addition, the pseudo-resolvents $R_{1,\,(\cdot)}$ and $R_{2,\,(\cdot)}$ are generated by the operators $A$ and $B$ respectively, then
\begin{equation*}
\bigl(\varphi_A\boxdot\varphi_B\bigr)(r_n)=\bigl((\lambda_0\mathbf1-A)^{n-1}\boxtimes\mathbf{1}
+\dots+\mathbf{1}\boxtimes(\lambda_0\mathbf1-B)^{n-1}\bigr)\bigl(R_{A,\,\lambda_0}^n\boxtimes R_{B,\,\lambda_0}^n\bigr).
\end{equation*}
 \end{corollary}
 \begin{proof}
It suffices to make use of Example~\ref{ex:f[1]}
and to apply Theorem~\ref{t:Theta} and Corollary~\ref{c:functional_calculas:infty in rho:R_1,R_2->C}.
 \end{proof}

 \begin{corollary}\label{c:mult} Let $g,h\in\mathoo O\bigl(\bar\sigma(R_{1,\,(\cdot)})\cup\bar\sigma(R_{2,\,(\cdot)})\bigr)$. Then
\begin{equation*}
\bigl(\varphi_1\boxdot\varphi_2\bigr)(gh)=\bigl[\bigl(\varphi_1\boxdot\varphi_2\bigr)(g)\bigr]\bigl(\mathbf1\boxtimes\varphi_2(h)\bigr)
+\bigl(\varphi_1(g)\boxtimes\mathbf1\bigr)\bigl[\bigl(\varphi_1\boxdot\varphi_2\bigr)(h)\bigr],
\end{equation*}
where $(gh)(\lambda)=g(\lambda)h(\lambda)$.
 \end{corollary}
 \begin{proof}
By Theorem~\ref{t:Theta}, this formula is equivalent to the identity
\begin{equation*}
\bigl(\varphi_1\boxtimes\varphi_2\bigr)(gh)^{[1]}=
\bigl[\bigl(\varphi_1\boxtimes\varphi_2\bigr)g^{[1]}\bigr]\bigl(\mathbf1\boxtimes\varphi_2(h)\bigr)
+\bigl(\varphi_1(g)\boxtimes\mathbf1\bigr)\bigl[\bigl(\varphi_1\boxtimes\varphi_2\bigr)h^{[1]}\bigr].
\end{equation*}
We have
\begin{align*}
\frac{g(\lambda)h(\lambda)-g(\mu)h(\mu)}{\lambda-\mu}&=
\frac{g(\lambda)h(\lambda)-g(\mu)h(\lambda)+g(\mu)h(\lambda)-g(\mu)h(\mu)}{\lambda-\mu}\\
&=\frac{g(\lambda)h(\lambda)-g(\mu)h(\lambda)}{\lambda-\mu}
+\frac{g(\mu)h(\lambda)-g(\mu)h(\mu)}{\lambda-\mu}\\
&=\frac{g(\lambda)-g(\mu)}{\lambda-\mu}h(\lambda)
+g(\mu)\frac{h(\lambda)-h(\mu)}{\lambda-\mu}.
\end{align*}
Taking into account the ability of passages to the limits as $\lambda-\mu\to0$ and $\lambda-\mu\to\infty$ we arrive at
\begin{equation*}
(gh)^{[1]}(\lambda,\mu)=g(\lambda)h^{[1]}(\lambda,\mu)+g^{[1]}(\lambda,\mu)h(\mu).
\end{equation*}
It remains to apply Theorems~\ref{t:functional_calculas:infty in rho:R_1,R_2->C},~\ref{t:functional_calculas:infty in rho+sigma:R_1,R_2->C}, and \ref{t:functional_calculas:infty in sigma:R_1,R_2->C}.
 \end{proof}

The function $\beta_{g,\,h}:\,U\times U\to\mathbb C$ defined by the formula
\begin{equation*}
\beta_{g,\,h}(\lambda,\mu)=\begin{cases}
\frac{g(\lambda)h(\mu)-h(\lambda)g(\mu)}{\lambda-\mu}, & \text{if $\lambda\neq\mu$},\\
g'(\lambda)h(\mu)-h'(\lambda)g(\mu), & \text{if $\lambda=\mu$},\\
0,& \text{if $\lambda=\infty$ or $\mu=\infty$},
 \end{cases}
\end{equation*}
is similar to the divided difference. It is generated by two analytic functions
$g,h:\,U\to\mathbb C$. By analogy with~\cite{Krein-Naimark,Heinig:eng,Heinig-Rost}, we call the function $\beta_{g,\,h}$ the \emph{Bezoutian}.
The Bezoutian is a difference-differential analogue of the Wronskian. For example, the Bezoutian of the functions $\sin$ and $\cos$ is $\text{sinc}(\lambda-\mu)=\frac{\sin(\lambda-\mu)}{\lambda-\mu}$.
We note that the Bezoutian can be expressed in terms of divided differences:
\begin{equation*}
\beta_{g,h}(\lambda,\mu)=g^{[1]}(\lambda,\mu)h(\mu)-h^{[1]}(\lambda,\mu)g(\mu).
\end{equation*}
(In particular, this formula and Proposition~\ref{e:f[1] is analytic} imply that $\beta_{g,\,h}$ is an analytic function.)
Conversely,
\begin{equation*}
g^{[1]}(\lambda,\mu)=\beta_{g,u}(\lambda,\mu),
\end{equation*}
where $u(\lambda)=1$.

 \begin{corollary}\label{c:div} Let $g,h\in\mathoo O\bigl(\bar\sigma(R_{1,\,(\cdot)})\cup\bar\sigma(R_{2,\,(\cdot)})\bigr)$ and let $h(\lambda)\neq0$ for $\lambda\in\bar\sigma(R_{1,\,(\cdot)})\cup\bar\sigma(R_{2,\,(\cdot)})$. Then
\begin{equation*}
\bigl(\varphi_1\boxdot\varphi_2\bigr)\Bigl(\frac gh\Bigr)=\bigl[\bigl(\varphi_1\boxtimes\varphi_2\bigr)(\beta_{g,h})\bigr]\bigl[\varphi_1(h)\boxtimes\varphi_2(h)\bigr]^{-1},
\end{equation*}
where $\Bigl(\dfrac gh\Bigr)(\lambda)=\frac{g(\lambda)}{h(\lambda)}$.
 \end{corollary}
 \begin{proof}
The proof is analogous to that of Corollary~\ref{c:mult} and follows from the formula
 \begin{equation*}
 \Bigl[\frac gh\Bigr]^{[1]}(\lambda,\mu)=\frac{\beta_{g,h}(\lambda,\mu)}{h(\lambda)h(\mu)}.\qed
 \end{equation*}
 \renewcommand\qed{}
 \end{proof}

 \begin{proposition}\label{p:phi(f')}
Let $X$ be a Banach space. We take for an extended tensor product the space $\mathoo B(X,X)$ {\rm(}see Example~\ref{ex:TP of cat}{\rm(e))}, and we take for $R_{1,\,(\cdot)}$ and $R_{2,\,(\cdot)}$ the resolvent $R_{(\cdot)}$ of the same operator $A\in\mathoo B(X)$. Let an operator $C$ commute with at least one value $R_\mu$ of the pseudo-resolvent $R_{(\cdot)}$.
Then
\begin{equation*}
\bigl[\bigl(\varphi_A\boxdot\varphi_A\bigr)f\bigr]C=f'(A)C=Cf'(A).
\end{equation*}
 \end{proposition}

 \begin{proof}
We note that, by virtue of Theorem~\ref{t:pseudoresovent}, $C$ commutes with all values $R_\lambda$ of the pseudo-resolvent $R_{(\cdot)}$.

By the definition and commutativity, we have
\begin{align*}
\bigl[\bigl(\varphi\boxdot\varphi\bigr)f\bigr](C)=\frac1{2\pi i}\int_{\Gamma}f(\lambda)R_\lambda CR_\lambda\,d\lambda
=\frac1{2\pi i}\int_{\Gamma}f(\lambda)R_{\lambda}^{2}C\,d\lambda
=\Bigl(\frac1{2\pi
i}\int_{\Gamma}f(\lambda)R_{\lambda}^{2}\,d\lambda\Bigr)C.
\end{align*}
Passing to the limit in the Hilbert identity~\eqref{e:Hilbert identity:pseudo} we obtain the relation
\begin{equation*}
R_{\lambda}^{2}=-R_{\lambda}',\qquad\lambda\notin\sigma\bigl(R_{(\cdot)}\bigr).
\end{equation*}
Substituting this identity into the previous equality and integrating by parts we obtain
\begin{equation*}
\bigl[\bigl(\varphi\boxdot\varphi\bigr)f\bigr](C)=-\Bigl(\frac1{2\pi
i}\int_{\Gamma}f(\lambda)R_{\lambda}'\,d\lambda\Bigr)C=\Bigl(\frac1{2\pi
i}\int_{\Gamma}f'(\lambda)R_{\lambda}\,d\lambda\Bigr)C=\varphi(f')C.\qed
\end{equation*}
\renewcommand\qed{}
 \end{proof}

We note that the divided differences $f^{[1]}(A,B)$ of the operators $A$ and $B$ are also closely related to the calculation of functions of block triangular matrices~\cite{Davis_Ch,Parlett76,Davies-Higham2003,Golub-Van_Loan96,Higham08}.

\section{The impulse response}\label{s:pulse char}
In subsequent sections, we discuss some applications.

In this Section, the previous results are applied to the representation of the impulse response of a second order differential equation.
Here we regard the space $\mathoo B(Y,X)$ (example~\ref{ex:TP of cat}(e)) as an extended tensor product.
Therefore, for example, the action of the transformator $\varphi_1(g)\boxtimes\varphi_2(h)$ on the operator $C\in\mathoo B(Y,X)$
results in the operator $\varphi_1(g)C\varphi_2(h)$.

Let $X$ and $Y$ be Banach spaces and $E,F,H\in\mathoo B(Y,X)$.
A function $\lambda\mapsto\lambda^2E+\lambda F+H$, where $\lambda\in\mathbb C$, is called~\cite{Marcus:eng,Krein-Langer:eng,Gantmakher:eng}
a \emph{square pencil}.
The \emph{resolvent set} of the pencil is the set $\rho(E,F,H)$ of all $\lambda\in\mathbb C$ such that the operator
$\lambda^2E+\lambda F+H$ is invertible. The \emph{spectrum} is the compliment $\sigma(E,F,H)=\mathbb C\setminus\rho(E,F,H)$ and the \emph{resolvent} is the function
\begin{equation}\label{e:resolvent2}
R_\lambda=(\lambda^2E+\lambda F+H)^{-1},\qquad\lambda\in\rho(E,F,H).
\end{equation}

The main sources~\cite{Marcus:eng,Tisseur-Meerbergen} of square pencils are the second order differential equations of the form
 \begin{equation}\label{e:DE 2nd order}
E\ddot y(t)+F\dot y(t)+Hy(t)=0,
 \end{equation}
where $y:\,\mathbb R\to Y$. In this Section, it is always assumed that the operator $E$ is
invertible\footnote{We note that even if $E$ is invertible, the multiplication of the equation~\eqref{e:DE
2nd order} by $E^{-1}$ is not always desirable. For example, the operators $E,F,H$ are often assumed~\cite{Krein-Langer:eng, Seshu-Reed:eng} to be
self-adjoint, but the multiplication by $E^{-1}$ may cause to the loss of this property.}. We recall the following proposition.

\begin{proposition}[{\rm see, for example,~\cite[Theorem 16]{Kurbatova-POMI:eng}}]\label{p:impulse char}
Let the operator $E$ be invertible. Then the solution of the initial value problem
\begin{align*}
E\ddot y(t)+F\dot y(t)+Hy(t)&=0, \\
y(0)&=y_0,\\
\dot y(0)&=y_1
\end{align*}
can be represented in the form
\begin{equation*}
y(t)=\dot{\mathcal T}(t)Ey_0+\mathcal T(t)(Ey_1+Fy_0),
\end{equation*}
where
\begin{align*}
\mathcal T(t)&=\frac1{2\pi i}\int_\Gamma\exp_t(\lambda)(\lambda^2E+\lambda F+H)^{-1}\,d\lambda, \\
\dot{\mathcal T}(t)&=\frac1{2\pi i}\int_\Gamma\exp_t^{(1)}(\lambda)(\lambda^2E+\lambda F+H)^{-1}\,d\lambda,
\end{align*}
$\Gamma$ is an oriented envelope of the pencil spectrum $\sigma(E,F,H)$, and
\begin{equation*}
\exp_t(\lambda)=e^{\lambda t},\qquad \exp_t^{(1)}(\lambda)=\lambda e^{\lambda t}.
\end{equation*}
\end{proposition}
It can be shown that the function $\mathcal T$ is the impulse response, and $\dot{\mathcal T}$ is its derivative.

A \emph{factorization} of the pencil is the representation of its resolvent in the form
 \begin{equation}\label{e:factor of res}
R_\lambda=R_{1,\,\lambda}CR_{2,\,\lambda},
 \end{equation}
where $R_{1,\,(\cdot)}$ and $R_{2,\,(\cdot)}$ are pseudo-resolvents acting in $X$ and $Y$ respectively, and
$C\in\mathoo B(Y,X)$. It is assumed that $\rho(R_{1,\,(\cdot)})\cap\rho(R_{2,\,(\cdot)})\supseteq\rho(E,F,H)$.

\begin{proposition}\label{e:C=E}
Let the operator $E$ be invertible. Then we have $C=E$ in formula~\eqref{e:factor of res}, and the pseudo-resolvents $R_{1,\,(\cdot)}$ and $R_{2,\,(\cdot)}$ are the resolvents of some operators $A_1$ and $A_2$.
\end{proposition}
\begin{proof}
By Proposition~\ref{c:pseudo in an infinite point}, we have
\begin{align*}
R_{1,\,\lambda}=-N_1+\frac {P_1}\lambda+\frac
{A_1}{\lambda^2}+\frac{A_1^2}{\lambda^3}+\frac{A_1^3}{\lambda^4}+\ldots, \\
R_{2,\,\lambda}=-N_2+\frac {P_2}\lambda+\frac
{A_2}{\lambda^2}+\frac{A_2^2}{\lambda^3}+\frac{A_2^3}{\lambda^4}+\ldots.
\end{align*}
Hence,
\begin{equation*}
R_{1,\,\lambda}CR_{2,\,\lambda}=N_1CN_2-\frac {P_1CN_2+N_1CP_2}\lambda+\frac{-A_1CN_2+P_1CP_2-N_1CA_2}{\lambda^2}
+\ldots.
\end{equation*}
On the other hand,
\begin{equation*}
R_\lambda=\frac{E}{\lambda^2}-\frac{E^{-1}FE{-1}}{\lambda^3}+\ldots.
\end{equation*}
Therefore,
\begin{equation*}
N_1CN_2=\mathbf{0},\qquad P_1CN_2+N_1CP_2=\mathbf{0},\qquad -A_1CN_2+P_1CP_2-N_1CA_2=E.
\end{equation*}
Multiplying the second equation on the left by $A_1P_1$ (keeping in mind the identities $P^{2}=P$, $AP=PA=A$ and $NP=PN=\mathbf0$ from Proposition~\ref{c:pseudo in an infinite point}), we arrive at
\begin{equation*}
A_1CN_2=\mathbf{0}.
\end{equation*}
Similarly, multiplying the second equation on the right by $A_2$, we have
\begin{equation*}
N_1CA_2=\mathbf{0}.
\end{equation*}
Substituting these results into the third equality, we obtain
\begin{equation*}
P_1CP_2=E.
\end{equation*}
Because of the invertibility of $E$, it follows that the projectors $P_1$ and $P_2$ coincide with $\mathbf{1}$, and $C=E$. Consequently (by the identity $NP=PN=\mathbf0$, see Proposition~\ref{c:pseudo in an infinite point}), we have $N_1=\mathbf0$ and $N_2=\mathbf0$. It follows that $\lim_{\lambda\to\infty}\lambda R_{1,\,\lambda}=\mathbf1$ and $\lim_{\lambda\to\infty}\lambda R_{2,\,\lambda}=\mathbf1$. By Proposition~\ref{p:1 in B_R}(c), these means that the pseudo-resolvents $R_{1,\,(\cdot)}$ and $R_{2,\,(\cdot)}$ are the resolvents
of the operators $A_1$ and $A_2$.
\end{proof}

 \begin{theorem}\label{t:pulse char:2}
Let the operator $E$ be invertible, and the square pencil admits factorization~\eqref{e:factor of res}.
Then the impulse response $\mathcal T$ and its derivative $\dot{\mathcal T}$ can be represented in the form
\begin{align*}
\mathcal T(t)&=\bigl(\varphi_1\boxdot\varphi_2\bigr)(\exp_t)C=\bigl(\varphi_1\boxtimes\varphi_2\bigr)(\exp_t^{[1]})C,\\
\dot{\mathcal T}(t)&=\bigl(\varphi_1\boxdot\varphi_2\bigr)(\exp_t^{(1)})C=\bigl(\varphi_1\boxtimes\varphi_2\bigr)(\exp_t^{(1)\,[1]})C,
\end{align*}
where $\exp_t^{(1)\,[1]}(\lambda,\mu)=\frac{\lambda e^{\lambda t}-\mu e^{\mu t}}{\lambda-\mu}$
for $\lambda\neq\mu$.
 \end{theorem}
 \begin{proof}
The proof follows from Proposition~\ref{p:impulse char} and Theorem~\ref{t:Theta}.
 \end{proof}

 \begin{corollary}\label{c:spectra of H}
The spectra of the transformators $C\mapsto \mathcal T(t)$ and $C\mapsto \dot{\mathcal T}(t)$ in the algebra $\mathoo B\bigl(\mathoo
B(Y,X)\bigr)$ are equal to
\begin{equation*}
\bigl\{\,\exp_t^{[1]}(\lambda,\mu):\,\lambda\in\sigma(A),\;\mu\in\sigma(B)\,\bigr\},\qquad
\bigl\{\,\exp_t^{(1)\,[1]}(\lambda,\mu):\,\lambda\in\sigma(A),\;\mu\in\sigma(B)\,\bigr\},
\end{equation*}
respectively.
 \end{corollary}
 \begin{proof}
The proof follows from Theorems~\ref{t:spectral map:R_1,R_2->C} and~\ref{t:pulse char:2}.
 \end{proof}

 \begin{remark}\label{r:Kenney-Laub98:tanh}
(a) In article~\cite{Kenney-Laub98}, for the approximate calculation of expressions of the type
$\bigl(\varphi_1\boxtimes\varphi_2\bigr)(\exp_t^{[1]})$ it
is suggested to use the following representation (written in other notations and verified directly)
\begin{equation*}
\exp_t^{[1]}(\lambda,\mu)=(e^{\lambda t}+e^{\mu t})\frac{\tanh
\bigl(\frac{\lambda-\mu}{2}t\bigr)}{\frac{\lambda-\mu}{2}t},\qquad\lambda\neq\mu.
\end{equation*}
By Theorem~\ref{t:functional_calculas:infty in rho:R_1,R_2->C},
$\bigl(\varphi_1\boxtimes\varphi_2\bigr)(e^{\lambda t}+e^{\mu t})$ is
$\varphi_1(\exp_t)\boxtimes\mathbf1+\mathbf1\boxtimes\varphi_2(\exp_t)$. By
Corollary~\ref{c:functional_calculas:infty in rho:R_1,R_2->C:lambda-mu}, the operator
$\bigl(\varphi_1\boxtimes\varphi_2\bigr)\Bigl(\frac{\tanh
\bigl(\frac{\lambda-\mu}{2}t\bigr)}{\frac{\lambda-\mu}{2}t}\Bigr)$ is the function
$\tau(z)=\frac{\tanh\bigl(\frac{z}{2}\bigr)}{\frac{z}{2}}$ of the transformator
$(A\boxtimes\mathbf1-\mathbf1\boxtimes B)t$. The function $\tau$ is
analytic in the circle $|z|<\pi$. In~\cite{Kenney-Laub98}, for its computation, it is suggested to use the
Taylor polynomials or rational approximations.

(b) Formulae
\begin{equation*}
\exp_t^{[1]}(\lambda,\mu)=e^{\frac{(\lambda+\mu)t}2}\frac{\sinh\bigl(\frac{\lambda-\mu}2t\bigr)}{\frac{\lambda-\mu}2t},\qquad
\exp_t^{[1]}(\lambda,\mu)=e^{\mu t}\frac{e^{(\lambda-\mu)t}-1}{\lambda-\mu},
\end{equation*}
assuming a similar usage, are suggested in book~\cite[formulae (10.17)]{Higham08}.

(c) We present a formula that enables one to apply similar ideas for the calculation of $\exp_t^{(1)\,[1]}$:
\begin{equation}\label{e:d exp_t(1):0}
\begin{split}
\exp_t^{(1)[1]}(\lambda,\mu)&=\frac{\lambda e^{\lambda t}-\mu e^{\mu t}}{\lambda-\mu}=\frac{\lambda e^{\lambda t}-\mu e^{\lambda t}+\mu e^{\lambda t}-\mu e^{\mu t}}{\lambda-\mu}
= e^{\lambda t}+\mu\frac{e^{\lambda t}-e^{\mu t}}{\lambda-\mu}=\\
&=e^{\lambda t}+\mu\exp_t^{[1]}(\lambda,\mu).
\end{split}
\end{equation}
 \end{remark}

 \begin{corollary}\label{c:semigroup property for T}
Let $E=\mathbf1$. Then
\begin{equation*}
\mathcal T(t+s)=\mathcal T_1(t)\mathcal T(s)+\mathcal T(t)\mathcal T_2(s),
\end{equation*}
where
\begin{align*}
\mathcal T_1(t)&=\varphi_1(\exp_t)=\frac1{2\pi i}\int_{\Gamma_1}\exp_t(\lambda)R_{1,\,\lambda}\,d\lambda, \\
\mathcal T_2(t)&=\varphi_2(\exp_t)=\frac1{2\pi i}\int_{\Gamma_2}\exp_t(\mu)R_{2,\,\mu}\,d\mu.
\end{align*}
 \end{corollary}
 \begin{proof}
This is a special case of Corollary~\ref{c:mult}.
 \end{proof}

Issues related to  factorization are widely discussed in the
literature~\cite{Keldysh1951:eng,Krein-Langer:eng,Marcus:eng,Davis81,Kostyuchenko-Shkalikov:eng,Shkalikov-Pliev:eng,Tisseur-Meerbergen,Kurbatov-Oreshina2003:eng}.
The factorization of an operator pencils of an arbitrary order is discussed in~\cite{Marcus:eng,Isaev:eng,Zayachkovskii-Pankov1:eng,Zayachkovskii-Pankov2:eng,Maroulas-Psarrakos,Krupnik:eng, Gohberg-Lancaster-Rodman:Matrix_polynomials,Guo-Higham-Tisseur,Harbarth-Langer,Gohberg-Lancaster-Rodman:Matrix_polynomials}.

Estimates of the norms of operators $\bigl(\varphi_1\boxtimes\varphi_2\bigr)(f)C$ are obtained in~\cite{Gil-11:EJLA,Gil-15:RMJM};
special attention is paid to $\mathcal T(t)$ and $\dot{\mathcal T}(t)$. Estimates of the norm of
$e^{(A\otimes\mathbf1+\mathbf1\otimes B)t}$ are given in~\cite{Benzi-Simoncini:1501.07376}.

\section{The transformator $Q$ and the Sylvester equation}\label{s:Sylvester}
It often arises the problem of calculating the transformator $Q=(\varphi_1\boxtimes\varphi_2)w$, where
\begin{equation*}
w(\lambda,\mu)=\frac1{\lambda-\mu}.
\end{equation*}
As a rule, it is equivalent to solving the Sylvester equation. In this Section, we discuss some properties of the transformator $Q$.

Let $X\boxtimes Y$ be an extended tensor product of Banach spaces $X$ and $Y$, and
$R_{1,\,(\cdot)}$ and $R_{2,\,(\cdot)}$ be pseudo-resolvents in the algebras $\mathoo B_0(X)$ and $\mathoo B_0(Y)$ respectively.
We assume that the extended singular sets $\bar\sigma(R_{1,\,(\cdot)})$ and $\bar\sigma(R_{2,\,(\cdot)})$ are disjoint. We consider the transformator $Q$ defined as
\begin{equation*}
Q=(\varphi_1\boxtimes\varphi_2)w,
\end{equation*}
where\footnote{The function $w$ is meromorphic with the point of indeterminacy $(\infty,\infty)$.}
\begin{equation*}
w(\lambda,\mu)=\frac1{\lambda-\mu}.
\end{equation*}
If necessary, we will use the more detailed notation $Q_{\varphi_1,\varphi_2}$ or $Q_{A,B}$.

\begin{proposition}\label{p:Q via boxdot}
We assume that the extended singular sets $\bar\sigma(R_{1,\,(\cdot)})$ and $\bar\sigma(R_{2,\,(\cdot)})$ of the pseudo-resolvents  $R_{1,\,(\cdot)}$ and $R_{2,\,(\cdot)}$ are disjoint.
Then
\begin{equation}\label{e:Q via boxdot}
Q=\frac12\bigl(\varphi_1\boxdot\varphi_2\bigr)(\sgn_{1|2}),
\end{equation}
where the function $\sgn_{1|2}$ is equal to 1 in a neighborhood of the extended singular set $\bar\sigma(R_{1,\,(\cdot)})$ and is equal to $-1$ in a neighborhood of the extended singular set $\bar\sigma(R_{2,\,(\cdot)})$.
The transformator $Q$ can be represented in the form
\begin{equation}\label{e:Q via contour}
Q=\frac1{2\pi i}\int_\Gamma R_{1,\,\lambda}\boxtimes R_{2,\,\lambda}\,d\lambda,
\end{equation}
where $\Gamma$ is an oriented envelope of $\bar\sigma(R_{1,\,(\cdot)})$ with respect to  $\bar\sigma(R_{2,\,(\cdot)})$.
\end{proposition}

 \begin{figure}[htb]
\unitlength=1.mm
\begin{center}
\begin{picture}(55,20)
\put(0,15){\circle*{10}} \put(0,15){\circle{25}} \put(2,12){$\sigma_1$} \put(5,6){$\Gamma_1$} \put(7,15){\vector(0,1){1}}
\put(25,15){\circle*{10}} \put(25,15){\circle{20}} \put(27,12){$\sigma_2$} \put(30,6){$\Gamma_2$} \put(32,15){\vector(0,1){1}}
\end{picture}\hfil
\begin{picture}(45,20)
\put(7,18){\circle*{10}}\put(8,14){\circle*{10}}\put(5,17){\circle*{10}}
\put(30,15){\circle{25}} \put(2,10){$\bar\sigma_1$} \put(46,6){$\Gamma_1$} \put(45,15){\vector(0,-1){1}}
\put(30,15){\circle*{10}} \put(30,15){\oval(30,22)} \put(32,12){$\sigma_2$} \put(35,6){$\Gamma_2$} \put(37,15){\vector(0,1){1}}
\end{picture}
\caption{Various options of an arrangement of the contours $\Gamma_1$ and $\Gamma_2$ and the extended singular sets}\label{f:Gamma3}
\end{center}
 \end{figure}
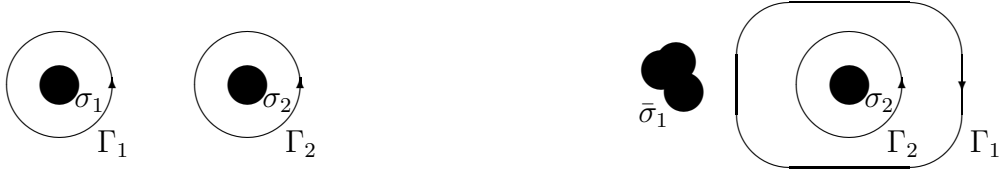

\begin{proof}
It is easy to verify that $\sgn_{1|2}^{[1]}=2w$. So, formula~\eqref{e:Q via boxdot} follows from Theorem~\ref{t:Theta}.

We calculate $(\varphi_1\boxtimes\varphi_2)w$. To be definite, we assume that $\infty\notin\bar\sigma(R_{2,\,(\cdot)})$. We assume that the oriented envelope $\Gamma_1$ of the set $\bar\sigma(R_{1,(\cdot)})$ and the oriented envelope
$\Gamma_2$ of the set $\sigma(R_{2,(\cdot)})$ are located as shown in Fig.~\ref{f:Gamma3}. In particular, $\lambda-\mu$ is not equal to zero for $\lambda\in\Gamma_1$ and $\mu\in\Gamma_2$.
We have (note that in representation~\eqref{e:delta1,delta2} for the function $w$, we have $\delta_1=\delta_2=0$)
\begin{align*}
Q&=(\varphi_1\boxtimes\varphi_2)w=\frac{1}{(2\pi i)^2}\int_{\Gamma_1}\int_{\Gamma_2}\frac1{\lambda-\mu}R_{1,\lambda}\boxtimes R_{2,\mu}\,d\mu\,d\lambda\\
&=\frac{1}{2\pi i}\int_{\Gamma_1}R_{1,\lambda}\boxtimes \Bigl(\frac{1}{2\pi i}\int_{\Gamma_2}\frac{R_{2,\mu}}{\lambda-\mu}\,d\mu\Bigr)\,d\lambda
=\frac{1}{2\pi i}\int_{\Gamma_1}R_{1,\lambda}\boxtimes R_{2,\lambda}\,d\lambda.
\end{align*}
Obviously, $\Gamma_1$ is an oriented envelope of $\bar\sigma(R_{1,(\cdot)})$ with respect to $\sigma(R_{2,(\cdot)})$.
\end{proof}

\begin{proposition}[{\rm\cite{Heinz,Bhatia-Rosenthal},~\cite[Lemma 2.2]{Kressner-Tobler}}]\label{p:Q+Im}
Let $A\in\mathoo B(X)$, $B\in\mathoo B(Y)$,
and the embeddings $\sigma(A)\subset\{\,\lambda\in\mathbb C:\,\Real\lambda<\rho\,\}$ and $\sigma(B)\subset\{\,\lambda\in\mathbb C:\,\Real\lambda>\rho\,\}$ hold for some $\rho\in\mathbb R$. Then
\begin{equation*}
Q=-\int_0^\infty e^{At}\boxtimes e^{-Bt}\,dt.
\end{equation*}
\end{proposition}

\begin{proof}
We begin with the identity
\begin{equation*}
w(\lambda,\mu)=-\int_0^\infty e^{\lambda t} e^{-\mu t}\,dt.
\end{equation*}
It is valid for $\lambda\in U$ and $\mu\in V$ provided the neighborhoods $U\supset\sigma(A)$ and $V\supset\sigma(B)$ are sufficiently small.
Moreover, we may assume that the integral converges uniformly for $\lambda\in U$ and $\mu\in V$. We substitute this integral into the formula
\begin{equation*}
\bigl(\varphi_1\boxtimes\varphi_2\bigr)w=\frac1{(2\pi i)^2}\int_{\Gamma_1}\int_{\Gamma_2}w(\lambda,\mu)R_{A,\,\lambda}\boxtimes
R_{B,\,\mu}\,d\mu\,d\lambda
\end{equation*}
from Theorem~\ref{t:functional_calculas:infty in rho:R_1,R_2->C} assuming that $\Gamma_1\subset U$ and $\Gamma_2\subset V$.
 By the uniform convergence of the last integral, we may change the order of integration:
\begin{align*}
\bigl(\varphi_1\boxtimes\varphi_2\bigr)(w)&=\frac1{(2\pi i)^2}\int_{\Gamma_1}\int_{\Gamma_2}
\Bigl(-\int_0^\infty e^{\lambda t} e^{-\mu t}\,dt\Bigr)R_{A,\,\lambda}\boxtimes
R_{B,\,\mu}\,d\mu\,d\lambda\\
&=-\int_0^\infty\Bigl(\frac1{(2\pi i)^2}\int_{\Gamma_1}\int_{\Gamma_2} e^{\lambda t}e^{-\mu t}
R_{A,\,\lambda}\boxtimes R_{B,\,\mu}\,d\mu\,d\lambda\Bigr)\,dt\\
&=-\int_0^\infty e^{At}\boxtimes e^{-Bt}\,dt.\qed
\end{align*}
\renewcommand\qed{}
\end{proof}

\begin{proposition}[{\rm\cite[Theorem 9.1]{Bhatia-Rosenthal}}]\label{p:Q+rho}
Let the embeddings $\bar\sigma(R_{1,\,(\cdot)})\subset\{\,\lambda\in\mathbb C:\,|\lambda|<\rho\,\}$ and $\bar\sigma(R_{2,\,(\cdot)})\subset\{\,\lambda\in\overline{\mathbb C}:\,|\lambda|>\rho\,\}$ hold for some $\rho>0$. Then
\begin{equation*}
Q=-\sum_{n=0}^\infty A^{n}\boxtimes R_{2,\,0}^{n+1},
\end{equation*}
where the operator $A\in\mathoo B(X)$ generates $R_{1,\,(\cdot)}$ according to Proposition~{\rm\ref{p:1 in B_R}}.
\end{proposition}
\begin{proof}
We consider the identity
\begin{equation*}
w(\lambda,\mu)=-\sum_{n=0}^\infty\frac{\lambda^n}{\mu^{n+1}}.
\end{equation*}
It is valid for $\lambda\in U$ and $\mu\in V$, where the neighborhoods $U\supset\sigma(A)$ and $V\supset\bar\sigma(R_{2,\,(\cdot)})$ are sufficiently small.
Moreover, we may assume that the series converges uniformly for $\lambda\in U$ and $\mu\in V$. We substitute this series into the formula
 \begin{equation*}
\bigl(\varphi_1\boxtimes\varphi_2\bigr)(w)=\frac1{(2\pi
i)^2}\int_{\Gamma_1}\int_{\Gamma_2}w(\lambda,\mu)R_{1,\,\lambda}\boxtimes
R_{2,\,\mu}\,d\mu\,d\lambda+\frac\delta{2\pi
i}\int_{\Gamma_1}w(\lambda,\infty)R_{1,\,\lambda}\boxtimes\mathbf 1_Y\,d\lambda
\end{equation*}
from Theorem~\ref{t:functional_calculas:infty in rho+sigma:R_1,R_2->C} assuming that $\Gamma_1\subset U$ and $\Gamma_2\subset V$. By the uniform convergence of the series (and by $w(\lambda,\infty)=0$), we have
\begin{align*}
\bigl(\varphi_1\boxtimes\varphi_2\bigr)(w)&=\frac1{(2\pi
i)^2}\int_{\Gamma_1}\int_{\Gamma_2}-\sum_{n=0}^\infty\frac{\lambda^n}{\mu^{n+1}}R_{A,\,\lambda}\boxtimes
R_{2,\,\mu}\,d\mu\,d\lambda\\
&=-\sum_{n=0}^\infty\frac1{(2\pi
i)^2}\int_{\Gamma_1}\int_{\Gamma_2}\frac{\lambda^n}{\mu^{n+1}}R_{A,\,\lambda}\boxtimes
R_{2,\,\mu}\,d\mu\,d\lambda=-\sum_{n=0}^\infty A^{n}\boxtimes R_{2,\,0}^{n+1}.\qed
\end{align*}
\renewcommand\qed{}
\end{proof}

Theorem~\ref{t:Theta via Q} below reduces the calculation of $\bigl[(\varphi_1\boxdot\varphi_2)f\bigr]C$ to the calculation of
$\varphi_1(f)$ and $\varphi_2(f)$ provided $Q(C)$ is known; it is a version of Theorem~\ref{t:Theta via Q:0}.

 \begin{theorem}\label{t:Theta via Q}
Let the extended singular sets $\bar\sigma(R_{1,\,(\cdot)})$ and $\bar\sigma(R_{2,\,(\cdot)})$ of the pseudo-re\-sol\-vents  $R_{1,\,(\cdot)}$ and $R_{2,\,(\cdot)}$ be disjoint,
and $f\in\mathoo O\bigl(\bar\sigma(R_{1,\,(\cdot)})\cup\bar\sigma(R_{2,\,(\cdot)})\bigr)$. Then
\begin{equation*}
\bigl[(\varphi_1\boxdot\varphi_2)f\bigr]C=\bigl[\varphi_1(f)\boxtimes\mathbf1-\mathbf1\boxtimes\varphi_2(f)\bigr]Q(C).
\end{equation*}
In the special case, where $\mathoo B(Y,X)$ is taken as extended tensor product {\rm(}see example~\ref{ex:TP of cat}{\rm(e))},
\begin{equation*}
\bigl[(\varphi_1\boxdot\varphi_2)f\bigr]C=\varphi_1(f)\cdot Q(C)-Q(C)\cdot\varphi_2(f).
\end{equation*}
 \end{theorem}
 \begin{proof}
By Theorem~\ref{t:Theta},
\begin{align*}
\bigl[(\varphi_1\boxdot\varphi_2)f\bigr]C
&=\bigl[(\varphi_1\boxtimes\varphi_2)f^{[1]}\bigr]C
=\bigl(\varphi_1\boxtimes\varphi_2\bigr)\bigl((f\otimes u)w\bigr)C
-\bigl(\varphi_1\boxtimes\varphi_2\bigr)\bigl((u\otimes f)w\bigr)C,
\end{align*}
where $(f\otimes u)(\lambda,\mu)=f(\lambda)$, $(u\otimes f)(\lambda,\mu)=f(\mu)$.
From Theorems \ref{t:functional_calculas:infty in rho:R_1,R_2->C},~\ref{t:functional_calculas:infty in rho+sigma:R_1,R_2->C}, and~\ref{t:functional_calculas:infty in sigma:R_1,R_2->C} it follows that
\begin{align*}
\bigl[(\varphi_1\boxdot\varphi_2)f\bigr]C&=(\varphi_1(f)\boxtimes\mathbf{1})\bigl[\bigl(\varphi_1\boxtimes\varphi_2\bigr)(w)C\bigr]
-(\mathbf{1}\boxtimes\varphi_2(f))\bigl[\bigl(\varphi_1\boxtimes\varphi_2\bigr)(w)C\bigr]\\
&=(\varphi_1(f)\boxtimes\mathbf{1})Q(C)-(\mathbf{1}\boxtimes\varphi_2(f))Q(C).\qed
\end{align*}
\renewcommand\qed{}
 \end{proof}

In Corollary~\ref{c:Theta via Q} below, we present a version of Theorem~\ref{t:Theta via Q}.
It suggests the reverse order of operations, which enables one to apply the transformator $Q$ only once; namely, at first, $\varphi_1(f)$ and $\varphi_2(f)$ are calculated, and then
$Q(\cdot)$ is applied.

 \begin{corollary}\label{c:Theta via Q}
Let the extended singular sets $\bar\sigma(R_{1,\,(\cdot)})$ and $\bar\sigma(R_{2,\,(\cdot)})$ be disjoint. Then
\begin{equation*}
\bigl[(\varphi_1\boxdot\varphi_2)f\bigr]C=Q\bigl(\varphi_1(f)\cdot C-C\cdot\varphi_2(f)\bigr).
\end{equation*}
 \end{corollary}
\begin{proof}
It is sufficient to note that the transformators $\varphi_1(f)\boxtimes\mathbf{1}$, $\mathbf{1}\boxtimes\varphi_2(f)$, and $Q=\bigl(\varphi_1\boxtimes\varphi_2\bigr)(w)$ commute, and then apply Theorem~\ref{t:Theta via Q}.
\end{proof}

Let $A\in\mathoo B(X)$ and $B\in\mathoo B(Y)$.
The equation
\begin{equation}\label{e:Sylvester}
AZ-ZB=C
\end{equation}
for the unknown $Z\in\mathoo B(Y,X)$ with the free term $C\in\mathoo B(Y,X)$
is called the (\emph{continuous{\rm)} Sylvester equation}~\cite{Ikramov1984:eng,Bhatia-Rosenthal,Simoncini1996}.
The Sylvester equation is connected with a number of applications~\cite{Antoulas,Konstantinov-Gu-Mehrmann-Petkov,Ikramov1991:eng,Simoncini2013,Bellman69,Daletskii-Krein:eng,Gantmakher:eng} and is widely discussed in the literature.

 \begin{theorem}\label{t:Sylvester}
Let $A\in\mathoo B(X)$ and $B\in\mathoo B(Y)$. The Sylvester equation~\eqref{e:Sylvester}
has a unique solution $Z\in\mathoo{B}(Y,X)$ for any $C\in\mathoo{B}(Y,X)$ if and only if the spectra of the operators $A$ and $B$ are disjoint.
This solution coincides with the operator $Q(C)$.
 \end{theorem}

 \begin{proof}
By Theorem~\ref{t:functional_calculas:infty in rho:R_1,R_2->C} and Corollary~\ref{c:functional_calculas:infty in rho:R_1,R_2->C}, the transformator $Z\mapsto AZ-ZB$ is equal to $\bigl(\varphi_1\boxtimes\varphi_2\bigr)f$, where $f(\lambda,\mu)=\lambda-\mu$. By Theorem~\ref{t:spectral map:R_1,R_2->C}, its spectrum is equal to $\sigma(A)-\sigma(B)$. Therefore the transformator is invertible if and only if $0\notin\sigma(A)-\sigma(B)$. By Theorem~\ref{t:functional_calculas:infty in rho:R_1,R_2->C}, the inverse transformator is $\bigl(\varphi_1\boxtimes\varphi_2\bigr)w$, where $w(\lambda,\mu)=\frac1{\lambda-\mu}$.
 \end{proof}

The equation
\begin{equation}\label{e:Sylvester:c}
Z-AZB=C
\end{equation}
is called the (\emph{discrete{\rm)} Sylvester equation}~\cite{Ikramov1984:eng} or the \emph{Stein equation}~\cite{Kressner:10.7153/oam-08-23}. Its theory is similar to the theory of equation~\eqref{e:Sylvester}.

 \begin{theorem}\label{t:Sylvester:c}
Let $A\in\mathoo B(X)$ and $B\in\mathoo B(Y)$. The Sylvester equation~\eqref{e:Sylvester:c}
has a unique solution $Z\in\mathoo{B}(Y,X)$ for any $C\in\mathoo{B}(Y,X)$ if and only if the product of the spectra of the operators $A$ and $B$ does not contain~1. This solution coincides with the operator $\bigl[(\varphi_1\boxtimes\varphi_2)(s)\bigr](C)$, where
\begin{equation*}
s(\lambda,\mu)=\frac1{1-\lambda\mu}.
\end{equation*}
 \end{theorem}

 \begin{proof}
By Theorem~\ref{t:functional_calculas:infty in rho:R_1,R_2->C} and Corollary~\ref{c:functional_calculas:infty in rho:R_1,R_2->C}, the transformator $Z\mapsto Z-AZB$ is equal to $\bigl(\varphi_1\boxtimes\varphi_2\bigr)f$, where $f(\lambda,\mu)=1-\lambda\mu$. By Theorem~\ref{t:spectral map:R_1,R_2->C}, its spectrum is equal to $1-\sigma(A)\sigma(B)$. Therefore the transformator is invertible if and only if $1\notin\sigma(A)\sigma(B)$. By Theorem~\ref{t:functional_calculas:infty in rho:R_1,R_2->C}, the inverse transformator is $\bigl(\varphi_1\boxtimes\varphi_2\bigr)s$.
 \end{proof}

\begin{remark}
Let us return to equation~\eqref{e:Sylvester} and discuss the case where $A$ and $B$ are unbounded operators or linear relations.
The natural hypothesis is as follows: if the extended singular sets of the resolvents of $A$ and $B$ are disjoint (and thus $A$ or $B$ is a bounded operator), then equation~\eqref{e:Sylvester} has a unique solution, which is determined by the transformator $Q$. The problem is: How one can interpret equation~\eqref{e:Sylvester}? We discuss some variants.

First, we assume that $A$ and $B$ are linear relations with non-empty resolvent sets, and the extended spectra of $A$ and $B$ are disjoint. We assume that $C\in\mathoo B(Y,X)$
and a solution $Z\in\mathoo{B}(Y,X)$ of equation~\eqref{e:Sylvester} is of interest.

To begin with, we show that without loss of generality one can assume that the inverse operators of $A$ and $B$ are
everywhere  defined bounded operators. Since the extended spectra of $A$ and $B$ are closed and disjoint, there exists $\nu\notin\bar\sigma(A)\cup\bar\sigma(B)$. We rewrite~\eqref{e:Sylvester} in the form
\begin{equation*}
-\nu Z+AZ+\nu Z-ZB=C,
\end{equation*}
and then in the form (with the invertible coefficients $\nu\mathbf1-A$ and $\nu\mathbf1-B$)
\begin{equation*}
-(\nu\mathbf1-A)Z+Z(\nu\mathbf1-B)=C.
\end{equation*}
See~\cite[Proposition A.1.1, p. 281]{Haase2006} or~\cite[Theorem 36]{Kurbatova-PMM09:eng} for a justification of the last equality in the case of linear relations.

We consider the case where $\infty\in\bar\sigma(A)$. Since the relation $A$ is invertible, its range coincides with the whole of $X$, and the image of the zero is zero
(otherwise the left side of equation~\eqref{e:Sylvester} is not an operator). So, $A$ is an operator (not a relation). We call an operator $Z\in\mathoo{B}(Y,X)$, whose range is contained in the domain of the operator $A$ (otherwise the domains of the left and the right sides of equation~\eqref{e:Sylvester} are different), a \emph{solution} of equation~\eqref{e:Sylvester} provided it satisfies the equation. Since $\infty\notin\bar\sigma(B)$, $B$ is a bounded linear operator, see Proposition~\ref{p:1 in B_R}. Multiplying~\eqref{e:Sylvester} by $A^{-1}$ we obtain
\begin{equation}\label{e:Sylvester:modified}
Z-A^{-1}ZB=A^{-1}C.
\end{equation}
By~\cite[Proposition A.3.1]{Haase2006}, $\sigma(A^{-1})=\{\,\frac1\lambda:\,\lambda\in\bar\sigma(A)\,\}$. By Theorem~\ref{t:Sylvester:c}, equation~\eqref{e:Sylvester:modified} has a unique solution  $Z$ for an arbitrary $A^{-1}C$ if
$0\notin\{\,1-\lambda\mu:\,\lambda\in\sigma(A^{-1}),\,\mu\in\sigma(B)\,\}
=\{\,1-\frac\mu\lambda:\,\lambda\in\bar\sigma(A),\,\mu\in\sigma(B)\,\}$,
which is the case, because the extended spectra of $A$ and $B$ are disjoint. We multiply~\eqref{e:Sylvester:modified} on the left by $A$ (taking into account that $AA^{-1}=\mathbf1$).
As a result we arrive at the original equation~\eqref{e:Sylvester}. Therefore $Z$ is a solution of equation~\eqref{e:Sylvester} as well.

We discuss the case where $\infty\in\bar\sigma(B)$. We assume that the domain of the relation $B$ coincides with the whole of $Y$ (otherwise the domains of the left and right sides of equation~\eqref{e:Sylvester} are different). Since the relation $B$ is invertible, its range coincides with the whole of $Y$, the kernel is zero, but the image of the zero $\Image_0B=\{\,x:\,(0,x)\in B\,\}$ may consist not only of zero. We call an operator $Z\in\mathoo{B}(Y,X)$, whose kernel contains $\Image_0B$ (otherwise the left side of equation~\eqref{e:Sylvester} is not an operator), a \emph{solution} of equation~\eqref{e:Sylvester} provided it satisfies the equation. Multiplying~\eqref{e:Sylvester} on the right by $B^{-1}$ we obtain
\begin{equation*}
AZB^{-1}-ZBB^{-1}=CB^{-1}.
\end{equation*}
According to~\cite[Theorem 16]{Kurbatova-PMM09:eng} we rewrite this equation in the form
\begin{equation*}
AZB^{-1}-Z\mathbf1_{Y:\Image_0B}=CB^{-1},
\end{equation*}
where $\mathbf1_{Y:\Image_0B}=\{\,(y_1,y_2)\in Y\times Y:\,(0,y_1-y_2)\in B\,\}$.
Since the kernel of the operator $Z$ contains $\Image_0B$, the last equation can be rewritten as
\begin{equation*}
AZB^{-1}-Z=CB^{-1}.
\end{equation*}
By Theorem~\ref{t:Sylvester:c}, this equation has a unique solution $Z$ for an arbitrary $CB^{-1}$ if $0\notin\{\,1-\lambda\mu:\,\lambda\in\sigma(A),\,\mu\in\sigma(B^{-1})\,\}
=\{\,1-\frac\lambda\mu:\,\lambda\in\sigma(A),\,\mu\in\bar\sigma(B)\,\}$, which is the case, because the extended spectra of $A$ and $B$ are disjoint. Multiplying the last equation on the right by $B$ we obtain
\begin{equation*}
AZB^{-1}B-ZB=CB^{-1}B,
\end{equation*}
or (according to~\cite[Theorem 16]{Kurbatova-PMM09:eng} and $\Ker B={0}$)
\begin{equation*}
AZ-ZB=C.
\end{equation*}
So, $Z$ is a solution of original equation~\eqref{e:Sylvester} as well.

We consider another case: let $B$ be an invertible unbounded operator with the dense domain $\Dom B$ (in particular, $\infty\in\bar\sigma(B)$). We call an operator $Z\in\mathoo{B}(Y,X)$ a \emph{solution} if
\begin{equation*}
AZy-ZBy=Cy
\end{equation*}
for all $y\in\Dom B$. Let us look for a solution of equation~\eqref{e:Sylvester} in the form $Z=VB^{-1}$, where $V\in\mathoo{B}(Y,X)$ is a new unknown operator. Substituting $Z=VB^{-1}$ into equation~\eqref{e:Sylvester} we obtain
\begin{equation}\label{e:AVB-1-VB-1B=C}
AVB^{-1}-VB^{-1}B=C,
\end{equation}
or
\begin{equation*}
AVB^{-1}-V\mathbf1_{\Dom B}=C,
\end{equation*}
where $\mathbf1_{\Dom B}$ is the identity operator with the domain $\Dom B$. By our definition of a solution, the last equation is  equivalent to the equation
\begin{equation*}
AVB^{-1}-V=C.
\end{equation*}
Obviously, it has a unique solution $V$. Returning to equivalent equation~\eqref{e:AVB-1-VB-1B=C} we see that the operator $Z=VB^{-1}$ is a solution of the original equation.
\end{remark}

Theorem~\ref{t:Sylvester} for matrices was first proved in~\cite{Sylvester}.
An independent proof of its sufficient part for the case of operators was obtained in~\cite{Krein64,Daleckii53:eng,Rosenblum56:Duke}.
For a Hilbert space, a necessary and sufficient condition for the solvability of the Sylvester equation was first obtained in~\cite{Davis_Ch-Rosenthal}, see also~\cite[c.~54]{Gohberg-Goldberg-Kaashoek}. An analogue of Theorem~\ref{t:Sylvester:c} for matrices is proved, for example, in~\cite{Kressner:10.7153/oam-08-23}.

The representation for the solution of the Sylvester equation in the form of contour integral~\eqref{e:Q via contour} was first published in~\cite{Rosenblum56:Duke}, see also Example~\ref{ex:tensor sum}.
Estimates of the solution of the Sylvester equation are given in~\cite{Gil-11:EJLA,Gil2015resolvents,Gil2015Hindawi}.

The Sylvester equation~\eqref{e:Sylvester} with unbounded operator coefficients $A$ and $B$ is considered in~\cite{Albeverio-Motovilov:eng,Albeverio-Makarov-Motovilov,
Kostrykin-Makarov-Motovilov03,Phong_Vu_Quoc,Freeman,Nguyen_Thanh_Lan,Shaw-Lin}.

A generalization of the transformator $Q$ ($Q$ corresponds to the function $w(\lambda,\mu)=\frac1{\lambda-\mu}$) is the inverse of the transformator
$v_{n+1}^{[1]}(A,B)$, where $v_n^{[1]}(\lambda,\mu)=\lambda^{n-1}+\lambda^{n-2}\mu+\dots+\mu^{n-1}$. It is discussed in~\cite{Bhatia-Uchiyama,Gil2015Hindawi,Furuta}.

\section{The differential of the functional calculus}\label{s:df}
Let $X$ be a Banach space.
The \emph{{\rm(}Fr\'echet{\rm)} differential} of a nonlinear transformator
$f\:D(f)\subseteq\mathoo B(X)\to\mathoo B(X)$ at a point $A\in\mathoo B(X)$ is defined to be a linear transformator
$df(\cdot,A)\:\mathoo B(X)\to\mathoo B(X)$ depending on the parameter $A$ that
possesses the property
\begin{equation}\label{e:Delta f}
f(A+\Delta A)=f(A)+df(\Delta A,A)+o(\Vert\Delta A\Vert).
\end{equation}
We assume that a neighborhood of $A$ is contained in the domain $D(f)$ of the transformator~$f$. We recall standard properties of the differential.

 \begin{proposition}[{\rm\cite[\S~2.2.2]{Alexeev-Tikhomirov-Fomin:eng}, \cite[8.2.1]{Dieudonne:eng}%
}]\label{p:properties of d}
Let a transformator $g\:\mathoo B(X)\to\mathoo B(X)$ be differentiable at a point $A\in\mathoo B(X)$ and
a transformator $f\:\mathoo B(X)\to\mathoo B(X)$ be differentiable at the point $g(A)\in\mathoo B(X)$. Then
the composition $f\circ g$ is differentiable at the point $A$, and
\begin{equation*}
d(f\circ g)(\cdot,A)=df\bigl[dg(\cdot,A),g(A)\bigr].
\end{equation*}
 \end{proposition}

 \begin{corollary}[{\rm\cite[\S~2.3.4]{Alexeev-Tikhomirov-Fomin:eng}, \cite[8.2.3]{Dieudonne:eng}}]\label{c:d of inverse}
Let a transformator $f\:\mathoo B(X)\to\mathoo B(X)$ be continuously differentiable {\rm(}i.~e. $df(\cdot,A)$ depends on $A$ continuously in norm{\rm)} in a neighborhood of a point $A\in\mathoo B(X)$ and let the transformator $df(\cdot,A)$ be invertible. Then the inverse transformator of $f$ is defined and differentiable in a neighborhood of the point $B=f(A)$, and the differential of the inverse transformator is equal to the inverse of the original differential:
\begin{equation*}
 df^{-1}(\cdot,B)=\bigl[df(\cdot,A)\bigr]^{-1}.
\end{equation*}
 \end{corollary}

Let $A\in\mathoo{B}(X)$ and $f\in\mathoo O\bigl(\sigma(A)\bigr)$.
For the transformator
\begin{equation}\label{e:def of functional calculus:A}
A\mapsto f(A)=\frac1{2\pi i}\int_\Gamma f(\lambda)(\lambda\mathbf1-A)^{-1}\,d\lambda
\end{equation}
definition~\eqref{e:Delta f} of a differential looks as follows:
\begin{equation*}
\frac1{2\pi i}\int_\Gamma f(\lambda)\bigl(\lambda\mathbf1-(A+\Delta
A)\bigr)^{-1}\,d\lambda=\frac1{2\pi i}\int_\Gamma
f(\lambda)(\lambda\mathbf1-A)^{-1}\,d\lambda+df(\Delta A,A)+o(\Vert\Delta A\Vert).
\end{equation*}

We note that
\begin{equation*}
\bigl(\lambda\mathbf1-(A+\Delta A)\bigr)^{-1}=\bigl((\lambda\mathbf1-A)-\Delta A\bigr)^{-1}
=R_\lambda(\mathbf1-\Delta A\cdot R_\lambda)^{-1}=(\mathbf1-R_\lambda\cdot\Delta A)^{-1}R_\lambda.
\end{equation*}
Based on this formula, we adopt the following definition.

Let $R_{(\cdot)}$ be a pseudo-resolvent in the algebra $\mathoo B(X)$.
We call \emph{the perturbation of $R_{(\cdot)}$ by an operator} $\Delta A\in\mathoo B(X)$ the function
\begin{equation}\label{e:def of Delta R}
T_\lambda=R_\lambda(\mathbf1-\Delta A\cdot R_\lambda)^{-1}=(\mathbf1-R_\lambda\cdot\Delta A)^{-1}R_\lambda.
\end{equation}

\begin{remark}
We note shortly an additional reasoning in favor of definition~\eqref{e:def of Delta R}. Let $R_{(\cdot)}$ be the resolvent of a linear relation $A$, i.~e. $R_\lambda=(\lambda\mathbf1-A)^{-1}$. We show that
\begin{equation*}
(\mathbf1-R_\lambda\cdot\Delta A)^{-1}R_\lambda=(\lambda\mathbf1-A-\Delta A)^{-1}.
\end{equation*}
Obviously (for details, see~\cite[Proposition A.1.1]{Haase2006} or~\cite[Proposition 12]{Kurbatova-PMM09:eng}),
$(\mathbf1-R_\lambda\cdot\Delta A)^{-1}R_\lambda=\bigl[R_\lambda^{-1}(\mathbf1-R_\lambda\cdot\Delta A)\bigr]^{-1}=\bigl[(\lambda\mathbf1-A)(\mathbf1-R_\lambda\cdot\Delta A)\bigr]^{-1}$. Further, since the image of the operator $R_\lambda\cdot\Delta A$ is contained in the image of $R_\lambda$, which is equal to the domain of the relation $A$, by virtue of~\cite[Theorem 36(a)]{Kurbatova-PMM09:eng}, we can develop the internal parentheses: $\bigl[\lambda\mathbf1-A-(\lambda\mathbf1-A)R_\lambda\cdot\Delta A)\bigr]^{-1}=\bigl[\lambda\mathbf1-A-(\lambda\mathbf1-A)(\lambda\mathbf1-A)^{-1}\cdot\Delta A\bigr]^{-1}$. We note that $(\lambda\mathbf1-A)(\lambda\mathbf1-A)^{-1}$ is equal to the relation $\mathbf1_{X:\Image_0A}=\{\,(x_1,x_2)\in X\times X:\,(0,x_1-x_2)\in A\,\}$. Obviously, $(\lambda\mathbf1-A-\mathbf1_{X:\Image_0A}\cdot\Delta A)^{-1}=(\lambda\mathbf1-A-\Delta A)^{-1}$.
\end{remark}

\begin{proposition}
For any perturbation $\Delta A\in\mathoo B(X)$ the function
\begin{equation*}
T_\lambda=R_\lambda(\mathbf1-\Delta A\cdot R_\lambda)^{-1}=(\mathbf1-R_\lambda\cdot\Delta A)^{-1}R_\lambda
\end{equation*}
is  a pseudo-resolvent.
\end{proposition}
\begin{proof}
We verify the Hilbert identity for all $\lambda$ and $\mu$ such that $T_\lambda$ and $T_\mu$ are defined. We have
\begin{align*}
T_\lambda-T_\mu+(\lambda-\mu)T_\lambda T_\mu
&=(\mathbf1-R_\lambda\cdot\Delta A)^{-1}R_\lambda-R_\mu(\mathbf1-\Delta A\cdot R_\mu)^{-1}\\
&+(\lambda-\mu)(\mathbf1-R_\lambda\cdot\Delta A)^{-1}R_\lambda R_\mu(\mathbf1-\Delta A\cdot R_\mu)^{-1}\\
&=(\mathbf1-R_\lambda\cdot\Delta A)^{-1}\bigl[R_\lambda(\mathbf1-\Delta A\cdot R_\mu)-(\mathbf1-R_\lambda\cdot\Delta A)R_\mu\\
&+(\lambda-\mu)R_\lambda R_\mu\bigr]\\
&=(\mathbf1-R_\lambda\cdot\Delta A)^{-1}\bigl[R_\lambda-R_\mu+(\lambda-\mu)R_\lambda R_\mu\bigr]=0.\qed
\end{align*}
\renewcommand\qed{}
\end{proof}

We define the \emph{differential $df(\cdot,R_{(\cdot)})$ of the mapping} (which is a generalization of~\eqref{e:def of functional calculus:A})
\begin{equation}\label{e:def of functional calculus:R}
R_{(\cdot)}\mapsto f(R_{(\cdot)})=\frac1{2\pi i}\int_\Gamma f(\lambda)R_\lambda\,d\lambda+\delta f(\infty)\mathbf1,
\end{equation}
by means of the formula
\begin{multline*}
\frac1{2\pi i}\int_\Gamma f(\lambda)R_\lambda(\mathbf1-\Delta A\cdot R_\lambda)^{-1}\,d\lambda+\delta f(\infty)\mathbf1\\
=\frac1{2\pi i}\int_\Gamma
f(\lambda)R_\lambda\,d\lambda+\delta f(\infty)\mathbf1+df(\Delta A,R_{(\cdot)})+o(\Vert\Delta A\Vert).
\end{multline*}

 \begin{theorem}\label{t:differential}
Let $R_{(\cdot)}$ be a pseudo-resolvent in the algebra $\mathoo B(X)$, and $f\in\mathoo O\bigl(\bar\sigma(R_{(\cdot)})\bigr)$. Then
the differential of mapping~\eqref{e:def of functional calculus:R}
admits the representation
\begin{equation*}
df(\Delta A,R_{(\cdot)})=\frac1{2\pi i}\int_\Gamma f(\lambda)R_\lambda\,\Delta A\,R_\lambda\,d\lambda.
\end{equation*}
In other words,
\begin{equation*}
df(\cdot,R_{(\cdot)})=\bigl(\varphi\boxdot\varphi\bigr)(f),
\end{equation*}
where $\varphi$ is the functional calculus generated by the pseudo-resolvent $R_{(\cdot)}$.
 \end{theorem}

 \begin{proof}
We assume that
\begin{equation*}
\Vert\Delta A\Vert\cdot\Vert R_\lambda\Vert<1.
\end{equation*}
By Theorem~\ref{t:Neumann}, we have
\begin{equation*}
\bigl\Vert\bigl[ R_\lambda(\mathbf1-\Delta A\cdot R_\lambda)^{-1}-R_\lambda\bigr]-R_\lambda\Delta A\cdot R_\lambda\bigr\Vert
\le\frac{\Vert R_\lambda\Vert^3\cdot\Vert\Delta A\Vert^2}{1-\Vert R_\lambda\Vert\cdot\Vert\Delta A\Vert}.
\end{equation*}
Therefore,
\begin{gather*}
\Biggl\Vert\frac1{2\pi i}\int_\Gamma f(\lambda)R_\lambda(\mathbf1-\Delta A\cdot R_\lambda)^{-1}\,d\lambda-\frac1{2\pi i}\int_\Gamma
f(\lambda)R_\lambda\,d\lambda-\frac1{2\pi i}\int_\Gamma
f(\lambda)R_\lambda\,\Delta A\,R_\lambda\,d\lambda\Biggr\Vert\\
=\Bigl\Vert\frac1{2\pi i}\int_\Gamma f(\lambda)\bigl[R_\lambda(\mathbf1-\Delta A\cdot R_\lambda)^{-1}-R_\lambda
-R_\lambda\,\Delta A\,R_\lambda\bigr]\,d\lambda\Bigr\Vert\\
\le\frac1{2\pi}\int_\Gamma
|f(\lambda)|\,\frac{\Vert R_\lambda\Vert^3\cdot\Vert\Delta A\Vert^2}{1-\Vert R_\lambda\Vert\cdot\Vert\Delta A\Vert}\,|d\lambda|=o(\Vert\Delta A\Vert).\qed
\end{gather*}
\renewcommand\qed{}
 \end{proof}

\begin{proposition}[\rm{\cite[Theorem 2.1]{Bhatia-Sinha99}}]\label{p:Theta via Q:1}
Let $A\in\mathoo{B}(X)$ and $f\in\mathoo O\bigl(\sigma(A)\bigr)$. Then
\begin{equation*}
df(A\,\Delta A-\Delta A\,A,A)=\varphi_A(f)\Delta A-\Delta A\varphi_A(f),
\end{equation*}
where the functional calculus $\varphi_A$ is generated by the operator $A$.
\end{proposition}
\begin{proof}
The proof follows from Theorems~\ref{t:Theta via Q:0} and~\ref{t:differential}.
\end{proof}

\begin{proposition}[{\cite[Theorem 3.3]{Higham08}}]\label{p:mult:df} Let $g,h\in\mathoo O\bigl(\bar\sigma(R_{(\cdot)})\bigr)$. Then
\begin{equation*}
d(gh)(\Delta A,R_{(\cdot)})=dg(\Delta A,R_{(\cdot)})\,h(R_{(\cdot)})+g(R_{(\cdot)})\,dh(\Delta A,R_{(\cdot)}),
\end{equation*}
where $(gh)(\lambda)=g(\lambda)h(\lambda)$.
\end{proposition}
\begin{proof}
The proof follows from Corollary~\ref{c:mult}.
\end{proof}

\begin{corollary}[{\rm\cite[Theorem 10.36]{Rudin:eng}, see also~\cite{Rinehart}}]\label{c:Rinehart}
Let $f\in\mathoo O\bigl(\bar\sigma(R_{(\cdot)})\bigr)$.
We assume that $\Delta A$ commutes with at least one value $R_\mu$ of the pseudo-resolvent $R_{(\cdot)}$.
Then
\begin{equation*}
df(\Delta A,R_{(\cdot)})=\varphi(f')\,\Delta A=\Delta A\,\varphi(f').
\end{equation*}
\end{corollary}
\begin{proof}
The proof follows from~\ref{p:phi(f')}.
\end{proof}

 \begin{theorem}\label{t:differential2}
Let $R_{(\cdot)}$ be a pseudo-resolvent in the algebra $\mathoo B(X)$ and $f\in\mathoo O\bigl(\bar\sigma(R_{(\cdot)})\bigr)$. Then
 the differential of  mapping~\eqref{e:def of functional calculus:R} admits the representation
\begin{equation}\label{e:df via boxtimes}
df(\Delta A,R_{(\cdot)})=\frac1{(2\pi
i)^2}\int_{\Gamma_1}\int_{\Gamma_2}f^{[1]}(\lambda,\mu)R_\lambda\,\Delta A\,R_\mu\,d\mu\,d\lambda,
\end{equation}
where the divided difference $f^{[1]}$ is defined by formula~\eqref{e:f[1]}, the contours $\Gamma_1$ and $\Gamma_2$ are  oriented
envelopes of the extended singular set $\bar\sigma(R_{(\cdot)})$ with respect to the complement $\overline{\mathbb C}\setminus U$, and  $U$ is the domain of the function $f$.

The spectrum of the transformator $df(\cdot,R_{(\cdot)})\:\mathoo B(X)\to\mathoo B(X)$ is given by the formula
\begin{equation}\label{e:sigma(df)}
\sigma\bigl[df(\cdot,R_{(\cdot)})\bigr]
=\bigl\{\,f^{[1]}(\lambda,\mu):\,\lambda,\mu\in\bar\sigma(R_{(\cdot)})\,\bigr\}.
\end{equation}
\end{theorem}

 \begin{proof}
Theorem~\ref{t:differential} shows that $df(\Delta A,R_{(\cdot)})$ is the operator
$\bigl[\bigl(\varphi\boxdot\varphi\bigr)(f)\bigr]\Delta A$. By Theorem~\ref{t:Theta},
\begin{equation*}
\bigl[\bigl(\varphi\boxdot\varphi\bigr)(f)\bigr]\Delta A=\bigl[\bigl(\varphi\boxtimes\varphi\bigr)f^{[1]}\bigr]\Delta A.
\end{equation*}
It remains to apply Theorem~\ref{t:functional_calculas:infty in rho:R_1,R_2->C}.

Formula~\eqref{e:sigma(df)} follows from Theorem~\ref{t:spectral map:R_1,R_2->C}.
 \end{proof}

 \begin{example}\label{ex:df(Delta A,A)}
Let $A\in\mathoo B(X)$.
The following corollaries are consequences of Example~\ref{ex:f[1]} and Theorem~\ref{t:differential2}.

The differential of the transformator $v_2(A)=A^2$ is given by the formula
\begin{equation}\label{e:dv_2}
dv_2(\Delta A,A)=A\cdot\Delta A+\Delta A\cdot A,
\end{equation}
and its spectrum at a point $A\in\mathoo B(X)$ is equal to
\begin{equation*}
\sigma\bigl[dv_2(\cdot,A)\bigr]=\bigl\{\,\lambda+\mu:\,\lambda,\mu\in\sigma(A)\,\bigr\}.
\end{equation*}

The differential of the mapping $r_1(R_{(\cdot)})=R_{\lambda_0}$ is given by the formula
\begin{equation*}
dr_1(\Delta A,R_{(\cdot)})=R_{\lambda_0}\cdot\Delta A\cdot R_{\lambda_0},
\end{equation*}
and its spectrum at a point $R_{(\cdot)}$ is equal to
\begin{equation*}
\sigma\bigl[dr_1(\cdot,R_{(\cdot)})\bigr]=\Bigl\{\,\dfrac1{(\lambda_0-\lambda)(\lambda_0-\mu)}:\,\lambda,\mu\in\bar\sigma(R_{(\cdot)})\,\Bigr\}.
\end{equation*}

The differential of the transformator $\exp_t(A)=e^{At}$ is given by the formula
\begin{equation*}
d\exp_t(\Delta A,A)=\sum_{n=0}^\infty\frac{t^n}{(n+1)!}\sum_{i=0}^nA^{n-i}\Delta A\,A^i,
\end{equation*}
and its spectrum at a point $A\in\mathoo B(X)$ is equal to
\begin{equation*}
\sigma\bigl[d\exp_t(\cdot,A)\bigr]
=\bigl\{\,\exp_t^{[1]}(\lambda,\mu):\,\lambda,\mu\in\sigma(A)\,\bigr\}.
\end{equation*}

The differential of the transformator $\exp_t^{(1)}(A)=A\,e^{At}$ is given by the formula
\begin{equation*}
d\exp_t^{(1)}(\Delta A,A)=\sum_{n=0}^\infty\frac{t^n}{n!}\sum_{i=0}^nA^{n-i}\Delta A\,A^i.
\end{equation*}
and its spectrum at a point $A\in\mathoo B(X)$ is equal to
\begin{equation*}
\sigma\bigl[d\exp_t^{(1)}(\cdot,A)\bigr]
=\bigl\{\,\exp_t^{(1)\,[1]}(\lambda,\mu):\,\lambda,\mu\in\sigma(A)\,\bigr\}.
\end{equation*}
 \end{example}

We note special formulae for the differentials of the transformators $\exp_t(A)=e^{At}$ and $\exp_t^{(1)}(A)=Ae^{At}$.
\begin{proposition}\label{p:d exp_t}
The differentials of the transformators $\exp_t(A)=e^{At}$ and $\exp_t(A)^{(1)}=Ae^{At}$ at a point $A\in\mathoo B(X)$ can be calculated by means of the formulae
\begin{align}
d\exp_t(\Delta A,A)&=\int_0^te^{(t-s) A}\Delta Ae^{s A}\,ds=\int_0^t\exp_{t-s}(A)\Delta A\exp_s(A)\,ds,\label{e:d exp_t:Karplus-Schwinger}\\
d\exp_t^{(1)}(\Delta A,A)&=\exp_t(A)\Delta A+\int_0^t\exp_{t-s}(A)\Delta A\,A\exp_s(A)\,ds.\label{e:d exp_t(1)}
\end{align}
\end{proposition}

\begin{proof}
For $\lambda\neq\mu$ we have
\begin{align*}
\frac{e^{\lambda t}-e^{\mu t}}{\lambda-\mu}&=\frac1{\lambda-\mu}e^{\lambda s}e^{\mu(t-s)}\Big|_{s=0}^{s=t}=
\frac1{\lambda-\mu}\int_0^t\frac d{ds}\Bigl[e^{\lambda s}e^{\mu(t-s)}\Bigr]\,ds=
\int_0^t e^{\lambda s}e^{\mu(t-s)}\,ds.
\end{align*}
By the continuity, the same representation of $\exp_t^{[1]}(\lambda,\mu)$ holds for all finite $\lambda$ and $\mu$. Hence, from Theorems~\ref{t:differential2},~\ref{t:functional_calculas:infty in rho:R_1,R_2->C}, and~\ref{t:spectral map:R_1,R_2->C} it follows~\eqref{e:d exp_t:Karplus-Schwinger}. Formula~\eqref{e:d exp_t(1)} follows from~\eqref{e:d exp_t(1):0}.
\end{proof}

The differentials of inverse functions are defined by the inverse transformators, see Corollary~\ref{c:d of inverse}.
To calculate them and their spectra, one can use the following theorem.

 \begin{theorem}\label{t:Frechet of inverse}
Let $R_{(\cdot)}$ be a resolvent of an operator $A\in\mathoo B(X)$, $f\in\mathoo{O}\bigl(\sigma(A)\bigr)$, and $0\notin\bigcup\bigl\{\,f^{[1]}(\lambda,\mu):\,\lambda,\mu\in\sigma(A)\,\bigr\}$. Then
the differential of the transformator $B\mapsto f^{-1}(B)$ at the point $B=f(A)$ is given by the formula
\begin{equation}\label{e:df{-1}}
df^{-1}(\Delta B,B)=\frac1{(2\pi
i)^2}\int_{\Gamma_1}\int_{\Gamma_2}\frac{1}{f^{[1]}(\lambda,\mu)}R_{\lambda}\,\Delta B\,R_{\mu}\,d\mu\,d\lambda,
\end{equation}
where the contours $\Gamma_1$ and $\Gamma_2$ are oriented envelopes of the spectrum $\sigma(A)$ with respect to the point $\infty$ and the complement of the domain of the function $f$.

The spectrum of the transformator $df^{-1}(\cdot,B)\:\mathoo B(X)\to\mathoo B(X)$ is given by the formula
\begin{equation}\label{e:df{-1}:sp}
\sigma\bigl[df^{-1}(\cdot,B)\bigr]
=\bigcup\Bigl\{\,\frac{1}{f^{[1]}(\lambda,\mu)}:\,\lambda,\mu\in\sigma(A)\,\Bigr\}.
\end{equation}
 \end{theorem}

\begin{proof}
By Theorem~\ref{t:differential2},
\begin{align*}
df(\Delta A,A)&=\frac1{(2\pi
i)^2}\int_{\Gamma_1}\int_{\Gamma_2}f^{[1]}(\lambda,\mu)R_\lambda\,\Delta A\,R_\mu\,d\mu\,d\lambda,\\
\sigma\bigl[df(\cdot,A)\bigr]
&=\bigl\{\,f^{[1]}(\lambda,\mu):\,\lambda,\mu\in\sigma(A)\,\bigr\}.
\end{align*}
Since $0\notin\bigcup\bigl\{\,f^{[1]}(\lambda,\mu):\,\lambda,\mu\in\sigma(A)\,\bigr\}$, the transformator
$df(\cdot,A)$ is invertible. By Corollary~\ref{c:d of inverse}, the inverse transformator is the differential of the mapping $B\mapsto f^{-1}(B)$, which is the inverse  mapping of $A\mapsto f(A)$.
By Theorem~\ref{t:functional_calculas:infty in rho:R_1,R_2->C}, the inverse transformator is given by formula~\eqref{e:df{-1}}. By Theorem~\ref{t:spectral map:R_1,R_2->C}, its spectrum is given by formula~\eqref{e:df{-1}:sp}.
\end{proof}

\begin{example}
Let the spectrum of an operator $B\in\mathoo B(X)$ be contained in $\mathbb C\setminus(-\infty,0]$.
The following corollaries are consequences of Example~\ref{ex:df(Delta A,A)} and Theorem~\ref{t:Frechet of inverse}.

From~\eqref{e:dv_2} it is clear that the differential $dv_{1/2}(\Delta B,B)$ of the transformator $v_{1/2}(B)=\sqrt B$ satisfies the Sylvester equation
\begin{equation*}
\sqrt{B}\cdot dv_{1/2}(\Delta B,B)+dv_{1/2}(\Delta B,B)\sqrt{B}=\Delta B.
\end{equation*}
Therefore,
\begin{equation*}
dv_{1/2}(\Delta B,B)=Q_{\sqrt B,\,-\sqrt B}(\Delta B),
\end{equation*}
where the transformator $Q_{\sqrt B,\,-\sqrt B}$ is constructed by means of the functional calculus generated by the operators $\sqrt{B}$ and $-\sqrt{B}$.
The spectrum of the differential of the transformator $v_{1/2}(B)=\sqrt B$ at the point $B$ is equal to
\begin{equation*}
\sigma\bigl[dv_{1/2}(\cdot,B)\bigr]
=\Bigl\{\,\frac1{\lambda+\mu}:\,\lambda,\mu\in\sqrt{\sigma(B)}\,\Bigr\}.
\end{equation*}

The spectrum of the differential of the transformator $f:\,B\mapsto\ln B$, where $\ln$ denotes the principal value, at the point $B$ is equal to
\begin{equation*}
\sigma\bigl[df(\cdot,B)\bigr]
=\Bigl\{\,\frac1{\exp_1^{[1]}(\lambda,\mu)}:\,\lambda,\mu\in\ln\bigl(\sigma(B)\bigr)\,\Bigr\}.
\end{equation*}
\end{example}

\begin{remark}
The equation
\begin{equation}\label{e:a-Riccati}
AZ+ZB+ZCZ+D=0
\end{equation}
with $A,B,C,D\in\mathoo B(X)$ and unknown $Z\in\mathoo B(X)$ is called~\cite{Ikramov1984:eng,Kostrykin-Makarov-Motovilov03,Kostrykin-Makarov-Motovilov05}
 the \emph{Riccati equation}.
It arises in control theory~\cite{Pervozvanskii:eng,Antoulas}.
The differential $dZ=dZ(\Delta A,\Delta B,\Delta C,\Delta D;Z)$ of the solution $Z$ of Riccati equation~\eqref{e:a-Riccati} satisfies~\cite[p. 135]{Ikramov1984:eng} the continuous Sylvester equation
$$
(A+ZC)dZ+dZ(B+CZ)=-\Delta D-\Delta A\,Z-Z\,\Delta B-Z\,\Delta C\,Z.
$$
So,
\begin{equation*}
dZ(\Delta A,\Delta B,\Delta C,\Delta D;Z)=Q_{A+ZC,\,-B-CZ}(-\Delta D-\Delta A\,Z-Z\,\Delta B-Z\,\Delta C\,Z).
\end{equation*}
%
\end{remark}

For pseudo-resolvents generated by bounded operators, Theorem~\ref{t:differential} is proved in~\cite[Theorem 10.38]{Rudin:eng}, see also \cite[formula (2.3)]{Bhatia-Sinha99},~\cite{Dal57:eng} and~\cite{Stickel}.
Representation~\eqref{e:df via boxtimes} for matrices is given in~\cite[Theorem 5.1]{Kressner:10.7153/oam-08-23}.
Formula~\eqref{e:sigma(df)} for matrices is proved in~\cite[Theorem 3.9]{Higham08}, its weaker version was previously obtained in~\cite[Lemma~2.1]{Kenney-Laub89}.
Formula~\eqref{e:d exp_t:Karplus-Schwinger} was first obtained in~\cite[formula (1.8)]{Karplus-Schwinger}, see also~\cite[ch.~10, \S~14]{Bellman69},~\cite[formula (10.15)]{Higham08}, \cite[example 2]{Kenney-Laub89},~\cite{Mathias92},~\cite{Najfeld-Havel},~\cite{van_Loan77}.
Differentials connected with some specific functions $f$ are investigated in~\cite{Kenney-Laub89,Kenney-Laub98,van_Loan77,Al-Mohy-Higham-Relton,Higham-Lin,Zhu-Xue-Gao,Dieci-Morini-Papini,Higham08}; a special attention is paid to estimates of their norms which helps to find the condition number of the transformator $A\mapsto f(A)$.
Properties of differentials of higher orders are investigated in~\cite{Bhatia-Sinha99}.

\end{document}